\documentclass[twoside]{report}
\usepackage[a4paper]{geometry}
\usepackage[latin1]{inputenc} 

\usepackage{tensor}
\usepackage{color}
\newcommand{\blue}{\color{blue}}

\usepackage{graphicx}
\usepackage{ upgreek }
\usepackage{graphicx, amsmath,  amssymb}
\usepackage{accents}

\newcommand*{\dt}[1]{%
  \accentset{\mbox{\large\bfseries .}}{#1}}

\newcommand*{\dtv}[1]{%
  \accentset{\mbox{\large\bfseries {\scriptsize$v$}}}{#1}}

\usepackage{amsthm}

\usepackage{algorithm}
\usepackage{algorithmicx}
\usepackage[noend]{algpseudocode}

\newtheorem{pr}{Proposition}[chapter]
\newtheorem{lm}{Lemma}[chapter]
\newtheorem{al}{Algorithm}[chapter]
\newtheorem{definition}{Definition}[chapter]

\newtheorem{theorem}{Theorem}[chapter]

\newcommand{\uideg}{\textnormal{deg}_w}

\newcommand{\giat}{\widehat{g}}

\newcommand{\E}{\mathcal{E}}

\begin{document}

\newcommand{\revtex}{REV\TeX\ }
\newcommand{\classoption}[1]{\texttt{#1}}
\newcommand{\macro}[1]{\texttt{\textbackslash#1}}
\newcommand{\m}[1]{\macro{#1}}
\newcommand{\env}[1]{\texttt{#1}}

\newcommand{\C}{\mathcal{C}}

\newcommand{\SUIO}{$\mathcal{SUIO}$}
\newcommand{\B}{$\mathcal{B}$}
\newcommand{\VS}{$\mathcal{V}$}
\newcommand{\IMU}{$\mathcal{I}$}
\newcommand{\tIMU}{\textnormal{\IMU}}
\newcommand{\tVS}{\textnormal{\VS}}
\newcommand{\tih}{\widetilde{h}}
\newcommand{\tobs}{\widetilde{\OBS}}
\newcommand{\NO}{s}
\newcommand{\ND}{r}

\newsavebox{\mybox}

\newcommand{\chrono}{Ecco cosa succede quando si passa dal chonospace allo spazio normale:~}

\newcommand{\M}{\mathcal{M}}
\newcommand{\Li}{\mathcal{L}}
\newcommand{\RM}{\mathcal{RM}}
\newcommand{\Hf}{\mathcal{H}}

\newcommand{\oplusn}{+}
\newcommand{\bigoplusn}{\sum}

\setlength{\textheight}{9.5in}

\newcommand{\Obs}{$\mathcal{O}$}
\newcommand{\OBS}{\mathcal{O}}




\title{Nonlinear Unknown Input Observability and Unknown Input Reconstruction: The General Analytical Solution}

\author{Agostino Martinelli
\thanks{A. Martinelli is with INRIA Rhone Alpes,
Montbonnot, France e-mail: {\tt agostino.martinelli@inria.fr}} }

\maketitle



\tableofcontents

\begin{abstract}
Observability is a fundamental structural property of any dynamic system
and describes the possibility of reconstructing the state that characterizes the system from observing its inputs and outputs.
Despite the huge effort made to study this property and to introduce analytical criteria able to check whether a dynamic system satisfies this property or not, there is no general analytical criterion to automatically check the state observability when the dynamics are also driven by unknown inputs. Here, we introduce the general analytical solution of this fundamental problem, often called the unknown input observability problem. 
This paper provides the general analytical solution of this problem, namely, it provides the systematic procedure, based on automatic computation (differentiation and matrix rank determination), that allows us to automatically check the state observability even in the presence of unknown inputs (Algorithm \ref{AlgoFull} in Chapter \ref{ChapterSolutionNonCanonic}).
A first solution of this problem was presented in the second part of the book: {\it Observability: A New Theory Based on the Group of Invariance} \cite{SIAMbook}. The solution presented by this paper 
completes the previous solution in \cite{SIAMbook}. In particular, the new solution exhaustively accounts for the systems that do not belong to the category of the systems that are
{\it canonic with respect to their unknown inputs}.
The new solution is also provided 
in the form of a new algorithm (Algorithm \ref{AlgoFull}). A further novelty with respect to the algorithm provided in \cite{SIAMbook} consists of a new convergence criterion that holds in all the cases (the convergence criterion of the algorithm provided in \cite{SIAMbook} can fail in some cases).
The analytical derivations largely exploit several new concepts and analytical results introduced in \cite{SIAMbook}.
We illustrate the implementation of the new algorithm by studying the observability properties of a nonlinear system in the framework of visual-inertial sensor fusion, whose dynamics are driven by two unknown inputs and one known input. In particular, for this system, we follow step by step the algorithm introduced by this paper, which solves the unknown input observability problem in the most general case.

\end{abstract}

\chapter{Introduction}\label{ChapterIntroduction}

\noindent Observability refers to the state that characterizes a dynamic system (e.g., if the system is an aerial drone, its state, under suitable conditions, can be its position and its orientation). A system is also characterized by one or more inputs, which drive its dynamics and one or more outputs (e.g., for the
drone, the inputs could be the speeds of its rotators and the outputs the ones provided by an on-board monocular camera and/or a GPS).  A state is observable if the knowledge of the system inputs and outputs, during a given time interval, allows us its determination.

\noindent The concept of observability was first introduced for linear systems \cite{Kalman61,Kalman63} and
the analytic condition to check if a linear system satisfies this property has also been obtained.
\noindent The nonlinear case is much more complex. First, this concept becomes local. 
Second,
in general in the nonlinear case we can at most reconstruct the state only if we a priori know that it belongs to a given open set. In other words, by using the inputs and the outputs during an interval of time, we cannot distinguish, in general, states that are not close\footnote{This difference between the linear and the nonlinear case should not surprise. The observability concept is strongly related to the inverse function problem \cite{PHX}. A nonlinear function, even when its Jacobian is nonsingular, can take the same value at two or more separated points and, consequently, cannot be inverted. Invertibility becomes a local property.}.
This is the reason why, for nonlinear systems, the concept of {\it weak local observability} has been introduced \cite{Her77}.
In this paper, with the word observability we actually mean the concept of weak local observability, as defined in \cite{Her77,Casti82}
(definitions 8, 9, 10, 11, in \cite{Casti82}).
Third, unlike the linear case, observability depends on the system inputs.

\noindent The analytic condition to check if a nonlinear system satisfies this property has also been introduced \cite{Her77,Casti82,Suss83,Isi95}. It is known as the {\it observability rank condition}. The observability rank condition is a fundamental result that has extensively been used in many application domains, ranging from computer vision (e.g., \cite{Paul18}), 
robotics (e.g., \cite{Hesch14,Schiano18,Huang19}),
calibration (e.g., \cite{Guo13a}), mechanical engineering (e.g., \cite{Chat15})
up to biology (e.g., \cite{Villaverde19,Diste15,Villaverde16a,Chis11,Miao11}) and chemistry (e.g., \cite{August09}). It has also been used in the context of estimation theory \cite{Huang10,Huang13,Villaverde16}.
It is a simple systematic procedure that allows us to give an answer to the previous fundamental question, i.e., whether the state is observable or not. It is based on very simple systematic computation (differentiation and matrix rank determination) on the functions that describe the system.
The observability rank condition can deal with any system, independently of its complexity and type of nonlinearity. It works automatically.

\noindent On the other hand, the observability rank condition presents an important limitation: 
it does not account for the presence of unknown inputs. 
This is a severe limitation. 
The dynamics of most real systems are driven by inputs that are usually unknown. This holds in robotics, in biology, in chemistry, in physics, in economics, etc. For instance, in the case of our drone, its dynamics could be driven also by the wind, which is in general unknown. The wind is a disturbance and acts on the dynamics as an (unknown) input. The problem of unknown input observability (i.e., the problem of obtaining the analytical and automatic criterion that extends the observability rank condition to the case when some of the inputs are unknown) was defined long time ago \cite{Basile69,Guido71} and remained unsolved for half a century. 
The control theory community has spent a huge effort to design observers for both linear and nonlinear systems in the presence of unknown inputs, in many cases in the context of fault diagnosis, e.g., \cite{Guido71,Wang75,Bha78,Yan88,Gua91,Hou92,Daro94,Koe01,DeP01,Ha04,Floq04,Che06,Floq07,Barb07,Koe08,Barb09,Ham10}.
In some of the previous works, interesting conditions for the existence of an unknown input observer were introduced. On the other hand, these conditions have the following very strong impediments:

\begin{itemize}

\item They refer to a restricted class of systems. In particular, the considered systems are often characterized by linearity (or some specific type of nonlinearity)  with respect to the state in some of the functions that characterize the dynamics\footnote{These are the functions that appear in (\ref{EquationSystemDefinitionUIO}), i.e., $g^0(x),g^1(x), \ldots,g^{m_w}(x),f^1(x), \ldots,f^{m_u}(x)$.} and/or the system outputs\footnote{These are the functions $h_1(x),\cdots,h_p(x)$ that appear in equation (\ref{EquationSystemDefinitionUIO}).}.
No condition refers to any type of nonlinearity in the aforementioned functions.

\item They cannot be implemented automatically, i.e., by following a systematic procedure that does not require human intervention (e.g., by the usage of a simple code that adopts a symbolic computation tool). Most of these conditions are alternative definitions of the existence of a given observer, with limited interest from a practical point of view.

\end{itemize}

\noindent These limitations do not affect the observability rank condition in \cite{Her77,Isi95}. However, as we mentioned, this condition cannot be used in the presence of unknown inputs.
The solution of the unknown input observability problem is the extension of the observability rank condition to the unknown input case, i.e., the analytic condition able to provide the state observability in the presence of unknown inputs that does not encounter  the two aforementioned limitations. Additionally, the condition must characterize the system observability and, in this sense, it will be more general than a condition that checks the existence of an unknown input observer that belongs to a special class of observers.

\noindent Recently, the unknown input observability problem has been approached by introducing an extended state that includes the original state together with the unknown inputs and their time derivatives up to a given order \cite{Belo10,MED15,Maes19,Villa19b}. All these works proposed automatic iterative algorithms able to study the observability properties of nonlinear systems driven by also unknown inputs. On the other hand, all these algorithms suffer from the following fundamental limitations:

\begin{itemize}

\item They do not converge, automatically. 
In particular, at each iterative step, the state is extended by including new time derivatives of the unknown inputs. Consequently, if the extended state is observable at a given step, convergence is achieved. However, if this were not the case, we can never exclude that, at a later step, the extended state becomes observable. Therefore, in the presence of unobservability, all these algorithms remain inconclusive.

\item Due to the previous state augmentation, the computational burden becomes prohibitive after a few steps.

\end{itemize}

This paper provides the general analytical solution of the unknown input observability problem.
The solution does not need to extend the state and, as a result, does not encounter the above limitations. It provides a complete answer to the problem of state observability.
~Specifically,  these answers are automatically obtained by the usage of an algorithm that converges in a finite number of steps.

\noindent This paper is organized as follows. In chapter \ref{ChapterSystem} we provide the general characterization of the systems here investigated, together with some reminders on basic algebraic operations. We also provide the new concepts of {\it canonic system with respect to the unknown inputs}
and {\it canonical form with respect to the unknown inputs} (Definition \ref{DefinitionCanonicUI} and Definition \ref{DefinitionCanonicalForm}, respectively).
Note that, the solution provided in \cite{SIAMbook} is based on the assumption that the system is canonic with respect to its unknown inputs. 
In this paper we provide the solution for any system (even not canonic and not even canonizable). 
Chapter \ref{ChapterSystem}
ends by discussing the problem of observability by referring to an elementary example and by providing the solution in absence of unknown input (Algorithm \ref{AlgoObsTI0}), i.e., the standard observability rank condition.

Chapters \ref{ChapterSolutionCanonic1} and \ref{ChapterSolutionCanonic} provide the solution for systems that are in canonical form with respect to their unknown inputs. This is the solution of the same case dealt with in \cite{SIAMbook}. Chapter \ref{ChapterSolutionCanonic1} refers to the case of driftless, time invariant systems with a single unknown input (from now on, the simplified system), while Chapter \ref{ChapterSolutionCanonic} refers to the general (but canonic) case. The discussion of the simplified system is only given for educational purposes. Its solution is clearly a special case of the solution for the general canonic case.

The solution to the simplified system is given by Algorithm \ref{AlgoAbel}.
This solution differs from the solution presented in \cite{SIAMbook} and \cite{TAC19}, which is Algorithm \ref{AlgoBookAbel}.
The convergence criterion of Algorithm \ref{AlgoBookAbel} introduced in \cite{TAC19} does not hold always and the condition to characterize the systems for which the criterion does not hold (equation (31) in \cite{TAC19}) can be annoying for a practical use.
In \cite{SARAFRAZI}, some examples of simplified systems that violate the above condition (equation (31) in \cite{TAC19}) were shown. In addition, based on a brilliant derivation, a new criterion for the convergence of Algorithm  \ref{AlgoBookAbel} was introduced. This new criterion holds always. 
In this paper, in Chapter \ref{ChapterSolutionCanonic1}, we introduce a criterion similar to the one presented in \cite{SARAFRAZI}. The advantage of our criterion, compared with the one introduced in \cite{SARAFRAZI}, is that it allows us its extension to deal with the general canonic case. In other words, its validity is not limited to the simplified systems. Additionally, instead of Algorithm \ref{AlgoBookAbel}, we introduce a new algorithm (Algorithm \ref{AlgoAbel}), where the initialization step includes all the terms of Algorithm \ref{AlgoBookAbel} that make the convergence criterion of Algorithm  \ref{AlgoBookAbel} non trivial. This is preferable for a practical implementation. In practice, thanks to this initialization step, the convergence criterion of Algorithm \ref{AlgoAbel} is the same of the case without unknown inputs (i.e., Algorithm \ref{AlgoObsTI0}).


Chapter  \ref{ChapterSolutionCanonic}  extends all the aforementioned novelties to the general (but still canonic) case. The solution is now Algorithm \ref{AlgoNonAbel}, while the solution
introduced in \cite{SIAMbook}, is Algorithm \ref{AlgoBook}.
The convergence criterion provided in \cite{SIAMbook}, which is based on the computation of the tensor $\mathcal{T}$, does not hold always\footnote{In this case, the book does not even provide the analogue of Equation (31) in \cite{TAC19}, i.e., the equation that characterizes the systems for which this criterion holds. It is erroneously said that the criterion always holds, with only the exception of the systems that meet the assumption of Lemma 8.11}. In this paper, we provide the convergence criterion that holds in all the cases. Again,
instead of Algorithm \ref{AlgoBook}, we introduce a new algorithm (Algorithm \ref{AlgoNonAbel}), where the initialization step includes all the terms of Algorithm \ref{AlgoBook} that make the convergence criterion of Algorithm \ref{AlgoBook} non trivial.

In chapter \ref{ChapterSolutionNonCanonic} we remove the assumption that the system is in canonical form with respect to its unknown inputs.
This chapter provides the automatic procedure that, in a finite number of steps, provides all the observability properties.
This automatic procedure uses iteratively Algorithm \ref{AlgoNonAbel} and it is Algorithm \ref{AlgoFull}.  This new algorithm is a fundamental result, not presented in \cite{SIAMbook}. To this regard, note that the {\it Canonization} procedure introduced in Appendix C of \cite{SIAMbook} can fail. It is incomplete and there are systems that are not canonizable. Algorithm \ref{AlgoFull} is the complete and automatic solution of the unknown input observability problem in the most general case. Surprisingly, it is simpler than the procedure given in Appendix C of \cite{SIAMbook}. 


Chapter \ref{ChapterUIReconstruction} provides a fundamental result that regards a problem strongly related with the problem of state observability
in the presence of unknown inputs. This is the problem of unknown input reconstruction.

Chapter \ref{ChapterApplication}  illustrates the implementation of Algorithm \ref{AlgoFull} by studying the observability properties of a nonlinear system in the framework of visual-inertial sensor fusion. The dynamics of this system are driven by two unknown inputs and one known input and they are also characterized by a nonlinear drift. The system is not in canonical form with respect to its unknown inputs. However, by following the steps of Algorithm \ref{AlgoFull}, it can be set in canonical form.

Chapter \ref{ChapterConclusion} provides our conclusion.

\vskip.2cm
\noindent
In summary, the novelties of the complete solution introduced by this paper, with respect to the solution given in \cite{SIAMbook} are:

\begin{enumerate}

\item Full characterization of the concept of {\it canonicity with respect to the unknown inputs}. This characterization, includes the following new fundamental definitions:

\begin{itemize}

\item Definition of the unknown input reconstructability matrix and the unknown input degree of reconstructability (Definitions \ref{DefinitionRM} and \ref{DefinitionUIDegReconstrFromF}, respectively).

\item Definition of canonic system with respect to its unknown inputs and system in canonical form with respect to its unknown inputs (Definitions \ref{DefinitionCanonicUI} and \ref{DefinitionCanonicalForm}, respectively).

\item Definition of the highest unknown input degree of reconstructability (Definition \ref{DefinitionHDegUIReconstr}).

%

\end{itemize}

\item Algorithm \ref{AlgoFull}, which is the general solution that holds even in the non canonic case and not even canonizable. 
In particular, 
when the system is not canonizable, Algorithm \ref{AlgoFull} returns a new system with the highest unknown input degree of reconstructability, together with the observability codistribution.

\item A new criterion of convergence of the solution in the canonic case. 
In particular, 
the criterion proposed in \cite{SIAMbook}, which is based on the computation of the tensor $\mathcal{T}$, can fail. The new criterion here introduced, which extends the one introduced in \cite{SARAFRAZI} to the general case with drift, multiple unknown inputs and TV, holds always (and is even simpler). In addition, 
the algorithm that solves the problem is written in a new manner, where the initialization step includes all the terms of Algorithm \ref{AlgoBook} that make the convergence criterion of Algorithm \ref{AlgoBook} non trivial.

Finally, the content of this paper has been published by the Journal of Information Fusion \cite{IF22}.

%
%
%

\end{enumerate}

\chapter{System Characterization and the problem of Observability}\label{ChapterSystem}

\section{Basic equations}\label{SectionSystemEquation}

\noindent A general characterization of a dynamic system, which includes the presence of known and unknown inputs, is given by the following equations:

\begin{equation}\label{EquationSystemDefinitionUIO}
\left\{\begin{array}{ll}
  \dot{x} &=   g^0(x, t)+\sum_{k=1}^{m_u}f^k (x, t) u_k +  \sum_{j=1}^{m_w}g^j (x, t) w_j  \\
  y &= [h_1(x, t),\ldots,h_p(x, t)], \\
\end{array}\right.
\end{equation}

\noindent where:

\begin{itemize}

\item $x\in\mathcal{M}$ is the state and $\mathcal{M}$ is a differential manifold of dimension $n$.

\item $u_1,\ldots,u_{m_u}$ are the known inputs. They are $m_u$ independent functions of time that can be assigned. In control theory they are called {\it controls} or {\it control inputs}.

\item $w_1,\ldots,w_{m_w}$ are the unknown inputs or {\it disturbances}. They are $m_w$ independent functions of time. In particular, in this paper we assume that they are analytic functions of time.

\item $f^1,\ldots,f^{m_u},g^0,g^1,\ldots,g^{m_w}$ are $m_u+m_w+1$ vector fields, which are assumed to be smooth functions of $x$ and $t$ (the time).

\end{itemize}

\noindent The vector field $g^0$ is often called the {\it drift} since, when all  the inputs vanish, the state evolution is still non-vanishing in the presence of $g^0$.

\noindent In many cases the system is time-invariant (from now on TI), namely it has not an explicit time-dependence and all the functions that appear in (\ref{EquationSystemDefinitionUIO}) do not depend explicitly on time.  Nevertheless, we also account for an explicit  time dependence to be as general as possible. From now on, we use the acronym TV to indicate a system with an explicit time-dependence (TV stands for Time-Variant). We also use the acronym UI to mean unknown input and UIO to mean unknown input observability.

\noindent Finally, 
we show that the characterization given in (\ref{EquationSystemDefinitionUIO}) is very general. In particular, we show that, apparently more general characterizations, can be easily transformed into the above characterization.

A first category of more general systems would also include the presence of the inputs in the output functions, i.e., 

\[
\left\{\begin{array}{ll}
  \dot{x} &=   g^0(x, t)+\sum_{k=1}^{m_u}f^k (x, t) u_k +  \sum_{j=1}^{m_w}g^j (x, t) w_j  \\
  y &= [h_1(x, t,u,w),\ldots,h_p(x, t,u,w)], \\
\end{array}\right.
\]

\noindent where $u=[u_1,\ldots,u_{m_u}]^T$ and  $w=[w_1,\ldots,w_{m_w}]^T$. However, this case can be easily converted to (\ref{EquationSystemDefinitionUIO}) by including the inputs in the state and by considering the system driven by the new inputs $\dot{u}$ and $\dot{w}$.

A second category of more general systems would take into account for a general nonlinear dependence of the dynamics with respect to the inputs (both known and unknown). In accordance with Equation (\ref{EquationSystemDefinitionUIO}), this dependence is affine. 
This second category is characterized by the new dynamics:

\[
  \dot{x} =  f(x, t,u,w)
\]

Again, also this case can be easily converted to (\ref{EquationSystemDefinitionUIO}) by including the inputs in the state and by considering the system driven by the new inputs $\dot{u}$ and $\dot{w}$.

\section{Reminders on basic algebraic operations}\label{SectionNotation}

In this section, we provide some basic algebraic operations that will be widely used in the rest of this paper. 

We remind the reader of the Lie derivative operation and of some related fundamental properties. The Lie derivative evaluates the change of a tensor field\footnote{Note that, here, tensor is with respect to a coordinate change in our manifold $\mathcal{M}$. Later, the same concept of tensor will be referred to another group of transformations, which is the group of invariance of observability.} along a given vector field. In this paper we only use this operation when the tensor field is a scalar (i.e., a tensor of rank $0$), when it is a vector  (i.e., a tensor of rank $1$ and type $(0,~1)$) and when it is a row vector, or covector (i.e., a tensor of rank $1$ and type $(1,~0)$).
We use the following notation:

\vskip .5cm
\begin{itemize}

\item We denote by $\nabla$ the differential operator. If $h(x)$ is a scalar field defined on the manifold $\mathcal{M}$, $\nabla h=\frac{\partial}{\partial x}h$. This is evidently a row vector of dimension $n$.
Note that, in some cases we work in extended spaces that include the space of the states ($x$) and some of the UIs together with their time derivatives up to a given order. In this paper, with the symbol $\nabla$ we always mean the differential with respect to all the coordinates of the considered space. When we wish to consider the differential with respect to a given subset of coordinates, we adopt the symbol $\partial_{*}$, where $*$ stands for the considered subset of coordinates. For instance, if we work in the extended space that includes the original state $x$ together with the unknown inputs $w_1,\ldots,w_{m_w}$, the differential with respect to the original state alone (i.e., $x$) will be indicated by the symbol  $\partial_x=\left[\partial_{x^1},\ldots,\partial_{x^n}\right]$ and, in this case, $\nabla=\left[\partial_x,\partial_{w_1},\ldots,\partial_{w_{m_w}}\right]$.

\item Given a vector field $f$ (defined on $\mathcal{M}$),  $\mathcal{L}_f$ denotes the Lie derivative along $f$. When applied to the scalar field $h(x)$ we obtain the following new scalar field

\[
\mathcal{L}_fh\triangleq\nabla h \cdot f
\]
 (the product of a row vector times a column vector is a scalar). When applied to the covector field 
$\omega$, we obtain the covector field

\[
\mathcal{L}_f\omega\triangleq f^T\left(\frac{\partial\omega^T}{\partial x}\right)^T + \omega\frac{\partial f}{\partial x},
\]

 where $f^T$ is the transpose of $f$ (i.e., a row vector) and 
$\left(\frac{\partial\omega^T}{\partial x}\right)^T$ is the transpose of the Jacobian of $\omega^T$, and it is an $n\times n$ matrix.
An important special case occurs when $\omega$ is the differential of a scalar field, i.e., $\omega=\nabla h$. In this case, the above equation simplifies as follows:

\begin{equation}
\label{EquationLieDifferential}
\mathcal{L}_f\nabla h=\nabla \mathcal{L}_fh.
\end{equation}

 Finally, when applied to a vector field, $g$, we obtain the following new vector field:

\[
\mathcal{L}_fg\triangleq\frac{\partial g}{\partial x}f - \frac{\partial f}{\partial x}g=[f,~g],
\]

 where the parenthesis $[\cdot,~\cdot]$ are called {\it Lie brackets}.

\item In accordance with the control theory literature, we use the term {\it distribution} to denote the span\footnote{Note that, here, span is over the ring of the scalar functions which are smooth in $\M$ (or in a given open set of $\M$, when it is specified). In other words, any element of $\Delta$ is a vector field, $f$, that can be expressed as follows: $f=\sum_{i=1}^dc_i(x)f^i$, with $c_1(x), \ldots, c_d(x)$ scalar and smooth functions in $\M$. When we work in a given extended space, the ring is extended accordingly.} of a set of $d\le n$ vector fields, $\Delta=$span$\{
f^1,~f^2,\ldots,f^d\}$. It is basically a vector space 
that depends on $x\in\mathcal{M}$.

\item Similarly, we also use the term {\it codistribution} to denote the span of a set of $s\le n$ covector fields, $\Omega=$span$\{
\omega_1,~\omega_2,\ldots,\omega_s\}$. Again, it is a vector space that depends on $x\in\mathcal{M}$.

\item Given the codistribution $\Omega$ and a vector field $f$, we set $\mathcal{L}_f \Omega$ the span of all the covectors $\mathcal{L}_f\omega$ for any $\omega\in\Omega$.

\item Given two vector spaces $V_1$ and $V_2$,  $V_1{+}V_2$ is their sum, i.e., the span of all the generators of $V_1$ and $V_2$. When we have $k(>2)$ vector spaces $V_1,\ldots, V_k$, we denote its sum in the compact notation $\sum_{i=1}^kV_i$.

\end{itemize}

\section{A simple illustrative example}\label{SectionSystemExample}

\begin{figure}[htbp]
\begin{center}
\includegraphics[width=.6\columnwidth]{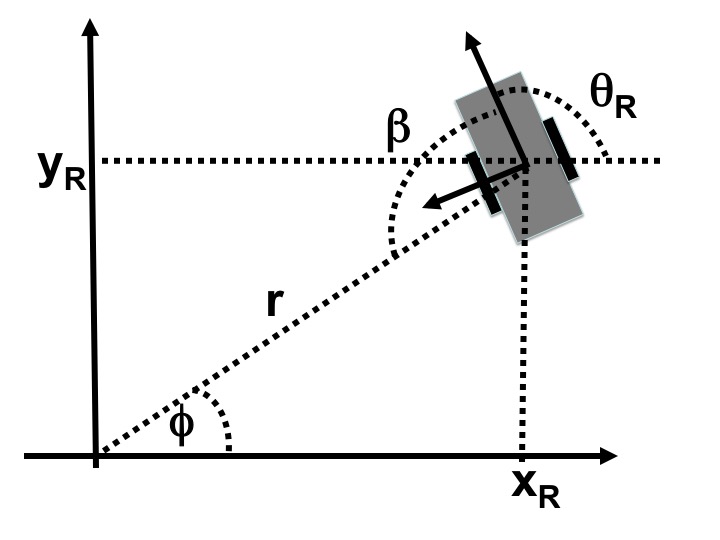}
\end{center}
\caption{Wheeled robot moving on a plane.} \label{Fig2Dvehicle}
\end{figure}

\noindent We provide a simple example to illustrate the above concepts and to highlight the goal of this paper.
We consider a wheeled robot that moves on a plane. By introducing on this plane a global frame, we can characterize the position and orientation  of the robot by the three parameters 
$x_R,~y_R,~\theta_R$ (see figure \ref{Fig2Dvehicle}). 

\noindent Under the assumption of the unicycle constraint, these parameters satisfy the following dynamics equations (unicycle dynamics)

\begin{equation}\label{EquationIntroductionExampleDynamics}
\left[\begin{array}{ll}
  \dot{x}_R &= v \cos\theta_R, \\
  \dot{y}_R &= v \sin\theta_R, \\
  \dot{\theta}_R &= \omega, \\
\end{array}\right.
\end{equation}

\noindent where $v$ and $\omega$ are the linear and the rotational
robot speed, respectively. 
We assume that a landmark is placed at the origin of the global frame. In addition, our robot is equipped with a sensor that perceives the landmark and provides its bearing angle in the local frame (i.e., the angle $\beta$ in figure \ref{Fig2Dvehicle}). 
This angle $\beta$ can be expressed in terms of the robot position and orientation. We have

\begin{equation}\label{EquationIntroductionExampleOutput}
\begin{array}{ll}
  \beta &= \pi-\theta_R+\textnormal{atan2}(y_R,~x_R).\\
\end{array}
\end{equation}

\noindent This system is a special case of the systems characterized by (\ref{EquationSystemDefinitionUIO}). In particular:

\begin{itemize}

\item $x=[x_R,~y_R,~\theta_R]^T$ is the state, it has dimension $n=3$ and it belongs to the manifold $\mathbb{R}^2\times\mathcal{S}^1$.

\item The vector $g^0$ is the zero 3-column vector.

\item  $m_u=2$ and $u_1=v$, $u_2=\omega$, $f^1=[\cos\theta_R,~\sin\theta_R,~0]^T$ and $f^2=[0,~0,~1]^T$.

\item $m_w=0$.

\item $y$ is the output and has a single component ($p=1$), i.e., $y=h_1(x)=\pi-\theta_R+\textnormal{atan2}(y_R,~x_R)$.

\end{itemize}

\begin{figure}[htbp]
\begin{center}
\begin{tabular}{ccc}
\includegraphics[width=.3\columnwidth]{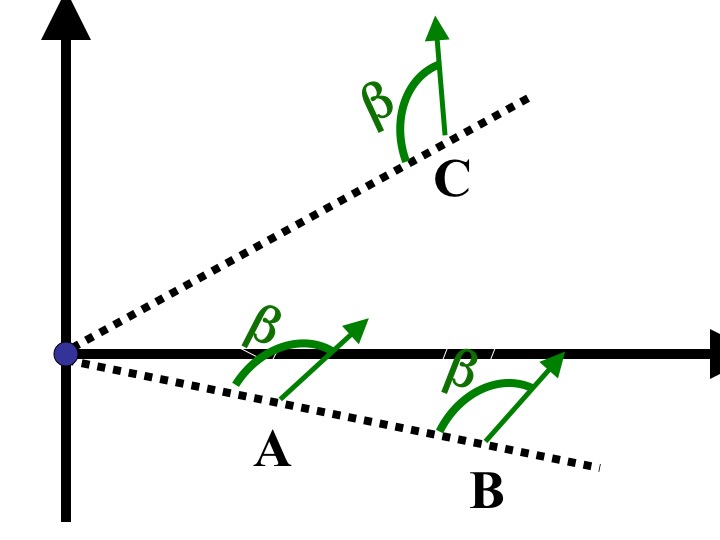}&
\includegraphics[width=.3\columnwidth]{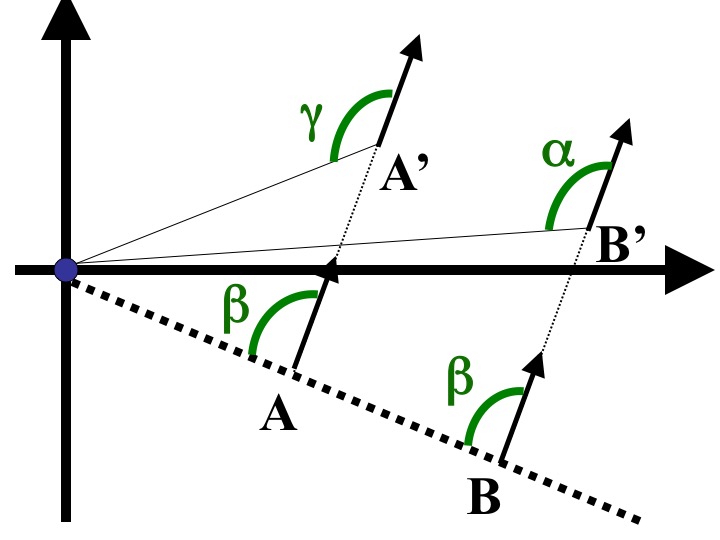}&
\includegraphics[width=.3\columnwidth]{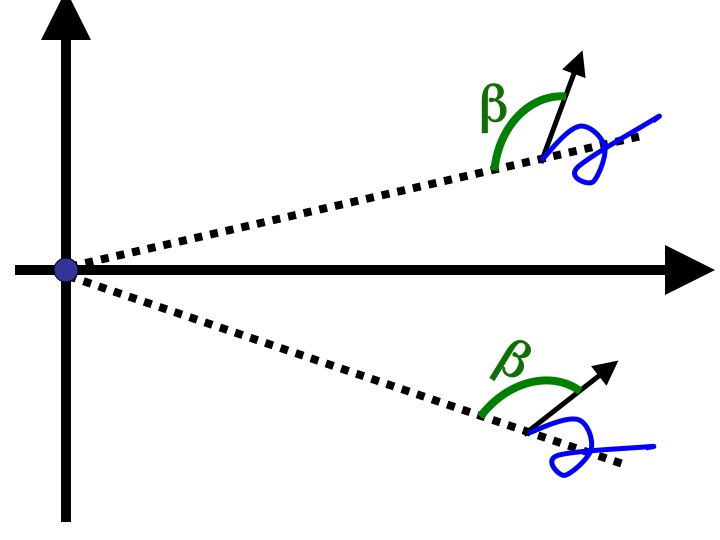}\\
$(a)$ & $(b)$& $(c)$
\end{tabular}
\caption{In $(a)$,
the three initial robot configurations are
compatible with the same initial observation ($\beta$).
In $(b)$, the two initial positions ($A$ and $B$) do not
reproduce the same observations ($\alpha\neq \gamma$). In $(c)$
the two indicated trajectories provide the same bearing
observations at every time.} \label{FigSimpleExamplebis}
\end{center}
\end{figure}

\noindent We study the observability properties of this system by following an intuitive procedure.
To check if the robot configuration $[x_R, ~y_R,
~\theta_R]^T$ is observable, we have to prove that it is
possible to uniquely reconstruct the initial robot configuration
by knowing the inputs and the outputs in a
given time interval. When, at the initial time, the bearing angle
$\beta$ of the origin is available, the robot can be everywhere
in the plane but, for each position, only one orientation
provides the right bearing $\beta$. In fig.
\ref{FigSimpleExamplebis}$a$ all the three positions $A$, $B$ and
$C$ are compatible with the observation $\beta$, provided that
the robot orientation satisfies
(\ref{EquationIntroductionExampleOutput}). In particular, the
orientation is the same for $A$ and $B$ but not for $C$.

\noindent Let us suppose that the robot moves according to the inputs
$v(t)$ and $\omega(t)$. With the exception of the special motion
consisting of a line passing by the origin, by only performing a
further bearing observation it is possible to distinguish all the
points belonging to the same line passing by the origin. In fig.
\ref{FigSimpleExamplebis}$b$ the two initial positions in $A$ and
$B$ do not reproduce the same observations after the movement, i.e., $\alpha \neq
\gamma$ (note that the segments $AA'$ and $BB'$ have the same length, which is known thanks to the knowledge of the system inputs). On the other hand, all the initial positions whose
distance from the origin is the same, cannot be distinguished
independently of the chosen trajectory. In fig.
\ref{FigSimpleExamplebis}$c$, the two indicated trajectories
provide the same bearing observations, at any time. Therefore,
the dimension of the {\it unobservable}
region is $1$. 
In particular, we introduce the following transformation

\begin{equation}\label{EquationIntroductionTransfomationIndSet}
\begin{array}{ll}
 x_R & \rightarrow x_R'= \cos\gamma ~x_R - \sin\gamma ~y_R,\\
 y_R & \rightarrow y_R'= \sin\gamma ~x_R + \cos\gamma ~y_R,\\
  \theta_R & \rightarrow \theta_R'=\theta_R+\gamma,\\
\end{array}
\end{equation}

\noindent where $\gamma\in[-\pi, ~\pi)$ is the parameter that defines the transformation.
The system inputs and output at any time are compatible with all the trajectories that differ because the initial state was transformed as above.

\noindent We wonder what is possible to reconstruct. The answer is immediate. All the physical quantities that are
 invariant with respect to the transform given in (\ref{EquationIntroductionTransfomationIndSet}). In other words, all the functions of the following two quantities:

\begin{equation}\label{EquationIntroductionObservableModes}
\begin{array}{ll}
  r&= r(x_R, ~y_R, ~\theta_R)=\sqrt{x_R^2+y_R^2},\\
\theta &=\theta(x_R, ~y_R, ~\theta_R)=\theta_R-\arctan2(y_R, x_R).\\
\end{array}
\end{equation}

\noindent Now, we assume that the first input, $v$, is unknown. In other words, we have $m_u=m_w=1$ and $g^1$ coincides with $f^1$ in the previous case and, the new $f^1$ is the old $f^2$. We immediately remark that, in addition to the invariance given by (\ref{EquationIntroductionTransfomationIndSet}), we have a new invariance, which is the scale. Indeed, both the known input ($\omega$) and the output ($\beta$) are angular measurements. The system has no source of metric information. Hence, the unobservable region is characterized by the further transformation

\begin{equation}\label{EquationIntroductionTransfomationIndSet2}
\begin{array}{ll}
 x_R & \rightarrow x_R'= \lambda ~x_R, \\
 y_R & \rightarrow y_R'= \lambda~y_R,\\
  \theta_R & \rightarrow \theta_R'=\theta_R,\\
\end{array}
\end{equation}

\noindent where $\lambda\in\mathbb{R}^+$ is the parameter that defines the transformation. In this case, we cannot reconstruct the distance $ r$ and we can only reconstruct any function of the only angle $\theta$ in (\ref{EquationIntroductionObservableModes}) since $\theta$ is invariant both with respect to (\ref{EquationIntroductionTransfomationIndSet}) and with respect to  (\ref{EquationIntroductionTransfomationIndSet2}) ($\arctan2(\lambda y_R, \lambda x_R)=\arctan2(y_R, x_R)$, for any  $\lambda\in\mathbb{R}^+$).

\noindent The above analysis was possible because the system is trivial: the state has dimension $3$ and, both the dynamics and the output have a very simple expression. We wish to provide the analytic tool able to automatically perform such analysis, i.e., by following a systematic analytical procedure. Note that, in the absence of unknown inputs, this procedure is the observability rank condition introduced in \cite{Her77} and summarized in Section \ref{SectionObservableFunction}. In the presence of a single unknown input and for driftless systems, the solution has very recently been obtained in \cite{TAC19}. Here, in Chapter \ref{ChapterSolutionCanonic1}, we provide a more efficient solution. Finally, in Chapters 
\ref{ChapterSolutionCanonic} and \ref{ChapterSolutionNonCanonic},
we provide the general solution (i.e., for driftless systems, with multiple unknown inputs and that holds even in the case TV). All these analytical systematic procedures use the algebraic operations summarized in section \ref{SectionNotation}.

\section{Observable function}\label{SectionObservableFunction}

When the state is not observable, there are functions of the state that can be observable, i.e., functions for which we have the possibility of reconstructing the value that they take at the initial state.
In the simple example discussed in section \ref{SectionSystemExample}, when all the inputs are known, we found that we cannot reconstruct the initial state. However, we can reconstruct the initial value of all the functions that depend on the state only via the two quantities $r$ and $\theta$ in (\ref{EquationIntroductionObservableModes}).
In other words, for this system all the observable functions are all the functions that depend on the state only via $r=\sqrt{x_R^2+y_R^2}$ and $\theta=\theta_R-\arctan2(y_R,~x_R)$ (e.g., a function of $x_R$ and/or of $\theta_R$ is not an observable function).
 When the input $v$ is unknown,  we can reconstruct the initial value of the functions that only depend on $\theta$. 
In this case, all the observable functions are all the functions that depend only on $\theta$.
The general definition of observable function is provided in \cite{TAC19}, where Definition 2 defines the observation space in the presence of unknown inputs, starting from the concept of indistinguishability. An observable function is precisely an element of this function space.
It is possible to define an observable function as follows\footnote{As we mentioned in the introduction, with observability we actually mean weak local observability.}.

\begin{definition}[Observable Function]\label{DefinitionObservableFunction}
The scalar function $\theta(x)$ is observable at $x_0\in\mathcal{M}$ if 
there exists an open neighbourhood $B$ of $x_0$ such that, by knowing that $x_0\in B$, 
it exists at least one choice of known inputs $u_1(t), \ldots,u_{m_u}(t)$ such that $\theta(x_0)$ can be  obtained from the knowledge of the output $y(t)$ and the inputs $u_1(t), \ldots,u_{m_u}(t)$ on the time interval $[t_0, ~t_0+T]$ for a given $T>0$. In addition, $\theta(x)$ is observable on a given set $\mathcal{A}\subseteq\mathcal{M}$, if it is observable at any $x\in \mathcal{A}$.
\end{definition}

This definition is the extension of the definition of state observability. Note that if we find $n$ independent observable functions (i.e., $n$ functions whose gradients are independent) the entire state is observable. More properties about the link of the above definition and the standard definition of state observability (based on the concept of indistinguishability) can be found in \cite{SIAMbook}.
The {\it observable codistribution}, or {\it observability codistribution}, is the span of the differentials of all the observable functions. By construction, the dimension of this codistribution cannot exceed the dimension of the state ($n$).
 If the dimension of this codistribution is $n$, every component of the state $x$ is observable.

\section{Canonic systems and canonical form  with respect to the unknown inputs}\label{SectionDefinitionCanonic}

We now introduce the new concepts of unknown input degree of reconstructability, canonic system with respect to its unknown inputs and canonical form with respect to the unknown inputs. 
This is obtained by introducing an important matrix that will be called the {\it Unknown Input Reconstructability Matrix} from a finite set of scalar functions.

\begin{definition}[Unknown input reconstructability matrix 
]\label{DefinitionRM}
Given the  $k$ scalar functions of the state, {  $\lambda_1(x),\ldots, \lambda_k(x)$}, the unknown input reconstructability matrix from {  $\lambda_1,\ldots, \lambda_k$} is defined as follows:

\begin{equation}\label{EquationRM}
\mathcal{RM}\left({  \lambda_1,\ldots, \lambda_k}\right)
\triangleq
\left[\begin{array}{cccc}
\Li_{g^1}{  \lambda_1} & \Li_{g^2}{  \lambda_1} & \ldots & \Li_{g^{m_w}}{  \lambda_1} \\
\Li_{g^1}{  \lambda_2} & \Li_{g^2}{  \lambda_2} & \ldots & \Li_{g^{m_w}}{  \lambda_2} \\
\ldots &\ldots &\ldots &\ldots \\
\Li_{g^1}{  \lambda_k} & \Li_{g^2}{  \lambda_k} & \ldots & \Li_{g^{m_w}}{  \lambda_k} \\
\end{array}
\right]
\end{equation}
\end{definition}

This matrix depends on $x$ and, for TV systems, on $t$. Based on this matrix, we introduce the following definition:

\begin{definition}[Unknown input degree of reconstructability
]\label{DefinitionUIDegReconstrFromF}
Given the state that satisfies the dynamics in (\ref{EquationSystemDefinitionUIO}), and the functions ${  \lambda_1,\ldots, \lambda_k}$, the unknown input degree of reconstructability from ${  \lambda_1,\ldots, \lambda_k}$ is the rank of $\mathcal{RM}\left({  \lambda_1,\ldots, \lambda_k}\right)$.
\end{definition}

By construction, the unknown input degree of reconstructability from any set of functions cannot exceed $m_w$. In addition, it depends in general on $x$ and, for TV systems, also on $t$.

Given the system characterized by (\ref{EquationSystemDefinitionUIO}), we consider two special cases, depending on the set of functions ${  \lambda_1,\ldots, \lambda_k}$:

\begin{itemize}

\item This set of functions consists of the output functions, i.e., $h_1,\ldots,h_p$.

\item This set of functions consists of all the independent observable functions\footnote{Here, with {\it independent} we mean that their differentials with respect to the state are independent. Hence, we can have at most $n$ independent observable functions (when there are $n$ independent observable functions the entire state is evidently observable).}.

\end{itemize}

In the first case, we refer to the {\it unknown input degree of reconstructability from the outputs}. In the second case, we omit to specify the functions and we refer to the  {\it unknown input degree of reconstructability}. 
%
We introduce the following definition:

\begin{definition}[Canonic system wrt the UIs]\label{DefinitionCanonicUI}
The system in (\ref{EquationSystemDefinitionUIO}) is canonic with respect to the unknown inputs if its unknown input degree of reconstructability is $m_w$.
\end{definition}

Let us consider the systems 
for which the unknown input degree of reconstructability from the outputs is $m_w$. Since the output functions are observable functions, these systems are canonic with respect to their unknown inputs. 
We say that these systems are in canonical form. In particular, we introduce the following definition:

\begin{definition}[System in canonical form wrt the UIs]\label{DefinitionCanonicalForm}
The system in (\ref{EquationSystemDefinitionUIO}) is in canonical form with respect to the unknown inputs if its unknown input degree of reconstructability from the outputs is $m_w$.
\end{definition}

Let us consider a system that satisfies (\ref{EquationSystemDefinitionUIO}). Its observability properties are the same of the new system defined as follows:

\begin{itemize}

\item It is characterized by the same state.

\item The time evolution of the state is the same (i.e., it satisfies the first equation of
(\ref{EquationSystemDefinitionUIO})).

\item The outputs are $h_1,\ldots,h_p,h_{p+1}$, where $h_{p+1}$ is an observable function for the first system.
\end{itemize}

Indeed, if the function $h_{p+1}$ is observable, its value can be reconstructed and, as a result, for the observability analysis, it can be set as an output of the system. For the simple example in Section \ref{SectionSystemExample}, it is easy to realize that, when both the inputs are known, we can add to the single output $h_1=\beta=\pi-\theta_R+\phi$, any function $h_{2}$ that depends on $r$ and/or $\theta_R-\phi$. The observability codistribution of the resulting system remains the same (e.g., we can add $h_2=r$).
Hence, for any system that is canonic with respect to its unknown inputs there exists a system that shares the same observability properties and that is in canonical form with respect to its unknown inputs. However, for a canonic system that is not in canonical form, a set of observable functions that makes the reconstructability matrix full rank, is not necessarily available. One of the outcomes of Algorithm \ref{AlgoFull}
is to automatically determine this set of functions.

Note that, as the rank of the matrix in (\ref{EquationRM}) depends in general on $x$ and, for TV systems, on $t$, all the above definitions are meant in a given open set where this rank takes the same value.

%
%
%
%

\section{The problem of Observability and the case without UI}\label{SectionAbsenceUI}

Given a dynamic system characterized by (\ref{EquationSystemDefinitionUIO}) the problem of observability is to obtain all the observable functions. 
In this paper we provide the general solution, which works automatically. In particular, we provide the algorithm that automatically computes a codistribution that includes the gradients of all the observable functions.
From now on, we call this codistribution the  {\it observablity codistribution}.
For educational purposes, instead of directly providing the algorithm that computes the observablity codistribution in the most general case, we first provide the algorithm in some special cases, which correspond to the characterization given in (\ref{EquationSystemDefinitionUIO}) when some of the parameters take special values. In particular, we consider separately the following cases:

\begin{enumerate}

\item Absence of unknown inputs. This case is characterized by (\ref{EquationSystemDefinitionUIO}) when $m_w=0$. The algorithm that computes the observablity codistribution is Algorithm \ref{AlgoObsTI0}, given in this section.

\item System in canonical form, in accordance with Definition \ref{DefinitionCanonicalForm}, TI, driftless, with a single UI and a single known input ($g^0$ equal to the zero column vector, $m_w=m_u=1$ and all the functions in (\ref{EquationSystemDefinitionUIO}) time-independent). The algorithm that computes the observablity codistribution is Algorithm \ref{AlgoAbel}, given in Chapter \ref{ChapterSolutionCanonic1}.

\item System in canonical form, in accordance with Definition \ref{DefinitionCanonicalForm}, both TI and TV, characterized by a nonlinear drift, any number of known and unknown inputs. The algorithm that computes the observablity codistribution is Algorithm \ref{AlgoNonAbel}, given in Chapter \ref{ChapterSolutionCanonic}.

\item System characterized by a nonlinear drift, any number of known and unknown inputs, but not necessarily in canonical form with respect to the unknwn inputs.
The algorithm that computes the observablity codistribution is Algorithm \ref{AlgoFull}. This Algorithm also computes the observablity codistribution when the system is not canonic with respect to its unknown inputs and not even canonizable. The analytical derivation is given in Chapter \ref{ChapterSolutionNonCanonic}.

\end{enumerate}

We denote by \Obs~ the codistribution computed by the algorithms \ref{AlgoObsTI0}, \ref{AlgoAbel}, \ref{AlgoNonAbel}, and  \ref{AlgoFull}, respectively in the four above cases.
The fact that \Obs~ is the observability codistribution is expressed by the following Theorem:

\begin{theorem}\label{TheoremBasic}
Let us consider a scalar function $\theta(x)$. If $\nabla\theta\in\OBS$ at $x_0\in\M$ then the function $\theta(x)$ is observable at $x_0$. Conversely, if the function $\theta$ is observable on an open set $\mathcal{A}\subseteq\M$ then $\nabla\theta\in\OBS$ in a dense set of $\mathcal{A}$.
\end{theorem}

\proof{

We distinguish the previous four cases.

\begin{itemize}

\item In the first case, the statement of the theorem is a well known result in the state of the art.  It can be easily obtained 
starting from theorem 3.1 and theorem 3.11 in \cite{Her77}. Its extension to the TV case is very simple and is available in \cite{ARCH20,TAC22} (it is also available in \cite{SIAMbook}).

\item In the second case, the proof of the statement is given in Chapter \ref{ChapterSolutionCanonic1}, starting from the results given in
\cite{TAC19} and in \cite{SIAMbook}. 

\item In the third case, 
the proof of the statement is given in Chapter \ref{ChapterSolutionCanonic}, starting from the results given in \cite{SIAMbook}.

\item In the fourth case, 
the proof of the statement is given in Chapter \ref{ChapterSolutionNonCanonic}.

\end{itemize}
$\blacktriangleleft$}

\noindent It is immediate to remark that, if the dimension of \Obs~ is $n$, every component of the state $x$ is observable.

\vskip.3cm

\noindent We conclude this section by providing the algorithm that computes \Obs~ for TI systems and in the absence of UIs. This is Algorithm \ref{AlgoObsTI0}.

\begin{al}[observability codistribution for TI nonlinear systems without unknown inputs]\label{AlgoObsTI0}
\begin{algorithmic}
  \State  $\Omega_0=\textnormal{span}\left\{\nabla  h_1,\ldots, \nabla h_p \right\}$ 
   \State  $\Omega_{k+1}=\Omega_k+\mathcal{L}_{g^0}\Omega_k+
   \sum_{i=1}^{m_u} \mathcal{L}_{f^i} \Omega_k$  
\end{algorithmic}
\end{al}

\noindent The above algorithm converges at the smallest integer $k$ such that $\Omega_k=\Omega_{k-1}=\OBS$. As the dimension of $\Omega$ cannot exceed $n$, this smallest integer cannot exceed $n$.

\vskip .3cm

\noindent We illustrate the use of this algorithm by referring to our illustrative example of section \ref{SectionSystemExample}. 

\noindent We run Algorithm \ref{AlgoObsTI0} to obtain \Obs. At the initialization we have:

\[
\Omega_0=\textnormal{span}\{\nabla h_1\}=\textnormal{span}\left\{\left[-\frac{y_R}{x_R^2+y_R^2}, ~\frac{x_R}{x_R^2+y_R^2},~-1\right]\right\}.
\]

\noindent By a direct computation we obtain that $ \mathcal{L}_{f_2}h_1=-1$ and, at the next step,  $\Omega=\textnormal{span}\{\nabla h_1, ~\nabla \mathcal{L}_{f_1}h_1\}$ (i.e., the dimension of $\Omega_1$ is $2$). At the next step, its dimension does not change and, consequently, the algorithm has converged. Therefore, \Obs$=\textnormal{span}\{\nabla h_1, ~\nabla \mathcal{L}_{f_1}h_1\}$.
We compute the orthogonal distribution (i.e., the span of the vectors orthogonal to $\nabla h_1$ and $\nabla \mathcal{L}_{f_1}h_1$, simultaneously). By a direct computation, we obtain: $\Delta=($\Obs$)^{\bot}=\textnormal{span}\{
[-y_R,~x_R,~1]^T\}$. It has dimension $1$. From the expression of its generator $[-y_R,~x_R,~1]^T$, we obtain the system invariance under the following infinitesimal transformation ($\epsilon$ is an infinitesimal parameter):

\begin{equation}\label{EquationIntroductionTransfomationIndSetinf}
\left[\begin{array}{c}
x_R \\
y_R \\
  \theta_R \\
\end{array}\right]
 \rightarrow
 \left[\begin{array}{c}
x_R' \\
y_R' \\
  \theta_R' \\
\end{array}\right]=\left[\begin{array}{c}
x_R \\
y_R \\
  \theta_R \\
\end{array}\right]+\epsilon
\left[\begin{array}{c}
 -y_R \\
 x_R\\
  1 \\
\end{array}\right]
\end{equation}

\noindent (see \cite{SIAMbook,TRO10} where we introduce, in this context, the concept of {\it continuous symmetry} that is any killing vector of \Obs). This is precisely the same invariance expressed by (\ref{EquationIntroductionTransfomationIndSet}) in the limit $\gamma\rightarrow\epsilon$ (the vector $[-y_R,~x_R,~1]^T$ is the generator of the Lie algebra associated to the one parameter Lie group described by (\ref{EquationIntroductionTransfomationIndSet})).

\noindent The above procedure can be easily extended to cope with TV systems \cite{SIAMbook,ARCH20,TAC22}. Specifically, in Algorithm \ref{AlgoObsTI0}, we simply need to introduce the following substitution:

\begin{equation}\label{EquationSubstitutionTV}
\mathcal{L}_{g^0}\rightarrow
\mathcal{L}_{g^0}+\frac{\partial}{\partial t}
\end{equation}

\noindent In the case of driftless systems, it means that, in the recursive step, we need to add the term $\oplusn \frac{\partial}{\partial t}\Omega$.


\chapter{Systems in Canonical Form: the Analytical Solution in the case of TI driftless systems  with a single UI}\label{ChapterSolutionCanonic1}

In this and in the next chapter we provide the analytical and automatic procedure that builds, in a finite number of steps, all the observable functions for systems characterized by (\ref{EquationSystemDefinitionUIO}) and that are in canonical form with respect to their unknown inputs, according to Definition \ref{DefinitionCanonicalForm}. 
In this chapter we provide the solution for driftless TI systems, with a single UI and a single known input. 
In other words, the systems studied in this chapter are characterized by the following equation:

\begin{equation}\label{EquationSystemDefinitionUIO1}
\left\{\begin{array}{ll}
  \dot{x} &=   f(x) u + g(x) w  \\
  y &= [h_1(x),\ldots,h_p(x)], \\
\end{array}\right.
\end{equation}

which is (\ref{EquationSystemDefinitionUIO}) with $m_w=m_u=1$, $f^1=f$, $g^1=g$, $g^0$ identically null, and all the functions time independent.
Then, in Chapter \ref{ChapterSolutionCanonic}, we relax all these restrictions. As we mentioned above, in these two chapters we assume that the systems are
in canonical form with respect to their unknown inputs, according to Definition \ref{DefinitionCanonicalForm}.
The general solution, i.e., for systems neither in canonical form nor canonic, will be given in Chapter \ref{ChapterSolutionNonCanonic}. The general solution uses, iteratively, the solution for systems in canonical form. This is the reason why we prefer to analyze this case before. 

\section{The solution}\label{SectionSolutionSolutionAbel}

In the specific case of this chapter, i.e., when the system is characterized by (\ref{EquationSystemDefinitionUIO1}), the fact that it is in canonical form means that there exists at least one output function among $h_1,\ldots,h_p$ such that its Lie derivative along $g$ does not vanish. 
We denote this function by $ h$. In addition, we denote by $L^1_g$ its Lie derivative along $g$,

\begin{equation}\label{EquationL1g}
L^1_g=\mathcal{L}_g  h.
\end{equation}

The solution of this case was obtained in \cite{SIAMbook,TAC19} and is returned by a simple algorithm that automatically computes \Obs~(Algorithm 7.2 in \cite{SIAMbook}) . However, in \cite{SIAMbook,TAC19}, the criterion of convergence of this algorithm is not general. Here, we provide a new algorithm, for which the convergence criterion is trivial (it is the same of Algorithm \ref{AlgoObsTI0}).

\noindent We denote by $\widehat{g}$ the following vector field in $\mathcal{M}$:

\begin{equation}\label{EquationHatg}
   \widehat{g}  = \frac{g}{L^1_g}.
\end{equation}

\noindent The observability codistribution is constructed by Algorithm \ref{AlgoAbel}. We denote by $\OBS_k$ the codistribution computed at the $k^{th}$ step, instead of $\Omega_k$, which is used in Algorithm 7.2 in \cite{SIAMbook}.

\begin{al}
\begin{algorithmic}
 \State

\State$\OBS_0=\tobs+\textnormal{span}\left\{\nabla  h_1,\ldots, \nabla h_p \right\}$ 
  \State$\OBS_{k+1}=
  \OBS_k+\mathcal{L}_{f} \OBS_k
+ \mathcal{L}_{\widehat{g}} \OBS_k$ 
\end{algorithmic}\label{AlgoAbel}
\end{al}

\noindent where $\tobs$ is an integrable codistribution. Its computation can be performed automatically and is provided in Section \ref{SubSectionTildeObs}. Algorithm \ref{AlgoAbel} converges at the smallest $k$ for which 

\begin{equation}\label{EquationConvergenceAlgo1}
\OBS_k=\OBS_{k-1},
\end{equation}
precisely as Algorithm \ref{AlgoObsTI0}. As the dimension of $\OBS$ cannot exceed $n$, this smallest integer cannot exceed $n$.

\subsection{Ingredients of the Solution}\label{SubSectionTildeObs}

The initialization step of Algorithm \ref{AlgoAbel} (i.e., the computation of $\OBS_0$) requires the computation of the codstribution $\tobs$. In this section, we define this codistribution and we provide the method to automatically compute it. It is based on the knowledge of two integers, denoted by $\NO$ and $\ND$. They are defined as follows:

\begin{itemize}

\item  $\NO$ is the smallest integer such that, when running Algorithm \ref{AlgoOmegagAbel}, the condition $\Omega^g_\NO=\Omega^g_{\NO-1}$ that characterizes its convergence, is achieved.
Note that $\NO\le n$.

\item  $\ND$ is the smallest integer such that, when running Algorithm \ref{AlgoDeltaAbel}, the condition $\Delta_{\ND+1}=\Delta_\ND$ that characterizes its convergence, is achieved.
Note that $\ND\le n-1$.

%

%
%

\end{itemize}

\subsubsection{The codistribution $\Omega^g$}\label{SubSubSectionOmegaG1}

$\Omega^g$ is the smallest codistribution that includes $\nabla h$ and is invariant with respect to $\Li_g$.
It is automatically constructed by Algorithm \ref{AlgoOmegagAbel}, which converges at the smallest  integer $\NO$ such that $\Omega^g_{\NO}=\Omega^g_{\NO-1}$($=\Omega^g$). 

\begin{al}[Codistribution $\Omega^g$]
\begin{algorithmic}
\State

\State$\Omega^g_0=\textnormal{span}\left\{\nabla h \right\}$
\State$\Omega^g_k=\Omega^g_{k-1}+\Li_g \Omega^g_{k-1}$
\end{algorithmic}\label{AlgoOmegagAbel}
\end{al}

Note that $\NO\le n$. In addition, $\Omega^g=\textnormal{span}\left\{\nabla h,~ \nabla\Li_gh,\ldots,\nabla\Li^{s-1}_gh\right\}$ and the functions $ h,~\Li_gh,\ldots,\Li^{s-1}_gh$ are called a basis of $\Omega^g$ (even if their gradients constitute, actually, a basis).
The computation of $\tobs$ needs to know the value of $\NO$. This is automatically obtained by running Algorithm \ref{AlgoOmegagAbel}.

%

\subsubsection{The distribution $\Delta$}\label{SubSubSectionDelta}

\noindent We introduce the following abstract operation:

\begin{definition}[Autobracket]\label{DefinitionAutobracket1}
Given the vector field $\phi$, we define its autobracket with respect to the system characterized by equation (\ref{EquationSystemDefinitionUIO1}) the following vector field
\begin{equation}\label{EquationAutobracket1}
 [\phi] = \frac{[g,~\phi]}{L^1_g},
\end{equation}
\end{definition}

\noindent Note that, in \cite{SIAMbook} and in \cite{TAC19} the same operation was defined by using $[\phi,~g]$ instead of $[g,~\phi]$ in (\ref{EquationAutobracket1})\footnote{A sign does not affect the observability properties}.

We denote by

\[
[\phi]^{(m)}
\]

the vector field obtained by applying the autobracket repetitively $m$ consecutive times, and

\[
[\phi]^{(0)}=\phi.
\]

We introduce the distribution $\Delta$. It is the smallest distribution that includes $f$ and is invariant under the autobracket operation.

\begin{al}[Distribution $\Delta$]
\begin{algorithmic}
\State
\State$\Delta_0=\textnormal{span}\left\{f \right\}$
\State$\Delta_k=\Delta_{k-1}+\left[\Delta_{k-1}\right]$
\end{algorithmic}\label{AlgoDeltaAbel}
\end{al}

with:

\[
\left[\Delta\right]\triangleq
\sum_{v\in\Delta}\textnormal{span}
\left\{
[v]
\right\}
\]

It is immediate to remark that, if $v_1,\ldots,v_l$ is a basis of $\Delta_{k-1}$, then $v_1,\ldots,v_l,[v_1],\ldots,[v_l]$ is a set of generators of $\Delta_{k-1}+\left[\Delta_{k-1}\right]$.

Algorithm \ref{AlgoDeltaAbel} converges at the smallest  integer $k$ such that $\Delta_k=\Delta_{k-1}$($=\Delta$). We set $k=\ND+1$. A basis of $\Delta$ is:

\[
\phi_0\triangleq f,~~\phi_1\triangleq [f],~~\phi_2\triangleq [f]^{(2)},~~\ldots,~~\phi_{\ND}\triangleq [f]^{(\ND)}
\]

Note that $\ND\le n-1$.
The computation of $\tobs$ needs to know the value of $\ND$. This is automatically obtained by running Algorithm \ref{AlgoDeltaAbel}.

\subsubsection{The codistribution $\tobs$}\label{SubSubSectiontobs1}

We are ready to provide $\tobs$. It is:

\begin{equation}\label{EquationTOBSDef1}
\tobs\triangleq\sum_{j=0}^{\NO+\ND}\textnormal{span}\left\{
\nabla\Li_{[f]^{(j)}}h
\right\}
\end{equation}

\section{Proof}

In this section we prove that Algorithm \ref{AlgoAbel}, with $\tobs$ defined in (\ref{EquationTOBSDef1}), provides the observability codistribution for systems characterized by (\ref{EquationSystemDefinitionUIO1}) and in canonical form with respect to the unknown input. 

We start from the solution given in \cite{SIAMbook} and \cite{TAC19}. We know that the observability codistribution is computed by Algorithm 7.2 in \cite{SIAMbook}, that is:

\begin{al}[Observability codistribution for TI driftless systems with a single UI]
\begin{algorithmic}
\State

\State$\Omega_0=\textnormal{span}\left\{\nabla  h_1,\ldots, \nabla h_p \right\}$ 
  \State$\Omega_{k+1}=
  \Omega_k+ \mathcal{L}_{f} \Omega_k
+ \mathcal{L}_{\widehat{g}} \Omega_k
+ \mathcal{L}_{[f]^k} \nabla h$ 
\end{algorithmic}\label{AlgoBookAbel}
\end{al}

Because of the last term at the recursive step, we cannot conclude that this algorithm converges. In particular, without this term, this algorithm would converge in at most $n-1$ steps, precisely as Algorithm \ref{AlgoObsTI0}.

The convergence of this algorithm is proved by the result stated by Theorem \ref{TheoremFundamentalNewAbel}. More importantly, this result states that we can obtain the observability codistribution by directly running Algorithm \ref{AlgoAbel}, which converges at the smallest integer $k$ such that $\OBS_{k-1}=\OBS_k$ and this occurs in at most $n-1$ steps.

The proof of Theorem \ref{TheoremFundamentalNewAbel} is based on the result stated by Proposition \ref{PropositionFundamentalNewAbel}. This is a technical result and it will be given, separately, in Section \ref{SectionPropositionFundamentalAbel}.

Finally, the proof of Theorem \ref{TheoremFundamentalNewAbel} requires the introduction of a new object, which will be denoted by $\C$ and is defined as follows.

By construction, we have:

\[
[\phi_{\ND}]\in\textnormal{span}\left\{
\phi_0,\phi_1,\ldots,\phi_{\ND}
\right\}
\]

We introduce the following quantity $\C$. It is characterized by two indices:

\[
\C_{km}
\]

and it is implicitly defined by the following:

\begin{equation}\label{EquationCAbelian}
[\phi_k]=\sum_{m=0}^{\ND} \C_{km}\phi_m, ~~~k\le \ND
\end{equation}

We trivially have:

\begin{equation}\label{EquationCAbelian1Exp}
\C_{km}=\delta_{k+1,m},~~k\le \ND-1
\end{equation}

In other words, only the last row of $\C$ ($\C_{\ND1},\C_{\ND2},\ldots,\C_{\ND\ND}$) is not trivial (in particular, the differential of any entry of $\C$ that does not belong to this row vanishes).

%
%
%

The following result proves that the codistribution \Obs, which is automatically computed by
Algorithm \ref{AlgoAbel}, is the observability codistribution. 

\begin{theorem}\label{TheoremFundamentalNewAbel}
There exists $\hat{k}$ such that, for any $k\ge\hat{k}$, $\OBS\subseteq\Omega_k$. Conversely, for any $k$, $\Omega_{k}\subseteq\OBS$.
\end{theorem}

\proof{Let us prove the first statement.
From the $(k+1)^{th}$ step of Algorithm \ref{AlgoBookAbel}, we obtain that $\mathcal{L}_{[f]^k} \nabla h\in\Omega_{k+1}$. 
From (\ref{EquationTOBSDef1}), we immediately obtain that
\[
\tobs\subseteq\Omega_{\ND+\NO+1}.
\]

By comparing the recursive step of Algorithm \ref{AlgoAbel} with the one of Algorithm \ref{AlgoBookAbel} it is immediate to conclude that, for any integer $q\ge0$:

\[
\OBS_q\subseteq\Omega_{\ND+\NO+1+q}
\]

and, as Algorithm \ref{AlgoAbel} converges in at most $n-1$ steps to $\OBS$, we conclude that:

\[
\OBS\subseteq\Omega_{\ND+\NO+n}
\]

and this proves the first statement with $\hat{k}=\ND+\NO+n\le 3n-1$.

Let us prove the second statement. We must prove that, for any $k$, $\Omega_{k}\subseteq\OBS$.
We proceed by induction. It is true at $k=0$ as $\Omega_0\subseteq\OBS_0\subseteq\OBS$.

Let us assume that
\[
\Omega_p\subseteq\OBS
\]

We must prove that 

\[
\Omega_{p+1}\subseteq\OBS
\]

As $\OBS$ is invariant under $\Li_{\widehat{g}}$ and $\Li_{f}$, we only need to prove that

\[
\nabla\Li_{[f]^{(p)}} h \in\OBS
\]

We distinguish the following two cases:

\begin{enumerate}

\item $p\le\ND+\NO$.

\item $p\ge \ND+\NO+1$.

\end{enumerate}

The first case is trivial as $\nabla\Li_{[f]^{(p)}} h\in\tobs\subseteq\OBS$.

Let us consider the second case. We set $p=\ND+\NO+l'$, with $l'\ge1$, and
we use Equation 
(\ref{EquationPropoAbel}), with $l=\NO+l'$. We have:

\begin{equation}\label{EquationTheoremAbel}
\nabla\Li_{[f]^{(p)}} h =
\gamma_{p-\ND}\sum_{m=0}^{\ND}\left(\nabla\C_{\ND m}\right)\Li_{\phi_m}\Li^{\NO+l'-1}_gh~~~\mod~\Omega_{p}+\Li_{\widehat{g}}\Omega_{p}
\end{equation}

with $\gamma_{p-\ND}\neq0$.


Let us denote by $t=\NO+l'-1\ge \NO$. We know that $\nabla\Li_g^th\in\Omega^g$ and, as a result,

\[
\nabla\Li_g^th=\sum_{i=1}^{\NO}c_i\nabla\Li_g^{i-1}h
\]

with $c_1,\ldots,c_{\NO}$ suitable scalar functions in the manifold $\M$. We have:

\[
\Li_{\phi_m}\Li^t_gh=\nabla\Li_g^th\cdot\phi_m=\left(\sum_{i=1}^{\NO}c_i\nabla\Li_g^{i-1}h\right)\cdot\phi_m=
\sum_{i=1}^{\NO}c_i\Li_{\phi_m}\Li^{i-1}_gh
\]

We substitute in (\ref{EquationTheoremAbel}) and we obtain:

\begin{equation}\label{EquationTheoremAbelAf}
\nabla\Li_{[f]^{(p)}} h =
\gamma_{p-\ND}
\sum_{i=1}^{\NO}c_i
\sum_{m=0}^{\ND}\left(\nabla\C_{\ND m}\right)\Li_{\phi_m}\Li^{i-1}_gh~~~\mod~\Omega_{p}+\Li_{\widehat{g}}\Omega_{p}
\end{equation}

Now we use again Equation 
(\ref{EquationPropoAbel}), with $l=i$. We have:

\[
\nabla\Li_{[f]^{(\ND+i)}} h =
\gamma_i\sum_{m=0}^{\ND}\left(\nabla\C_{\ND m}\right)\Li_{\phi_m}\Li^{i-1}_gh~~~\mod~\Omega_{\ND+i}+\Li_{\widehat{g}}\Omega_{\ND+i}
\]

As $\gamma_i\neq0$ and $\nabla\Li_{[f]^{(\ND+i)}} h\in\Omega_{\ND+i+1}$ , we have:

\[
\sum_{m=0}^{\ND}\left(\nabla\C_{\ND m}\right)\Li_{\phi_m}\Li^{i-1}_gh\in\Omega_{\ND+i+1},
\]

and

\[
\sum_{i=1}^{\NO}c_i\sum_{m=0}^{\ND}\left(\nabla\C_{\ND m}\right)\Li_{\phi_m}\Li^{i-1}_gh\in\Omega_{\ND+\NO+1}.
\]

As $p\ge{\ND+\NO+1}$, we have $\Omega_{\ND+\NO+1}\subseteq\Omega_p$, and we obtain:

\[
\sum_{i=1}^{\NO}c_i\sum_{m=0}^{\ND}\left(\nabla\C_{\ND m}\right)\Li_{\phi_m}\Li^{i-1}_gh\in\Omega_p.
\]

By using this in (\ref{EquationTheoremAbelAf}), we obtain:

\[
\nabla\Li_{[f]^{(p)}} h \in\Omega_p+\Li_{\widehat{g}}\Omega_{p}.
\]

From the inductive assumption, we know that $\Omega_p\subseteq\OBS$. In addition, $\OBS$ is invariant under $\Li_{\widehat{g}}$.
Therefore, $\nabla\Li_{[f]^{(p)}} h\in\OBS$.
 $\blacktriangleleft$}

\subsection{Intermediate technical results}\label{SectionPropositionFundamentalAbel}

The goal of this section is to prove the validity of Equation (\ref{EquationPropoAbel}), which has been used several times in the proof of Theorem \ref{TheoremFundamentalNewAbel}.
To achieve this goal, we basically need to proceed along two directions:

\begin{itemize}

\item Generalize the quantities $\C$ by introducing the quantities $\C^q$, according to Equation (\ref{EquationCAbelianGen}), and obtaining the recursive law stated by Lemma \ref{LemmaLawAbel}.

\item First proving the validity of Equation (\ref{EquationLemmaAbel}), which is more general than Equation (\ref{EquationPropoAbel}). In particular, this latter is obtained from the former in a special setting. 

\end{itemize}

We generalize the quantities $\C$ by introducing the quantities $\C^q$ ($q$ being an integer), as follows:

\begin{equation}\label{EquationCAbelianGen}
[\phi_k]^{(q)}=\sum_m \C_{km}^q\phi_m,~~~k\le r
\end{equation}

where $\sum_m$ stands for $\sum_{m=0}^r$.
The validity of the above expression comes from the fact that $[\phi_k]^{(q)}\in\textnormal{span}\left\{\phi_0,\ldots,\phi_r
\right\}$.

By construction,

\[
\C^1=\C
\]

We have the following result:

\begin{lm}\label{LemmaLawAbel}
The following equation holds:

\begin{equation}\label{EquationLawAbel}
\C^{l+1}_{km}=\Li_{\widehat{g}}\C^l_{km}+\sum_{m'}\C^l_{km'}\C_{m'm},
\end{equation}
and its differential expression:
\begin{equation}\label{EquationLawDiffAbel}
\nabla\C^{l+1}_{km}=\Li_{\widehat{g}}\nabla\C^l_{km}+\sum_{m'}\nabla\C^l_{km'}\C_{m'm}+\sum_{m'}\C^l_{km'}\nabla\C_{m'm}
\end{equation}

\end{lm}

\proof{
By definition, we have:

\begin{equation}\label{EquationLemmaStep1}
[\phi_k]^{(l+1)}= \sum_m\C_{km}^{l+1}\phi_m
\end{equation}

On the other hand,

\[
[\phi_k]^{(l+1)}= \left[[\phi_k]^{(l)}\right]=
\left[
\sum_m\C_{km}^{l}\phi_m
\right]=
\sum_m\left(\Li_{\widehat{g}}\C_{km}^{l}\right)\phi_m+
\sum_m\C_{km}^l[\phi_m]=
\]

\[
\sum_m\left(\Li_{\widehat{g}}\C_{km}^{l}\right)\phi_m+
\sum_{mm'}\C_{km}^l\C_{mm'}\phi_{m'}=
\sum_m\left(\Li_{\widehat{g}}\C_{km}^{l}\right)\phi_m+
\sum_{m'm}\C_{km'}^l\C_{m'm}\phi_{m}
\]

By comparing with (\ref{EquationLemmaStep1})
we obtain (\ref{EquationLawAbel}) and, by differentiating, (\ref{EquationLawDiffAbel}).
$\blacktriangleleft$}

\begin{lm}\label{LemmaFundamentalNewAbel}
Given $\phi_k$ ($k=0,\ldots,r$), for any integer $l\ge1$, and any integer $0\le j\le l-1$, we have
\begin{equation}\label{EquationLemmaAbel}
\nabla\Li_{[\phi_k]^{(l)}} h = \beta_{l,j}\sum_m\left(\nabla\C^{l-j}_{km}\right)\Li_{\phi_m}\Li^j_gh~~~\mod~\Omega_{r+l}+\Li_{\widehat{g}}\Omega_{r+l}
\end{equation}
with $\beta_{l,j}\neq0$.

\end{lm}

\proof{
We proceed by induction on $l$. Let us set $l=1$. We only have $j=0$. We have:

\[
\nabla\Li_{[\phi_k]^{(1)}} h=
\nabla\Li_{\sum_m\C_{km}\phi_m} h=\nabla\sum_m\C_{km}\Li_{\phi_m}h=
\sum_m\C_{km}\nabla\Li_{\phi_m} h+\sum_m\left(\nabla\C_{km}\right)\Li_{\phi_m}h
\]

On the other hand, for $m\le r$

\[
\nabla\Li_{\phi_m} h\in\Omega_{r+1},
\]

Hence,

\[
\nabla\Li_{[\phi_k]^{(1)}} h=
\sum_m\left(\nabla\C_{km}\right)\Li_{\phi_m}h,~~~~~\mod~ \Omega_{r+1},
\]

which proves (\ref{EquationLemmaAbel}) when $l=1$, $j=0$ and $\beta_{1,0}=1\neq0$.

Let us consider the recursive step. We assume that (\ref{EquationLemmaAbel}) holds at a given $l=l^*>1$. We have:

\begin{equation}\label{EquationLemmaAbel*}
\nabla\Li_{[\phi_k]^{(l^*)}} h = \beta_{l^*,j}\sum_m\left(\nabla\C^{l^*-j}_{km}\right)\Li_{\phi_m}\Li^j_gh~~~\mod~\Omega_{r+l^*}+\Li_{\widehat{g}}\Omega_{r+l^*}
\end{equation}

for any $j=0,1,\ldots,l^*-1$ and with $\beta_{l^*,j}\neq0$. In addition, (\ref{EquationLemmaAbel}) holds at any $l\le l^*$ for any $j\le l-1$. In particular, it holds for any 
$l\le l^*$ and $j=l-1$. Hence we also have:

\begin{equation}\label{EquationLemmaAbellm1}
\nabla\Li_{[\phi_k]^{(l)}} h = \beta_{l,l-1}\sum_m\left(\nabla\C_{km}\right)\Li_{\phi_m}\Li^{l-1}_gh~~~\mod~\Omega_{r+l}+\Li_{\widehat{g}}\Omega_{r+l}
\end{equation}

for any $l\le l^*$ with $\beta_{l,l-1}\neq0$.

We must prove the validity of (\ref{EquationLemmaAbel}) at $l^*+1$ and for any $j=0,\ldots,l^*$, i.e.:

\begin{equation}\label{EquationLemmaAbel*p1}
\nabla\Li_{[\phi_k]^{(l^*+1)}} h = \beta_{l^*+1,j} \sum_m\left(\nabla\C^{l^*-j+1}_{km}\right)\Li_{\phi_m}\Li^j_gh~~~\mod~\Omega_{r+l^*+1}+\Li_{\widehat{g}}\Omega_{r+l^*+1}
\end{equation}

We proceed by induction on $j$. Let us consider $j=0$. We have:

\[
\nabla\Li_{[\phi_k]^{(l^*+1)}} h = \nabla\Li_{\sum_m\C_{km}^{l^*+1}\phi_m} h=\nabla\sum_m\C_{km}^{l^*+1}\Li_{\phi_m}h=
\sum_m\C_{km}^{l^*+1}\nabla\Li_{\phi_m} h+\sum_m\left(\nabla\C_{km}^{l^*+1}\right)\Li_{\phi_m}h
\]

On the other hand, for $m\le r$

\[
\sum_m\C_{km}^{l^*+1}\nabla\Li_{\phi_m} h\in\Omega_{r+1}\subseteq\Omega_{r+l^*+1}.
\]

This proves the validity of (\ref{EquationLemmaAbel*p1}) at $j=0$ with $\beta_{l^*+1,0}=1$.

Let us assume that (\ref{EquationLemmaAbel*p1}) holds at a given $j=j^*\le l^*-1$. We must prove that it also holds at $j=j^*+1$. We prove the equality, $\mod~\Omega_{r+l^*+1}+\Li_{\widehat{g}}\Omega_{r+l^*+1}$, of the two expressions at $j^*$ and $j^*+1$, namely, we must prove the following:

\[
\beta_{l^*+1,j^*} \sum_m\left(\nabla\C^{l^*-j^*+1}_{km}\right)\Li_{\phi_m}\Li^{j^*}_gh=
\beta_{l^*+1,j^*+1} \sum_m\left(\nabla\C^{l^*-j^*}_{km}\right)\Li_{\phi_m}\Li_g\Li^{j^*}_gh
~~~\mod~\Omega_{r+l^*+1}+\Li_{\hat{g}}\Omega_{r+l^*+1}
\]

From the inductive assumption at $j^*$ we know that $\beta_{l^*+1,j^*}\neq0$. We denote $\rho\triangleq\frac{\beta_{l^*+1,j^*+1}}{\beta_{l^*+1,j^*}}$. As a result, we must prove:

\[
\sum_m\left(\nabla\C^{l^*-j^*+1}_{km}\right)\Li_{\phi_m}\Li^{j^*}_gh=
\rho \sum_m\left(\nabla\C^{l^*-j^*}_{km}\right)\Li_{\phi_m}\Li_g\Li^{j^*}_gh
~~~\mod~\Omega_{r+l^*+1}+\Li_{\widehat{g}}\Omega_{r+l^*+1}
\]

with $\rho\neq0$.

We adopt the notation:

\[
\lambda\triangleq\Li^{j^*}_gh,~~~~~~
\C^-=\C^{l^*-j^*},~~~~~~
\C^+=\C^{l^*-j^*+1}
\]

We must prove:

\begin{equation}\label{EquationLemmaFundamentalAbelProof}
\sum_m\left(\nabla\C^+_{km}\right)\Li_{\phi_m}\lambda=
\rho \sum_m\left(\nabla\C^-_{km}\right)\Li_{\phi_m}\Li_g\lambda
~~~\mod~\Omega_{r+l^*+1}+\Li_{\widehat{g}}\Omega_{r+l^*+1}
\end{equation}

with $\rho\neq0$.

From (\ref{EquationLawDiffAbel}) we have:

\[
\nabla\C^+_{km}=\Li_{\widehat{g}}\nabla\C^-_{km}+\sum_{m'}\nabla\C^-_{km'}\C_{m'm}+\sum_{m'}\C^-_{km'}\nabla\C_{m'm}
\]

We substitute in the above expression and we must prove:

\[
\sum_m\left(
\Li_{\widehat{g}}\nabla\C^-_{km}+\sum_{m'}\nabla\C^-_{km'}\C_{m'm}+\sum_{m'}\C^-_{km'}\nabla\C_{m'm}
\right)\Li_{\phi_m}\lambda=
\]
\[
\rho \sum_m\left(\nabla\C^-_{km}\right)\Li_{\phi_m}\Li_g\lambda
~~~\mod~\Omega_{r+l^*+1}+\Li_{\widehat{g}}\Omega_{r+l^*+1}
\]

with $\rho\neq0$.

Let us consider the first term on the left hand side.

\[
\sum_m
\left(\Li_{\widehat{g}}\nabla\C^-_{km}
\right)
\Li_{\phi_m}\lambda=
\Li_{\widehat{g}}
\left(
\sum_m
\nabla\C^-_{km}
\Li_{\phi_m}\lambda\right)-
\sum_m
\nabla\C^-_{km}
\Li_{\widehat{g}}
\left(
\Li_{\phi_m}\lambda
\right)
\]

Note that $j^*\le l^*-1$. Hence, we are allowed to use (\ref{EquationLemmaAbel*}) at $j=j^*$.
In the new notation, it tells us that

\[
\nabla\Li_{[\phi_k]^{(l^*)}} h = \beta_{l^*,j^*}\sum_m\left(\nabla\C^-_{km}\right)\Li_{\phi_m}\lambda~~~\mod~\Omega_{r+l^*}+\Li_{\widehat{g}}\Omega_{r+l^*}
\]

On the other hand,

\[
\nabla\Li_{[\phi_k]^{(l^*)}} h  \in\Omega_{k+l^*+1}
\subseteq\Omega_{r+l^*+1}
\]

As $\beta_{l^*,j^*}\neq0$ we obtain:

\[
\Li_{\widehat{g}}
\left(
\sum_m
\nabla\C^-_{km}
\Li_{\phi_m}\lambda\right)\in\Li_{\widehat{g}}\Omega_{r+l^*+1}
\]
 
and we have:

\[
\sum_m\left(
\Li_{\widehat{g}}\nabla\C^-_{km}+\sum_{m'}\nabla\C^-_{km'}\C_{m'm}+\sum_{m'}\C^-_{km'}\nabla\C_{m'm}
\right)\Li_{\phi_m}\lambda=
\]
\[-
\sum_m
\nabla\C^-_{km}
\Li_{\widehat{g}}
\left(
\Li_{\phi_m}\lambda
\right)
+\sum_{mm'}\nabla\C^-_{km'}\C_{m'm}\Li_{\phi_m}\lambda+\sum_{mm'}\C^-_{km'}\nabla\C_{m'm}\Li_{\phi_m}\lambda
~\mod~\Omega_{r+l^*+1}+\Li_{\widehat{g}}\Omega_{r+l^*+1}
\]

Let us consider the last term of the above. We consider (\ref{EquationLemmaAbellm1}). It holds for any $l\le l^*$. As $j^*\le l^*-1$, we are  allowed to set $l=j^*+1$. We have:

\[
\nabla\Li_{[\phi_k]^{(j^*+1)}} h = \beta_{j^*+1,j^*}\sum_m\left(\nabla\C_{km}\right)\Li_{\phi_m}\Li^{j^*}_gh
= \beta_{j^*+1,j^*}\sum_m\left(\nabla\C_{km}\right)\Li_{\phi_m}\lambda~~~\mod~\Omega_{r+j^*+1}+\Li_{\widehat{g}}\Omega_{r+j^*+1}
\]

As $\beta_{j^*+1,j^*}\neq0$,

\[
\sum_m\left(\nabla\C_{km}\right)\Li_{\phi_m}\lambda\in\Omega_{r+j^*+2}
\]

As $j^*\le l^*-1$, this proves that:

\[
\sum_m\left(\nabla\C_{km}\right)\Li_{\phi_m}\lambda\in\Omega_{r+l^*+1}
\]

We have:

\[
-
\sum_m
\nabla\C^-_{km}
\Li_{\widehat{g}}
\left(
\Li_{\phi_m}\lambda
\right)
+\sum_{mm'}\nabla\C^-_{km'}\C_{m'm}\Li_{\phi_m}\lambda
+\sum_{mm'}\C^-_{km'}\nabla\C_{m'm}\Li_{\phi_m}\lambda=
\]
\[
-
\sum_m
\nabla\C^-_{km}
\Li_{\widehat{g}}
\left(
\Li_{\phi_m}\lambda
\right)
+\sum_{mm'}\nabla\C^-_{km'}\C_{m'm}\Li_{\phi_m}\lambda
~~~\mod~\Omega_{r+l^*+1}
\]

We have

\[
-
\sum_m
\nabla\C^-_{km}
\Li_{\widehat{g}}
\left(
\Li_{\phi_m}\lambda
\right)
+\sum_{mm'}\nabla\C^-_{km'}\C_{m'm}\Li_{\phi_m}\lambda=
-
\sum_m
\nabla\C^-_{km}
\Li_{\widehat{g}}
\left(
\Li_{\phi_m}\lambda
\right)
+\sum_{mm'}\nabla\C^-_{km}\C_{mm'}\Li_{\phi_{m'}}\lambda=
\]

\[
\sum_m
\nabla\C^-_{km}
\left(
-\Li_{\widehat{g}}
\left(
\Li_{\phi_m}\lambda
\right)
+\sum_{m'}\C_{mm'}\Li_{\phi_{m'}}\lambda
\right)=
\sum_m
\nabla\C^-_{km}
\left(
-\frac{1}{L^1_g}\Li_g
\Li_{\phi_m}\lambda
+\Li_{[\phi_{m}]}\lambda
\right)=
\]
\[
\frac{1}{L^1_g}
\sum_m
\nabla\C^-_{km}
\left(
-\Li_g
\Li_{\phi_m}\lambda
+\Li_{[g,~\phi_{m}]}\lambda
\right)=
-\frac{1}{L^1_g}
\sum_m
\nabla\C^-_{km}
\Li_{\phi_{m}}\Li_g\lambda,
\]

which coincides with the right side of (\ref{EquationLemmaFundamentalAbelProof}) with $\rho=-\frac{1}{L^1_g}\neq0$.
$\blacktriangleleft$}

We have the following fundamental result:

\begin{pr}\label{PropositionFundamentalNewAbel}
For any integer $l\ge1$, we have
\begin{equation}\label{EquationPropoAbel}
\nabla\Li_{[\phi_r]^{(l)}} h =
\nabla\Li_{[f]^{(r+l)}} h =
\gamma_l\sum_{m=0}^{r}\left(\nabla\C_{rm}\right)\Li_{\phi_m}\Li^{l-1}_gh~~~\mod~\Omega_{r+l}+\Li_{\widehat{g}}\Omega_{r+l}
\end{equation}
with $\gamma_l\neq0$.
\end{pr}

\proof{The above equality is obtained from Lemma \ref{LemmaFundamentalNewAbel}, with $j=l-1$, $k=r$, and $\gamma_l=\beta_{l,l-1}$.
$\blacktriangleleft$}

\section{Application}

\noindent We conclude this chapter by discussing our illustrative example introduced in section \ref{SectionSystemExample}, when the linear speed $v$ is unknown. Note that in \cite{TAC19} we provide the same example. On the other hand, in \cite{TAC19} we computed the observability codistribution by using Algorithm \ref{AlgoBookAbel} and by using a criterion of convergence that does not hold always. Here, we use Algorithm \ref{AlgoAbel}, whose convergence criterion is the same as in the classical case (i.e., without unknown inputs).

As in \cite{TAC19}, it is more convenient to adopt polar coordinates, in which the expression of the output becomes very simple (the result by using Cartesian coordinates would be the same but its derivation more laborious).

Our system will be characterized by the state (see Figure  \ref{Fig2Dvehicle}):

\[
x=[r, ~\phi, ~\theta_R]^T
\]

\noindent where $r=\sqrt{x_R^2+y_R^2}$, and $\phi=\tan^{-1}\left(\frac{y_R}{x_R}\right)$

Note that, in these coordinates, we compute the gradient by differentiating with respect to $\rho$, $\phi$ and $\theta_R$. By using the dynamics in (\ref{EquationIntroductionExampleDynamics}) we obtain the following dynamics:

\begin{equation}\label{EquationObsExampleDynamics}
\left[\begin{array}{ll}
  \dot{r} &= v \cos(\theta_R-\phi) \\
  \dot{\phi} &= \frac{v}{r} \sin(\theta_R-\phi) \\
  \dot{\theta}_R &= \omega \\
  y&=\phi-\theta_R\\
\end{array}\right.
\end{equation}

The preliminary step is obtained by comparing (\ref{EquationObsExampleDynamics}) with (\ref{EquationSystemDefinitionUIO1}). We obtain: $p=1, w=v,~u=\omega$,

\[
g(x)=\left[\begin{array}{c}
\cos(\theta_R-\phi)\\
\frac{\sin(\theta_R-\phi)}{r}\\
0\\
\end{array}
\right],~~~~~~
f(x)=\left[\begin{array}{c}
0\\
0\\
1\\
\end{array}
\right],~~~~~~h(x)=\phi-\theta_R
\]

\noindent We compute the observability codistribution by running algorithm \ref{AlgoAbel}. We need, first of all, to compute $L^1_g=\Li_gh$ and $\widehat{g}$. We have:

\[
L^1_g=-\frac{\sin(\phi - \theta_R)}{r},~~~~~~\widehat{g}(x)=\left[\begin{array}{c}
-\frac{r\cos(\phi - \theta_R)}{\sin(\phi - \theta_R)}\\
1\\
0\\
\end{array}
\right].
\]

Then, we must compute the codistributions $\tobs$. 
We need to compute the two integers $s$ and $r$. 
Let us start by computing $s$.
We run Algorithm \ref{AlgoOmegagAbel}. We obtain $\Omega^g_2=\Omega^g_1$. Hence, $s=2$.

Let us compute $r$. We run 
Algorithm \ref{AlgoDeltaAbel} and we obtain: $\Delta_3=\Delta_2$. Hence, $r+1=3$ and $r=2$.

We compute $\tobs$ by using (\ref{EquationTOBSDef1}). We need to compute $\phi_0$, $\phi_1$, $\phi_2$, $\phi_3$, and$\phi_4$. We obtain:

\[
\phi_0=f=\left[\begin{array}{c}
0\\
0\\
1\\
\end{array}
\right],~~~
\phi_1=\left[\begin{array}{c}
r\\
\frac{\cos(\phi - \theta_R)}{\sin(\phi - \theta_R)}\\
0\\
\end{array}
\right],~~~
\phi_2=\left[\begin{array}{c}
-\frac{2r\cos(\phi - \theta_R)}{\sin(\phi - \theta_R)}\\
\frac{2\sin^2(\phi - \theta_R) - 2}{\sin^2(\phi - \theta_R)}\\
0\\
\end{array}
\right],
\]
\[
\phi_3=\left[\begin{array}{c}
\frac{6r-4r\sin^2(\phi - \theta_R)}{\sin^2(\phi - \theta_R)}\\
\frac{2\cos(\phi - \theta_R) + 4\cos^3(\phi - \theta_R)}{\sin^3(\phi - \theta_R)}\\
0\\
\end{array}
\right],~~~
\phi_4=\left[\begin{array}{c}
-\frac{8r\cos(\phi - \theta_R)(\cos(\phi - \theta_R)^2 + 2)}{\sin^3(\phi - \theta_R)}\\
8\frac{\cos(2\phi - 2\theta_R) - 10\cos(4\phi - 4\theta_R) - \cos(6\phi - 6\theta_R) + 10}{15\cos(2\phi - 2\theta_R) - 6\cos(4\phi - 4\theta_R) + \cos(6\phi - 6\theta_R) - 10}\\
0\\
\end{array}
\right].
\]

Hence:

\[
\Li_{\phi_0}h=-1,~~
\Li_{\phi_1}h=\frac{\cos(\phi - \theta_R)}{\sin(\phi - \theta_R)},~~
\Li_{\phi_2}h=\frac{2\sin^2(\phi - \theta_R) - 2}{\sin^2(\phi - \theta_R)},
\]
\[
\Li_{\phi_3}h=\frac{2\cos(\phi - \theta_R) + 4\cos^3(\phi - \theta_R)}{\sin^3(\phi - \theta_R)},~~
\Li_{\phi_4}h=8\frac{\cos(2\phi - 2\theta_R) - 10\cos(4\phi - 4\theta_R) - \cos(6\phi - 6\theta_R) + 10}{15\cos(2\phi - 2\theta_R) - 6\cos(4\phi - 4\theta_R) + \cos(6\phi - 6\theta_R) - 10}.
\]

By computing their differential, it is immediate to obtain $\tobs$. We can verify that $\tobs\subseteq\textnormal{span}\left\{\nabla h\right\}$.
This concludes the initialization step of Algorithm \ref{AlgoAbel}. We obtain:

\[
\OBS_0=\textnormal{span}\left\{\nabla h
\right\}=\textnormal{span}\left\{
[0,~1,~-1]
\right\}
\]

By executing the recursive step of Algorithm \ref{AlgoAbel}, we obtain $\OBS_1=\OBS_0$ and we conclude that the algorithm has converged and the observability codistribution is:

\[
\OBS=\textnormal{span}\left\{
[0,~1,~-1]
\right\}
\]

Its orthogonal distribution is:

\[
(\textnormal{\Obs})^{\bot}=\textnormal{span}\{
[0,~1,~1]^T,~~[1,~0, ~0]^T\}.
\]

\noindent It has dimension $2$. 
The first generator, $[0,~1,~1]^T$, characterizes the system invariance described by the infinitesimal transformation in (\ref{EquationIntroductionTransfomationIndSetinf}), but expressed in the polar coordinates. The second generator, $[1,~0, ~0]^T$, characterizes the following invariance:

\begin{equation}\label{EquationIntroductionTransfomationIndSetinfScale}
\left[\begin{array}{c}
r \\
\phi \\
  \theta_R \\
\end{array}\right]
 \rightarrow
 \left[\begin{array}{c}
r' \\
\phi' \\
  \theta_R' \\
\end{array}\right]=\left[\begin{array}{c}
r \\
\phi \\
  \theta_R \\
\end{array}\right]+\epsilon
\left[\begin{array}{c}
 1 \\
 0\\
  0 \\
\end{array}\right].
\end{equation}

\noindent This is precisely the same invariance expressed by (\ref{EquationIntroductionTransfomationIndSet2}) in the limit $\lambda\rightarrow\epsilon$ and in the polar coordinates (the vector $[1,~0,~0]^T$ is the generator of the Lie algebra associated to the one parameter Lie group described by (\ref{EquationIntroductionTransfomationIndSet2})).


\chapter{Systems in Canonical Form: the Analytical Solution in the case of TV systems, in the presence of a drift, and multiple UIs}\label{ChapterSolutionCanonic}

In this chapter we provide the analytical and automatic procedure that builds, in a finite number of steps, all the observable functions for systems characterized by (\ref{EquationSystemDefinitionUIO}) and that are in canonical form with respect to their unknown inputs, according to Definition \ref{DefinitionCanonicalForm}. With respect to Chapter \ref{ChapterSolutionCanonic1}, the solution is not limited to driftless systems with a single unknown input and a single known input and the functions in (\ref{EquationSystemDefinitionUIO}) can also explicitly depend on time.

\section{The solution}\label{SectionSolutionNonAbel}

We are dealing with systems in canonical form with respect to their unknown inputs. From Definition \ref{DefinitionRM} and Definition \ref{DefinitionCanonicalForm},
we know that we can extract from the output $m_w$ functions, which we denote by $\tih_1,\ldots,\tih_{m_w}$, such that the reconstructability matrix is full rank. We set:

\begin{equation}\label{EquationTensorMSynchro}
\mu^i_j = \mathcal{L}_{g^i}\tih_j, ~~~i,~ j=1,\ldots,m_w
\end{equation}

\noindent Note that, the above object ($\mu$) is characterized by two indices. On the other hand, these indices are unrelated to the dimension ($n$) of our manifold ($\mathcal{M}$) and, consequently, they cannot be tensorial indices with respect to a coordinates' change in $\mathcal{M}$. In addition, with respect to a coordinates' change, $\mu^i_j$ behaves as a scalar field, for any $i,~j$.
On the other hand, and this is the fundamental key to achieving a profound understanding, $\mu$ is a two index tensor of type $(1,~1)$ with respect to a new group of transformations. This is the group of invariance of observability, introduced in \cite{SIAMbook} and denoted by \SUIO.
Actually, the indices of the tensors associated to \SUIO~take $m_w+1$ values. The extra value is set to $0$ and, in the rest of this paper, we adopt the Einstein notation where, Latin indices take the values $1,2,\ldots,m_w$ and Greek indices take the values $0,1,2,\ldots,m_w$ (the reader can find a brief summary of tensor calculus in the second chapter of \cite{SIAMbook}). The tensor $\mu$ must be completed as follows:

\begin{equation}\label{EquationTensorMInertial2}
\mu^0_0 = 1,~
\mu^i_0=0,~
\mu^0_i = \frac{\partial\tih_i}{\partial t}  +
\mathcal{L}_{g^0}\tih_i,~  i=1,\ldots,m_w.
\end{equation}

\noindent We denote by $\nu$ the inverse of $\mu$. In other words, we have

\begin{equation}\label{EquationTensorNe}
\mu^{\alpha}_{\gamma}\nu^{\gamma}_{\beta}=\delta^{\alpha}_{\beta},~~~~~~~~\alpha,~\beta=0,1,\ldots,m_w,
\end{equation}

\noindent where $\delta^{\alpha}_{\beta}$ is the delta Kronecker and, in accordance with the Einstein notation, the dummy Greek index $\gamma$ is summed over $\gamma=0,1,\ldots,m_w$. We have:

\begin{equation}\label{EquationTensorNInertial}
\nu^0_0 = 1,~
\nu^i_0 = 0,~
\nu^0_i = -\mu^0_k\nu^k_i, ~
\nu^i_k \mu^k_j=\mu^i_k \nu^k_j=\delta^i_j.
\end{equation}

\noindent As for Greek indices, when Latin indices are dummy, they imply a sum. For Latin indices the sum is over $1,\ldots,m_w$ (this is the case of $k$ in the above equation). Since $\nu$ is the inverse of a tensor field of type $(1,~1)$ associated to \SUIO, it is a tensor of type $(1,~1)$ .

\noindent Note that, in our theory, we deal with two types of tensors, simultaneously.
This because observability is invariant both with respect to a coordinates' change in $\mathcal{M}$ and with respect to \SUIO\footnote{As we mentioned in the introduction, this is similar 
to what happens in the Standard Model 
of particle physics.
Also in that case, two types of tensors coexist, simultaneously: tensors with respect to the global Poincar\'e symmetry and tensors with respect to the local 
$SU(3)\times SU(2)\times U(1)$ 
gauge symmetry.}. To avoid confusion, we use a matrix format for the tensors associated to a coordinates' change in $\mathcal{M}$ and we use indices only for \SUIO. In this notation, $g^\alpha$, for a given $\alpha=0,1,\ldots,m_w$, is a vector (i.e., tensor of type $(0,~1)$) with respect to a coordinates' change and the component $\alpha$ of a tensor of type $(0,~1)$ with respect to \SUIO. 

\noindent We set

\begin{equation}\label{Equationh0t}
\tilde{h}_0=t,
\end{equation}

\noindent where $t$ is the time. $\tilde{h}_\alpha$ , for a given $\alpha=0,1,\ldots,m_w$, is a scalar (i.e., tensor of type $(0,~0)$) with respect to a coordinates' change and the component $\alpha$ of a tensor of type $(1,~0)$ with respect to \SUIO.

\noindent We denote by $\widehat{g}^\alpha$ the following vector fields in $\mathcal{M}$:

\begin{equation}\label{Equationgalpha}
   \widehat{g}^{\alpha}  = \nu^{\alpha}_{\beta} g^{\beta},~~\alpha=0,1,\ldots,m_w.
\end{equation}

\noindent Since they are obtained by an index contraction of two tensor fields of \SUIO, 
$\widehat{g}^{\alpha}$ is the component $\alpha$ of a tensor of type $(0,~1)$ with respect to \SUIO.

\noindent The observability codistribution is constructed by Algorithm \ref{AlgoNonAbel}.

\begin{al}[\Obs~ for systems characterized by (\ref{EquationSystemDefinitionUIO}), in canonical form with respect to the UIs]
\begin{algorithmic}
\State
\State$\OBS_0=\tobs+\textnormal{span}\left\{\nabla  h_1,\ldots, \nabla h_p \right\}$ 
  \State$\OBS_{k+1}=
  \OBS_k+\sum_{i=1}^{m_u}\mathcal{L}_{f^i} \OBS_k
+ \sum_{\beta=0}^{m_w}\mathcal{L}_{\widehat{g}^\beta} \OBS_k$ 
\end{algorithmic}\label{AlgoNonAbel}
\end{al}

\noindent where $\tobs$ is an integrable codistribution. Its computation can be performed automatically and is provided in Section \ref{SubSectionTildeObsNonAbel}. 
Algorithm \ref{AlgoNonAbel} can be easily extended to cope with TV systems. Specifically, we simply need to introduce the following substitution at the first term in the sum at the recursive step (the term with $\beta=0$):

\begin{equation}\label{EquationSubstitutionTGV}
\mathcal{L}_{\widehat{g}^0}\rightarrow
\dt{\mathcal{L}}_{\widehat{g}^0}
\triangleq
\mathcal{L}_{\widehat{g}^0}+\frac{\partial}{\partial t}
\end{equation}


Algorithm \ref{AlgoNonAbel} converges at the smallest $k$ for which 

\begin{equation}\label{EquationConvergenceAlgo}
\OBS_k=\OBS_{k-1}
\end{equation}

 and $k\le n-m_w+1$ (note that the dimension of $\OBS_0$ is at least $m_w$). The limit codistribution will be denoted by $\OBS$.

\subsection{Ingredients of the Solution}\label{SubSectionTildeObsNonAbel}

The initialization step of Algorithm \ref{AlgoNonAbel} (i.e., the computation of $\OBS_0$) requires the computation of the codstribution $\tobs$. In this section, we define this codistribution and we provide the method to automatically compute it. It is based on the knowledge of two integers, denoted by $\NO$ and $\ND$. They are defined as follows:

\begin{itemize}

\item  $\NO$ is the smallest integer such that, when running Algorithm \ref{AlgoOmegagNonAbel}, the condition $\Omega^g_\NO=\Omega^g_{\NO-1}$ that characterizes its convergence, is achieved.
Note that $\NO\le n-m_w+1$.

\item  $\ND$ is the smallest integer such that, when running Algorithm \ref{AlgoDeltaNonAbel}, the condition $\Delta_{\ND+1}=\Delta_\ND$ that characterizes its convergence, is achieved.
Note that $\ND\le n-1$ (it is even
$\ND\le n-\textnormal{dim}\left\{\textnormal{span}\left\{f^1,\ldots,f^{m_u} \right\}\right\}$).

%
%

%
%

\end{itemize}

\subsubsection{The codistribution $\Omega^g$}\label{SubSubSectionOmegaG}

$\Omega^g$ is the smallest codistribution that includes $\nabla \tih_1,\ldots,\nabla \tih_{m_w}$ and is invariant with respect to $\Li_{g^\beta}$ ($\beta=0,1,\ldots,m_w$). It is automatically constructed by Algorithm \ref{AlgoOmegagNonAbel}, which converges at the smallest  integer $\NO$ such that $\Omega^g_{\NO}=\Omega^g_{\NO-1}$($=\Omega^g$).

\begin{al}[Codistribution $\Omega^g$]
\begin{algorithmic}
\State
\State$\Omega^g_0=\sum_{j=1}^{m_w}\textnormal{span}\left\{\nabla\tih_j \right\}$
\State$\Omega^g_k=\Omega^g_{k-1}+\sum_{\beta=0}^{m_w}\Li_{g^\beta} \Omega^g_{k-1}$
\end{algorithmic}\label{AlgoOmegagNonAbel}
\end{al}

In the TV case, we simply need to replace:  $\Li_{g^0}\rightarrow\dt{\mathcal{L}}_{g^0}\triangleq
\mathcal{L}_{g^0}+\frac{\partial}{\partial t}
$.

Note that $\NO\le n-m_w+1$.
The computation of $\tobs$ needs to know the value of $\NO$. This is automatically obtained by running Algorithm \ref{AlgoOmegagNonAbel}.
%
%
%
%

%

\subsubsection{The distribution $\Delta$}\label{SubSubSectionDelta}

\noindent We introduce the following abstract operation:

\begin{definition}[Autobracket]\label{DefinitionAutobracket}
Given the quantity $\phi^{\alpha_1,\ldots,\alpha_k}$, which is a vector field with respect to a change of coordinates in $\mathcal{M}$ and one component of a tensor of type $(0,~k)$ with respect to \SUIO, we define its autobracket with respect to the system characterized by equation (\ref{EquationSystemDefinitionUIO}) the following object:
\[
 [\phi^{\alpha_1,\ldots,\alpha_k}]^{\alpha_{k+1}} = \nu^{\alpha_{k+1}}_{\beta} [g^{\beta},~\phi^{\alpha_1,\ldots,\alpha_k}]+\nu^{\alpha_{k+1}}_0\frac{\partial\phi^{\alpha_1,\ldots,\alpha_k}}{\partial  t}=
\]
\begin{equation}\label{EquationAutobracket}
 \nu^{\alpha_{k+1}}_{\beta} [g^{\beta},~\phi^{\alpha_1,\ldots,\alpha_k}]+\delta^{\alpha_{k+1}}_0\frac{\partial\phi^{\alpha_1,\ldots,\alpha_k}}{\partial  t},
\end{equation}
\noindent where, the square brackets on the right hand side are the Lie brackets and the dummy Greek index $\beta$ is summed over $\beta=0,1,\ldots,m_w$ (in accordance with the Einstein notation).
\end{definition}

\noindent The output of this operation is still a vector field with respect to a coordinates' change and one component of a tensor of type $(0,~k+1)$ with respect to \SUIO. The autobracket increases by one the tensor rank with respect to \SUIO.

\noindent Note that, in the case of driftless systems with a single UI (dealt with in Chapter \ref{ChapterSolutionCanonic1}), the tensor $\mu$ has only one non trivial component which is $\mu^1_1=L^1_g$ and $\nu^1_1=\frac{1}{L^1_g}$. In addition,  $\widehat{g}^1$ becomes $\frac{g}{L^1_g}$, which is the vector field in (\ref{EquationHatg}).
Finally, the autobracket operation reduces to the Autobracket defined in Definition \ref{DefinitionAutobracket1}  (in this case, we do not have indices with respect to \SUIO~ and, as it is shown in \cite{SIAMbook}, this is because \SUIO~becomes an Abelian group).

We can apply the autobracket repetitively. We denote by

\[
[\phi]^{(\alpha_1,\ldots,\alpha_m)}
\]

the vector field obtained by applying the autobracket repetitively $m$ consecutive times, first along $\alpha_1$ and last along $\alpha_m$ (these operations are not commutative).

We introduce the distribution $\Delta$. It is the smallest distribution that includes $f^1,\ldots,f^{m_u}$ and is invariant under the autobracket operation.

\begin{al}[Distribution $\Delta$]
\begin{algorithmic}
\State
\State$\Delta_0=\textnormal{span}\left\{f^1,\ldots,f^{m_u} \right\}$
\State$\Delta_k=\Delta_{k-1}+\sum_\beta\left[\Delta_{k-1}\right]^\beta$
\end{algorithmic}\label{AlgoDeltaNonAbel}
\end{al}
with:

\[
\left[\Delta\right]^\beta\triangleq
\sum_{v\in\Delta}\textnormal{span}
\left\{
[v]^\beta
\right\},
\]

for any $\beta=0,1,\ldots,m_w$.
It is immediate to remark that, if $v_1,\ldots,v_l$ is a basis of $\Delta_{k-1}$, then $v_1,\ldots,v_l,[v_1]^\beta,\ldots,[v_l]^\beta$ is a set of generators of $\Delta_{k-1}+\left[\Delta_{k-1}\right]^\beta$.

Algorithm \ref{AlgoDeltaNonAbel} converges at the smallest  integer $k$ such that $\Delta_k=\Delta_{k-1}$($=\Delta$). We set $k=\ND+1$.
Note that $\ND\le n-1$ (it is even
$\ND\le n-\textnormal{dim}\left\{\Delta_0\right\}$).
The computation of $\tobs$ needs to know the value of $\ND$. This is automatically obtained by running Algorithm \ref{AlgoDeltaNonAbel}.

\subsubsection{The codistribution $\tobs$}\label{SubSubSectiontobs}

We are ready to provide $\tobs$. It is:

\begin{equation}\label{EquationTOBSDef}
\tobs\triangleq
\sum_{q=1}^{m_w}
\sum_{j=0}^{\NO+\ND}
\sum_{\alpha_1,\ldots,\alpha_j}
\sum_{i=1}^{m_u}
\textnormal{span}\left\{
\nabla\Li_{[f^i]^{(\alpha_1,\ldots,\alpha_j)}}\tih_q
\right\}
\end{equation}

where the sum $\sum_{\alpha_1,\ldots,\alpha_j}$ stands for $\sum_{\alpha_1=0}^{m_w}\ldots\sum_{\alpha_j=0}^{m_w}$

\section{Proof}\label{SectionProofNonAbel}

In this section we prove that Algorithm \ref{AlgoNonAbel}, with $\tobs$ defined in (\ref{EquationTOBSDef}), provides the observability codistribution for systems characterized by (\ref{EquationSystemDefinitionUIO}) and in canonical form with respect to the unknown inputs.

We start from the solution given in \cite{SIAMbook}. We know that the observability codistribution is computed by Algorithm 8.2 in \cite{SIAMbook}, that is:

\begin{al}
\begin{algorithmic}
\State
  \State$\Omega_0=\textnormal{span}\left\{\nabla  h_1,\ldots\nabla h_p\right\}$ 

 \State$\Omega_{k+1}=\Omega_k+ \sum_{i=1}^{m_u}\Li_{f^i} \Omega_k + \sum_{\beta=0}^{m_w}\Li_{\widehat{g}^\beta} \Omega_k +
  \sum_{q=1}^{m_w}\sum_{i=1}^{m_u}\sum_{\alpha_1=0}^{m_w}\ldots\sum_{\alpha_k=0}^{m_w}
 \textnormal{span} \left \{\mathcal{L}_{[f^i]^{(\alpha_1,\ldots,\alpha_k)}} \nabla \tih_q \right \}$
\end{algorithmic}\label{AlgoBook}
\end{al}

Because of the last term at the recursive step, we cannot conclude that this algorithm converges. In particular, without this term, this algorithm would converge in at most $n-m_w$ steps, precisely as Algorithm \ref{AlgoObsTI0}.

The convergence of this algorithm is proved by the result stated by Theorem \ref{TheoremFundamentalNewNonAbel}. More importantly, this result states that we can obtain the observability codistribution by directly running Algorithm \ref{AlgoNonAbel}, which converges at the smallest integer $k$ such that $\OBS_{k-1}=\OBS_k$ and $k\le n-m_w+1$.

The proof of Theorem \ref{TheoremFundamentalNewNonAbel} is based on the result stated by Proposition \ref{PropositionFundamentalNewNonAbel}. This is a technical result and it will be given, separately, in Section \ref{SectionPropositionFundamentalNonAbel}.

Finally, the proof of Theorem \ref{TheoremFundamentalNewNonAbel} requires the introduction of
two fundamental new ingredients:

\begin{itemize}

\item The set of generators $\psi$.

\item The quantities $\C^\alpha_{\textbf{k}\textbf{m}}$.

\end{itemize}

\subsubsection{The set of generators $\psi$}

We denote by $\psi$ the vectors computed by Algorithm \ref{AlgoDeltaNonAbel} before its convergence. 
Note that these vectors generate the entire distribution $\Delta$ and they will be not necessarily independent. They will be listed by using a double index
$\mathbf{i}\triangleq (i,j_i)$. The first index is the index of the step of Algorithm \ref{AlgoDeltaNonAbel}. It takes the values $i=0,1,\ldots,r$. The second index lists all the vectors which are computed by Algorithm \ref{AlgoDeltaNonAbel}, when computing $\Delta_i$.
For this proof, we do not need to specify how we list these vectors. We only emphasize that at the initialization (computation of $\Delta_0$) we have $m_u$ vectors: $\psi_{(0,1)}=f^1,~\psi_{(0,2)}=f^2,~\ldots,~\psi_{(0,m_u)}=f^{m_u}$. In other words, the second index $j_0$ take the values $1,2,\ldots,m_u$. Then, when computing $\Delta_i$, we have, for every $\psi_{i-1,j_{i-1}}$, $m_w+1$ vectors with the first index equal to $i$ (each of them is obtained by computing the autobracket $[\psi_{i-1,j_{i-1}}]^\alpha$, for $\alpha=0,1,\ldots,m_w$). 
Therefore, the index $j_i$, take the values: $1,2,\ldots,(m_w+1)^im_u$.

%
%
%
%
%
%
%
%
%
%

The vectors $\psi_\mathbf{k}=\psi_{(k,j_k)}$ for $k=0,1,\ldots,r$, generate the entire distribution $\Delta$. However, they are not independent, in general.

\subsubsection{The quantities $\C^\alpha_{\textbf{k}\textbf{m}}$}

We introduce the following quantities $\C$. They are characterized by three indices:

\[
\C_{\mathbf{kz}}^\alpha
\]

The lower indices $\mathbf{k}$ and $\mathbf{z}$ are double indices, as explained above ($\mathbf{k}=(k,j_k)$ and $\mathbf{z}=(z,j_z)$, with $k,z=0,1,\ldots,r$). The upper index $\alpha$ takes the $m_w+1$ values: $0,1,\ldots,m_w$.
These quantities are implicitly defined by the following equation:

\begin{equation}\label{EquationDefinitionC}
[\psi_\mathbf{k}]^\alpha=\sum_\mathbf{z}\C_{\mathbf{kz}}^\alpha\psi_\mathbf{z}
\end{equation}
where $\sum_\mathbf{z}=\sum_{(z,j_z)}$ stands for $\sum_{z=0}^r\sum_{j_z=1}^{(m_w+1)^zm_u}$.
Indeed, the term at the left hand side of the above equation belongs to $\Delta$ and can be expressed in terms of the set of vectors $\psi$ (we remind the reader that $\Delta$ is invariant under the autobracket operation, which is the operation $[\cdot]^\alpha$ at the left side of (\ref{EquationDefinitionC})).
As the vectors $\psi$ are not necessarily independent, the choice of $\C$ is not unique. When the index $\mathbf{k}=(k,j_k)$ is such that $k<r$ we trivially choose $\C$ such that it selects the vector $\psi$ with the first index incremented by one.
%
%
%
For $\mathbf{k}=(r,j_r)$ the expression of $\C_{\mathbf{kz}}^\alpha$ is not trivial and provides $[\psi_\mathbf{k}]^\alpha$ in terms of the vectors $\psi_\mathbf{z}$ with $\mathbf{z}=(z,j_z)$ and $z=0,1,\ldots,r$. The choice is not unique, as the vectors are not independent. The result is independent of the choice (in other words, any choice can be done).

\subsection{Equivalence between Algorithm \ref{AlgoNonAbel} and \ref{AlgoBook}}

The following result proves that the codistribution \Obs, which is automatically computed by
Algorithm \ref{AlgoNonAbel}, is the observability codistribution.

\begin{theorem}\label{TheoremFundamentalNewNonAbel}
There exists $\widehat{k}$ such that, for any $k\ge\widehat{k}$, $\OBS\subseteq\Omega_k$. Conversely, for any $k$, $\Omega_{k}\subseteq\OBS$.
\end{theorem}

\proof{Let us prove the first statement.
From the $(k+1)^{th}$ step of Algorithm \ref{AlgoBook}, we obtain that $ \sum_{q=1}^{m_w}\sum_{i=1}^{m_u}\sum_{\alpha_1=0}^{m_w}\ldots\sum_{\alpha_k=0}^{m_w}
 \textnormal{span} \left \{\mathcal{L}_{[f^i]^{(\alpha_1,\ldots,\alpha_k)}} \nabla \tih_q \right \}\in\Omega_{k+1}$. 
From (\ref{EquationTOBSDef}), we immediately obtain that
\[
\tobs\subseteq\Omega_{\ND+\NO+1}.
\]

By comparing the recursive step of Algorithm \ref{AlgoNonAbel} with the one of Algorithm \ref{AlgoBook} it is immediate to conclude that, for any integer $v\ge0$:

\[
\OBS_v\subseteq\Omega_{\ND+\NO+1+v}
\]

and, as Algorithm \ref{AlgoNonAbel} converges in at most $n-m_w$ steps to $\OBS$, we conclude that:

\[
\OBS\subseteq\Omega_{\ND+\NO+1+n-m_w}
\]

and this proves the first statement with $\widehat{k}=\ND+\NO+1+n-m_w\le 
3n-2m_w+1$.

Let us prove the second statement. We must prove that, for any $k$, $\Omega_{k}\subseteq\OBS$.
We proceed by induction. It is true at $k=0$ as $\Omega_0\subseteq\OBS_0\subseteq\OBS$.

Let us assume that
\[
\Omega_p\subseteq\OBS
\]

We must prove that 

\[
\Omega_{p+1}\subseteq\OBS
\]

As $\OBS$ is invariant under $\Li_{\widehat{g}^\beta}$ and $\Li_{f^i}$, we only need to prove that


\[
\nabla\Li_{[f^i]^{(\alpha_1',\ldots,\alpha_p')}} \tih_q \in\OBS
\]

for any choice of $\alpha_1',\ldots,\alpha_p'=0,1,\ldots,m_w$, and every $q=1,\ldots,m_w$.
We distinguish the following two cases:

\begin{enumerate}

\item $p\le\ND+\NO$.

\item $p\ge \ND+\NO+1$.

\end{enumerate}

The first case is trivial as $\nabla\Li_{[f^i]^{(\alpha_1',\ldots,\alpha_p')}} \tih_q\in\tobs\subseteq\OBS$.

Let us consider the second case. We start by remarking that, when $p\ge r+s+1$ we have:

\[
[f^i]^{(\alpha_1',\ldots,\alpha_p')}=[\psi_{(r,j_r)}]^{(\alpha_{r+1}',\ldots,\alpha_{p}')}
\]

where $j_r$ is such that $[f^i]^{(\alpha_1',\ldots,\alpha_r')}=\psi_{(r,j_r)}$.

We set $p=\ND+\NO+l'$, with $l'\ge1$, and
we use Equation 
(\ref{EquationPropoNonAbel}), with $l=\NO+l'$, and $\alpha_1,\ldots,\alpha_l=\alpha_{r+1}',\ldots,\alpha_p'$. We have:

\begin{equation}\label{EquationProofFin}
\nabla\Li_{[f^i]^{(\alpha_1',\ldots,\alpha_p')}} \tih_q=
\nabla\Li_{[\psi_{(r,j_r)}]^{(\alpha_{r+1}',\ldots,\alpha_{p}')}} \tih_q = 
\end{equation}
\[
\sum_{\beta(l-1)} M_{\beta_1,\ldots,\beta_{l-1}}^{ \alpha_{r+2}',\ldots,\alpha_{p}'}\sum_\mathbf{z}\left(\nabla\C^{ \alpha_{r+1}'}_{(r,j_r)\mathbf{z}}\right)
\Li_{\psi_\mathbf{z}}\Li_{g^{\beta_{p-r-1}}}\ldots\Li_{g^{\beta_1}} \tih_q
\]
\[
~~\mod~\Omega_{p}+\Li_{\widehat{g}}\Omega_{p}
\]

where:
\begin{itemize}

\item The sum $\sum_{\beta(l-1)}$ is defined as follows:

\[
\sum_{\beta(l-1)}\triangleq\sum_{\beta_1=0}^{m_w}\sum_{\beta_2=0}^{m_w}\ldots\sum_{\beta_{l-1}=0}^{m_w},
\]

\item The sum $\sum_\mathbf{z}$ is defined as follows:

\[
\sum_\mathbf{z}=\sum_{(z,j_z)}\triangleq \sum_{z=0}^r\sum_{j_z=1}^{(m_w+1)^zm_u}.
\]

\item $M$ is a suitable multi-index object, which is non singular (i.e., it can be inverted with respect to all its indices).

\item $\Li_{\widehat{g}}\Omega_{p}$ stands for $\sum_{\gamma=0}^{m_w}\Li_{\widehat{g}^\gamma}\Omega_{p}$


\end{itemize}

The validity of the above result is proved in Section \ref{SectionPropositionFundamentalNonAbel} (see Proposition \ref{PropositionFundamentalNewNonAbel}).

Let us consider now the function $\Li_{g^{\beta_{p-r-1}}}\ldots\Li_{g^{\beta_1}} \tih_q$. Its differential belongs to $\Omega^g$.

We denote by $t$ the dimension of $\Omega^g$ and by $h^g_1,\ldots,h^g_t$ a basis of $\Omega^g$, i.e.:

\[
\Omega^g=\textnormal{span}\left\{
\nabla h^g_1,\ldots,\nabla h^g_t
\right\}
\]

In addition, we choose these generators ($h^g_1,\ldots,h^g_t$) directly among the functions built by Algorithm \ref{AlgoOmegagNonAbel}. 
%
%
As a result, 

\begin{equation}\label{EquationPerDopo1}
\nabla\Li_{g^{\beta_{p-r-1}}}\ldots\Li_{g^{\beta_1}} \tih_q=\sum_{j=1}^tc^p_j\nabla h^g_j
\end{equation}

with $c^p_1,\ldots,c^p_t$ suitable coefficients (scalar functions in the manifold $\M$). 
We have:

\begin{equation}\label{EquationPerDopo2}
\Li_{\psi_\mathbf{z}}\Li_{g^{\beta_{p-r-1}}}\ldots\Li_{g^{\beta_1}} \tih_q=
\end{equation}
\[
\nabla\Li_{g^{\beta_{p-r-1}}}\ldots\Li_{g^{\beta_1}} \tih_q \cdot \psi_\mathbf{z}=
\left(\sum_{j=1}^tc^p_j\nabla h^g_j\right)\cdot\psi_\mathbf{z}=
\sum_{j=1}^tc^p_j\Li_{\psi_\mathbf{z}}h^g_j
\]

By replacing this expression in (\ref{EquationProofFin}) we obtain:

\begin{equation}\label{EquationProofFin1}
\nabla\Li_{[f^i]^{(\alpha_1',\ldots,\alpha_p')}} \tih_q=
\end{equation}
\[
\sum_{j=1}^tc^p_j
\sum_{\beta(l-1)} M_{\beta_1,\ldots,\beta_{l-1}}^{ \alpha_{r+2}',\ldots,\alpha_{p}'}\sum_\mathbf{z}\left(\nabla\C^{\alpha_{r+1}'}_{(r,j_r)\mathbf{z}}\right)
\Li_{\psi_\mathbf{z}}h^g_j
~~\mod~\Omega_{p}+\Li_{\widehat{g}}\Omega_{p}
\]

Now we use again Equation 
(\ref{EquationPropoNonAbel}) with $l=s$. We have:

\[
\nabla\Li_{[\psi_\mathbf{r}]^{(\alpha_1,\ldots,\alpha_s)}} \tih_q = 
\]
\[
\sum_{\beta(s-1)} M_{\beta_1,\ldots,\beta_{s-1}}^{ \alpha_{2},\ldots,\alpha_s}\sum_\mathbf{z}\left(\nabla\C^{ \alpha_1}_\mathbf{rz}\right)\Li_{\psi_\mathbf{z}}\Li_{g^{\beta_{s-1}}}\ldots\Li_{g^{\beta_1}} \tih_q
\]
\[
\mod~\Omega_{r+s}+\Li_{\widehat{g}}\Omega_{r+s}
\]

We remark that

\[
\nabla\Li_{[\psi_\mathbf{r}]^{(\alpha_1,\ldots,\alpha_s)}} \tih_q \in\Omega_{r+s+1}
\]

for any choice of $\alpha_1,\ldots,\alpha_s$ and $q$.
In addition, because of the non singularity of $M$, we obtain that 

\[
\sum_\mathbf{z}\left(\nabla\C^{ \alpha_1}_\mathbf{rz}\right)\Li_{\psi_\mathbf{z}}\Li_{g^{\beta_{s-1}}}\ldots\Li_{g^{\beta_1}} \tih_q\in\Omega_{r+s+1}
\]

for any choice of $\beta_1,\ldots,\beta_{s-1}$, $q$ and $\alpha_1$. On the other hand, the generators of $\Omega^g$ ($\nabla h^g_j$) are
built by Algorithm \ref{AlgoOmegagNonAbel},and they are the differentials of the functions that belong to the function space that includes the outputs ($\tih_q$) together with their Lie derivatives along $g^\beta$ of order that does not exceed $s-1$. 
Therefore, by setting $\alpha_1=\alpha_{r+1}'$, we obtain:

\[
\sum_\mathbf{z}\left(\nabla\C^{\alpha_{r+1}'}_\mathbf{rz}\right)\Li_{\psi_\mathbf{z}}h^g_j\in\Omega_{r+s+1}, ~~j=1,\ldots,t.
\]

By using this in (\ref{EquationProofFin1}) and by knowing that $p\ge r+s+1$, we immediately obtain that

\[
\nabla\Li_{[f^i]^{(\alpha_1',\ldots,\alpha_p')}} \tih_q\in\Omega_p+\Li_{\widehat{g}}\Omega_{p}
\]

From the inductive assumption, we know that $\Omega_p\subseteq\OBS$.
In addition, as $\OBS$ is invariant under $\Li_{\widehat{g}^0},\Li_{\widehat{g}^1},\ldots,\Li_{\widehat{g}^{m_w}}$, 
we have $\Omega_p+\Li_{\widehat{g}}\Omega_{p}\subseteq\OBS$ and, consequently, $\nabla\Li_{[f^i]^{(\alpha_1',\ldots,\alpha_p')}} \tih_q\in\OBS$.
$\blacktriangleleft$}

\newpage

\subsection{Intermediate technical results}\label{SectionPropositionFundamentalNonAbel}

The goal of this section is to prove the validity of Equation (\ref{EquationPropoNonAbel}), which has been used twice in the proof of Theorem \ref{TheoremFundamentalNewNonAbel}.
To achieve this goal, we basically need to proceed along two directions:

\begin{itemize}

\item Introduce the quantities $\C_{\mathbf{kz}}^{\alpha_1,\ldots,\alpha_l}$, according to Equation (\ref{EquationCNonAbelianGen}), and obtaining the recursive law in (\ref{EquationLawNonAbel}) and (\ref{EquationLawDiffNonAbel}).

\item First proving the validity of Equation (\ref{EquationLemmaNonAbel}), which is more general than Equation (\ref{EquationPropoNonAbel}). In particular, this latter is obtained from the former in a special setting.

\end{itemize}


We introduce the quantities $\C_{\mathbf{kz}}^{\alpha_1,\ldots,\alpha_l}$ by extending the quantities $\C_{\mathbf{kz}}^\alpha$ in (\ref{EquationDefinitionC}) to the case when the autobracket is applied repetitively multiple times. $\C_{\mathbf{kz}}^{\alpha_1,\ldots,\alpha_l}$ are implicitly defined by the following equation:

\begin{equation}\label{EquationCNonAbelianGen}
[\psi_\mathbf{k}]^{(\alpha_1,\ldots,\alpha_l)}=\sum_\mathbf{z}\C_{\mathbf{kz}}^{\alpha_1,\ldots,\alpha_l}\psi_\mathbf{z}
\end{equation}
where $\sum_\mathbf{z}=\sum_{(z,j_z)}$ stands for $\sum_{z=0}^r\sum_{j_z=1}^{(m_w+1)^zm_u}$, and the operation at the left side,
$[\cdot]^{(\alpha_1,\ldots,\alpha_l)}$, is the multiple autobracket defined after Definition \ref{DefinitionAutobracket}.
As in the case of (\ref{EquationDefinitionC}), the term at the left hand side of the above equation belongs to $\Delta$ and can be expressed in terms of the set of vectors $\psi$ (again, we remind the reader that $\Delta$ is invariant under the autobracket operation and the vectors $\psi$ generate the entire distribution).
Equation (\ref{EquationCNonAbelianGen}) does not provide uniquely the quantities $\C_\mathbf{kz}^{\alpha_1,\ldots,\alpha_l}$. This because the vectors $\psi$ are not necessarily independent. On the other hand, it is possible to set them in such a way that they satisfy the following recursive law. The law is stated by the following Lemma:

\begin{lm}\label{LemmaLawNonAbel}
It is possible to set the quantities $\C_\mathbf{kz}^{\alpha_1,\ldots,\alpha_l}$ such that they satisfy the following recursive law:
\begin{equation}\label{EquationLawNonAbel}
\C^{\alpha_1,\ldots,\alpha_l,\alpha_{l+1}}_\mathbf{kz}=\Li_{\widehat{g}^{\alpha_{l+1}}}\C^{\alpha_1,\ldots,\alpha_l}_\mathbf{kz}+\sum_\mathbf{z'}\C^{\alpha_1,\ldots,\alpha_l}_\mathbf{kz'}\C^{\alpha_{l+1}}_\mathbf{z'z}
\end{equation}
and its differential form:
\begin{equation}\label{EquationLawDiffNonAbel}
\nabla\C^{\alpha_1,\ldots,\alpha_l,\alpha_{l+1}}_\mathbf{kz}=
\end{equation}
\[
\Li_{\widehat{g}^{\alpha_{l+1}}}\nabla\C^{\alpha_1,\ldots,\alpha_l}_\mathbf{kz}+
\sum_\mathbf{z'}\nabla\C^{\alpha_1,\ldots,\alpha_l}_\mathbf{kz'}\C^{\alpha_{l+1}}_\mathbf{z'z}+\sum_\mathbf{z'}\C^{\alpha_1,\ldots,\alpha_l}_\mathbf{kz'}\nabla\C^{\alpha_{l+1}}_\mathbf{z'z}
\]
where $\sum_\mathbf{z'}=\sum_{(z',j_{z'})}$ stands for $\sum_{z'=0}^r\sum_{j_{z'}=1}^{(m_w+1)^{z'}m_u}$.

\end{lm}

\proof{This proof needs the following equality, which is a property of the autobracket. Given a vector field $\phi$ and a scalar field $a$, we have:

\begin{equation}\label{EquationAutobracketProperty}
[a\phi]^\alpha=\left(\Li_{\widehat{g}^\alpha}a\right)~\phi + a[\phi]^\alpha
\end{equation}

for any $\alpha=0,1,\ldots,m_w$.
This is obtained by a direct computation. We have:
\[
[a\phi]^\alpha=\sum_{\beta=0}^{m_w}\nu^\alpha_\beta[g^\beta,~a\phi]=\sum_{\beta=0}^{m_w}\nu^\alpha_\beta a[g^\beta,~\phi]+\sum_{\beta=0}^{m_w}\left(\nu^\alpha_\beta\Li_{g^\beta}a\right)\phi=
\]
\[
 a[\phi]^\alpha+\left(\Li_{\widehat{g}^\alpha}a\right)~\phi.
\]

Now, let us prove the validity of (\ref{EquationLawNonAbel}).
By definition we have:

\begin{equation}\label{EquationLemmaNonAbelStep1}
[\psi_\mathbf{k}]^{(\alpha_1,\ldots,\alpha_l,\alpha_{l+1})}= \sum_\mathbf{z}\C_\mathbf{kz}^{\alpha_1,\ldots,\alpha_l,\alpha_{l+1}}\psi_\mathbf{z}
\end{equation}

On the other hand,

\[
[\psi_\mathbf{k}]^{(\alpha_1,\ldots,\alpha_l,\alpha_{l+1})}= 
\left[
[\psi_\mathbf{k}]^{(\alpha_1,\ldots,\alpha_l)}
\right]^{\alpha_{l+1}}=
\sum_\mathbf{z}
\left[
\C^{\alpha_1,\ldots,\alpha_l}_\mathbf{kz}\psi_\mathbf{z}
\right]^{\alpha_{l+1}}
\]

and, by using (\ref{EquationAutobracketProperty}), we obtain
\[
=\sum_\mathbf{z}
\left(\Li_{\widehat{g}^{\alpha_{l+1}}}\C^{\alpha_1,\ldots,\alpha_l}_\mathbf{kz}
\right)\psi_\mathbf{z}+
\sum_\mathbf{z}
\C^{\alpha_1,\ldots,\alpha_l}_\mathbf{kz}[\psi_\mathbf{z}]^{\alpha_{l+1}}=
\]
\[
\sum_\mathbf{z}
\left(\Li_{\widehat{g}^{\alpha_{l+1}}}\C^{\alpha_1,\ldots,\alpha_l}_\mathbf{kz}
\right)\psi_\mathbf{z}+
\sum_\mathbf{zz'}
\C^{\alpha_1,\ldots,\alpha_l}_\mathbf{kz}\C^{\alpha_{l+1}}_\mathbf{zz'}\psi_\mathbf{z'}=
\]
\[
\sum_\mathbf{z}
\left(
\Li_{\widehat{g}^{\alpha_{l+1}}}\C^{\alpha_1,\ldots,\alpha_l}_\mathbf{kz}+
\sum_\mathbf{z'}
\C^{\alpha_1,\ldots,\alpha_l}_\mathbf{kz'}\C^{\alpha_{l+1}}_\mathbf{z'z}
\right)\psi_\mathbf{z}
\]

By comparing with (\ref{EquationLemmaNonAbelStep1})
we obtain (\ref{EquationLawNonAbel}) and, by differentiating, (\ref{EquationLawDiffNonAbel}).
$\blacktriangleleft$}

We have the following fundamental result:

\begin{lm}\label{LemmaFundamentalNewNonAbel}
Given $\psi_\mathbf{k}$ ($\mathbf{k}=(k,j_k)$), for any integer $l\ge1$, and any set of integers $\alpha_1,\ldots,\alpha_l$ that take the values $0,1,\ldots,m_w$, for any integer $0\le j\le l-1$, and any $q=1,\ldots,m_w$, we have:
\begin{equation}\label{EquationLemmaNonAbel}
\nabla\Li_{[\psi_\mathbf{k}]^{(\alpha_1,\ldots,\alpha_l)}} \tih_q = 
\end{equation}
\[
\sum_{\beta_1=0}^{m_w}\sum_{\beta_2=0}^{m_w}\ldots\sum_{\beta_j=0}^{m_w} M_{\beta_1,\ldots,\beta_j}^{ \alpha_{l-j+1},\ldots,\alpha_l}\sum_\mathbf{z}\left(\nabla\C^{ \alpha_1,\ldots,\alpha_{l-j}}_\mathbf{kz}\right)\Li_{\psi_\mathbf{z}}\Li_{g^{\beta_j}}\ldots\Li_{g^{\beta_1}} \tih_q
\]
\[
\mod~\Omega_{r+l}+\Li_{\widehat{g}}\Omega_{r+l}
\]
where $M$ is a suitable multi-index object, which is non singular (i.e., it can be inverted with respect to all its indices).
\end{lm}

\proof{For notation simplicity we omit the index $q$, i.e., we denote by $h$ the function $\tih_q$. In addition, we use the notation:

\[
\sum_{\beta(j)}\triangleq \sum_{\beta_1=0}^{m_w}\sum_{\beta_2=0}^{m_w}\ldots\sum_{\beta_j=0}^{m_w}.
\]

We proceed by induction on $l$. Let us set $l=1$. We only have $j=0$. We have:

\[
\nabla\Li_{[\psi_\mathbf{k}]^{\alpha_1}} h=
\sum_\mathbf{z}
\nabla\Li_{\C_\mathbf{kz}^{\alpha_1}\psi_\mathbf{z}} h=
\sum_\mathbf{z}\nabla\left(
\C_\mathbf{kz}^{\alpha_1}\Li_{\psi_\mathbf{z}}h
\right)=
\]
\[
\sum_\mathbf{z}
\C_\mathbf{kz}^{\alpha_1}\nabla\Li_{\psi_\mathbf{z}} h+
\sum_\mathbf{z}
\left(\nabla\C_\mathbf{kz}^{\alpha_1}\right)\Li_{\psi_\mathbf{z}}h
\]

On the other hand,

\[
\nabla\Li_{\psi_\mathbf{z}} h\in\Omega_{r+1},
\]

Hence,

\[
\nabla\Li_{[\psi_\mathbf{k}]^{\alpha_1}} h=
\sum_\mathbf{z}
\left(\nabla\C_\mathbf{kz}^{\alpha_1}\right)\Li_{\psi_\mathbf{z}}h,~~~~~\mod~ \Omega_{r+1},
\]

which proves (\ref{EquationLemmaNonAbel}) when $l=1$, and $j=0$.

Let us consider the recursive step. We assume that (\ref{EquationLemmaNonAbel}) holds at a given $l=l^*>1$. We have:

\begin{equation}\label{EquationLemmaNonAbel*}
\nabla\Li_{[\psi_\mathbf{k}]^{(\alpha_1,\ldots,\alpha_{l^*})}} h = 
\end{equation}
\[
\sum_{\beta(j)} M_{\beta_1,\ldots,\beta_j}^{ \alpha_{l^*-j+1},\ldots,\alpha_{l^*}}\sum_\mathbf{z}\left(\nabla\C^{ \alpha_1,\ldots,\alpha_{l^*-j}}_\mathbf{kz}\right)\Li_{\psi_\mathbf{z}}\Li_{g^{\beta_j}}\ldots\Li_{g^{\beta_1}}h
\]
\[
\mod~\Omega_{r+l^*}+\Li_{\widehat{g}}\Omega_{r+l^*}
\]

for any $j=0,1,\ldots,l^*-1$ and with $M$ non singular. In addition, (\ref{EquationLemmaNonAbel}) holds at any $l\le l^*$ for any $j\le l-1$. In particular, it holds for any 
$l\le l^*$ and $j=l-1$. Hence we also have:

\begin{equation}\label{EquationLemmaNonAbellm1}
\nabla\Li_{[\psi_\mathbf{k}]^{(\alpha_1,\ldots,\alpha_l)}} h = 
\end{equation}
\[
\sum_{\beta(l-1)} M_{\beta_1,\ldots,\beta_{l-1}}^{ \alpha_{2},\ldots,\alpha_l}\sum_\mathbf{z}\left(\nabla\C^{ \alpha_1}_\mathbf{kz}\right)\Li_{\psi_\mathbf{z}}\Li_{g^{\beta_{l-1}}}\ldots\Li_{g^{\beta_1}} h
\]
\[
\mod~\Omega_{r+l}+\Li_{\widehat{g}}\Omega_{r+l}
\]

for any $l\le l^*$ with $M$ non singular.

We must prove the validity of (\ref{EquationLemmaNonAbel}) at $l^*+1$ and for any $j=0,\ldots,l^*$, i.e.:

\begin{equation}\label{EquationLemmaNonAbel*p1}
\nabla\Li_{[\psi_\mathbf{k}]^{(\alpha_1,\ldots,\alpha_{l^*},\alpha_{l^*+1})}} h = 
\end{equation}
\[
\sum_{\beta(j)} M_{\beta_1,\ldots,\beta_j}^{ \alpha_{l^*-j+2},\ldots,\alpha_{l^*+1}}\sum_\mathbf{z}\left(\nabla\C^{ \alpha_1,\ldots,\alpha_{l^*+1-j}}_\mathbf{kz}\right)\Li_{\psi_\mathbf{z}}\Li_{g^{\beta_j}}\ldots\Li_{g^{\beta_1}}h 
\]
\[
\mod~\Omega_{r+l^*+1}+\Li_{\widehat{g}}\Omega_{r+l^*+1}
\]

for a suitable non singular $M$.

We proceed by induction on $j$. Let us consider $j=0$. We have:

\[
\nabla\Li_{[\psi_\mathbf{k}]^{(\alpha_1,\ldots,\alpha_{l^*},\alpha_{l^*+1})}} h = 
\nabla \Li_{\sum_\mathbf{z}\C_\mathbf{kz}^{ \alpha_1,\ldots,\alpha_{l^*+1}}\psi_\mathbf{z}} h=
\]
\[
\sum_\mathbf{z}
\nabla\left(
\C_\mathbf{kz}^{ \alpha_1,\ldots,\alpha_{l^*+1}}\Li_{\psi_\mathbf{z}}h
\right)=
\]
\[
\sum_\mathbf{z}
\C_\mathbf{kz}^{ \alpha_1,\ldots,\alpha_{l^*+1}}\nabla\Li_{\psi_\mathbf{z}} h+
\sum_\mathbf{z}
\left(\nabla\C_\mathbf{kz}^{ \alpha_1,\ldots,\alpha_{l^*+1}}\right)\Li_{\psi_\mathbf{z}}h
\]

On the other hand,

\[
\nabla\Li_{\psi_\mathbf{z}} h\in\Omega_{r+1}\subseteq \Omega_{r+l^*+1}+\Li_{\widehat{g}}\Omega_{r+l^*+1},
\]

Hence,

\[
\nabla\Li_{[\psi_\mathbf{k}]^{(\alpha_1,\ldots,\alpha_{l^*},\alpha_{l^*+1})}} h = 
\]
\[
\sum_\mathbf{z}\left(\nabla\C_\mathbf{kz}^{ \alpha_1,\ldots,\alpha_{l^*+1}}\right)\Li_{\psi_\mathbf{z}}h
~~\mod~\Omega_{r+l^*+1}+\Li_{\widehat{g}}\Omega_{r+l^*+1}
\]

which proves (\ref{EquationLemmaNonAbel}) when $l=l^*+1$, and $j=0$.

Let us assume that (\ref{EquationLemmaNonAbel*p1}) holds at a given $j=j^*\le l^*-1$. 
\[
\nabla\Li_{[\psi_\mathbf{k}]^{(\alpha_1,\ldots,\alpha_{l^*},\alpha_{l^*+1})}} h = 
\]
\[
\sum_{\beta(j^*)} M_{\beta_1,\ldots,\beta_{j^*}}^{ \alpha_{l^*-j^*+2},\ldots,\alpha_{l^*+1}}\sum_\mathbf{z}\left(\nabla\C^{ \alpha_1,\ldots,\alpha_{l^*+1-j^*}}_\mathbf{kz}\right)\Li_{\psi_\mathbf{z}}\Li_{g^{\beta_{j^*}}}\ldots\Li_{g^{\beta_1}}h
 \]
\[
\mod~\Omega_{r+l^*+1}+\Li_{\widehat{g}}\Omega_{r+l^*+1}
\]

We must prove:

\[
\nabla\Li_{[\psi_\mathbf{k}]^{(\alpha_1,\ldots,\alpha_{l^*},\alpha_{l^*+1})}} h = 
\]
\[
\sum_{\beta(j^*+1)} M_{\beta_1,\ldots,\beta_{j^*},\beta_{j^*+1}}^{ \alpha_{l^*-j^*+1},\ldots,\alpha_{l^*+1}}\sum_\mathbf{z}\left(\nabla\C^{ \alpha_1,\ldots,\alpha_{l^*-j^*}}_\mathbf{kz}\right)\Li_{\psi_\mathbf{z}}\Li_{g^{\beta_{j^*+1}}}\Li_{g^{\beta_{j^*}}}\ldots\Li_{g^{\beta_1}}h
\]
\[
\mod~\Omega_{r+l^*+1}+\Li_{\widehat{g}}\Omega_{r+l^*+1}
\]

We adopt the notation:

\[
\lambda\triangleq\Li_{g^{\beta_{j^*}}}\ldots\Li_{g^{\beta_1}}h
\]
\[
\C^-=\C^{ \alpha_1,\ldots,\alpha_{l^*-j^*}}
\]
\[
\C^{+\alpha_{l^*-j^*+1}}=\C^{ \alpha_1,\ldots,\alpha_{l^*-j^*},\alpha_{l^*-j^*+1}}
\]

We have:

\begin{equation}\label{Equation37}
\nabla\Li_{[\psi_\mathbf{k}]^{(\alpha_1,\ldots,\alpha_{l^*},\alpha_{l^*+1})}} h = 
\end{equation}
\[
\sum_{\beta(j^*)} M_{\beta_1,\ldots,\beta_{j^*}}^{ \alpha_{l^*-j^*+2},\ldots,\alpha_{l^*+1}}\sum_\mathbf{z}\left(\nabla\C^{+\alpha_{l^*-j^*+1}}_\mathbf{kz}\right)\Li_{\psi_\mathbf{z}}\lambda
\]
\[
\mod~\Omega_{r+l^*+1}+\Li_{\widehat{g}}\Omega_{r+l^*+1}
\]

We must prove:

\begin{equation}\label{Equation38}
\nabla\Li_{[\psi_\mathbf{k}]^{(\alpha_1,\ldots,\alpha_{l^*},\alpha_{l^*+1})}} h = 
\end{equation}
\[
\sum_{\beta(j^*+1)} M_{\beta_1,\ldots,\beta_{j^*},\beta_{j^*+1}}^{ \alpha_{l^*-j^*+1},\ldots,\alpha_{l^*+1}}\sum_\mathbf{z}\left(\nabla\C^{-}_\mathbf{kz}\right)\Li_{\psi_\mathbf{z}}\Li_{g^{\beta_{j^*+1}}}\lambda 
\]
\[
\mod~\Omega_{r+l^*+1}+\Li_{\widehat{g}}\Omega_{r+l^*+1}
\]

From (\ref{EquationLawDiffNonAbel}) we have:

\[
\nabla\C^{+\alpha_{l^*-j^*+1}}_\mathbf{kz}=\Li_{\widehat{g}^{\alpha_{l^*-j^*+1}}}\nabla\C^{-}_\mathbf{kz}+
\sum_\mathbf{z'}\nabla\C^{-}_\mathbf{kz'}\C^{\alpha_{l^*-j^*+1}}_\mathbf{z'z}+
\]
\[
\sum_\mathbf{z'}\C^{-}_\mathbf{kz'}\nabla\C^{\alpha_{l^*-j^*+1}}_\mathbf{z'z}
\]

We substitute in (\ref{Equation37}) and we obtain:

\[
\nabla\Li_{[\psi_\mathbf{k}]^{(\alpha_1,\ldots,\alpha_{l^*},\alpha_{l^*+1})}} h = 
\]
\[
\sum_{\beta(j^*)} M_{\beta_1,\ldots,\beta_{j^*}}^{ \alpha_{l^*-j^*+2},\ldots,\alpha_{l^*+1}}
\sum_\mathbf{z}\left(
\Li_{\widehat{g}^{\alpha_{l^*-j^*+1}}}\nabla\C^{-}_\mathbf{kz}+
\sum_\mathbf{z'}\nabla\C^{-}_\mathbf{kz'}\C^{\alpha_{l^*-j^*+1}}_\mathbf{z'z}+
\right.
\]
\[
\left.
\sum_\mathbf{z'}\C^{-}_\mathbf{kz'}\nabla\C^{\alpha_{l^*-j^*+1}}_\mathbf{z'z}
\right)
\Li_{\psi_\mathbf{z}}\lambda
~~\mod~\Omega_{r+l^*+1}+\Li_{\widehat{g}}\Omega_{r+l^*+1}
\]

Let us consider the first term on the right hand side.

\[
\sum_{\beta(j^*)} M_{\beta_1,\ldots,\beta_{j^*}}^{ \alpha_{l^*-j^*+2},\ldots,\alpha_{l^*+1}}
\sum_\mathbf{z}\left(
\Li_{\widehat{g}^{\alpha_{l^*-j^*+1}}}\nabla\C^{-}_\mathbf{kz}\right)\Li_{\psi_\mathbf{z}}\lambda=
\]
\[
\sum_{\beta(j^*)} M_{\beta_1,\ldots,\beta_{j^*}}^{ \alpha_{l^*-j^*+2},\ldots,\alpha_{l^*+1}}
\left[
\sum_\mathbf{z}
\Li_{\widehat{g}^{\alpha_{l^*-j^*+1}}}\left(
\left(
\nabla\C^{-}_\mathbf{kz}\right)\Li_{\psi_\mathbf{z}}\lambda
\right)-
\right.
\]
\[
\left.
\sum_\mathbf{z}
\left(
\nabla\C^{-}_\mathbf{kz}\right)\Li_{\widehat{g}^{\alpha_{l^*-j^*+1}}}\Li_{\psi_\mathbf{z}}\lambda
\right]
\]

Note that $j^*\le l^*-1$. Hence, we are allowed to use (\ref{EquationLemmaNonAbel*}) at $j=j^*$.
In the new notation, it tells us that

\[
\nabla\Li_{[\psi_\mathbf{k}]^{(\alpha_1,\ldots,\alpha_{l^*})}} h = \sum_{\beta(j^*)} M_{\beta_1,\ldots,\beta_{j^*}}^{ \alpha_{l^*-j^*+1},\ldots,\alpha_{l^*}}
\sum_\mathbf{z}\left(\nabla\C^{-}_\mathbf{kz}\right)\Li_{\psi_\mathbf{z}}\lambda
\]
\[
\mod~\Omega_{r+l^*}+\Li_{\widehat{g}}\Omega_{r+l^*}
\]

On the other hand,

\[
\nabla\Li_{[\psi_\mathbf{k}]^{(\alpha_1,\ldots,\alpha_{l^*})}} h   \in\Omega_{k+l^*+1}
\subseteq\Omega_{r+l^*+1}
\]

and, from the invertibility of $M$, it follows that

\[
\sum_\mathbf{z}\left(\nabla\C^{-}_\mathbf{kz}\right)\Li_{\psi_\mathbf{z}}\lambda\in\Omega_{r+l^*+1}
\]

As a result,

\[
\sum_{\beta(j^*)} M_{\beta_1,\ldots,\beta_{j^*}}^{ \alpha_{l^*-j^*+2},\ldots,\alpha_{l^*+1}}
\sum_\mathbf{z}
\Li_{\widehat{g}^{\alpha_{l^*-j^*+1}}}\left(
\left(
\nabla\C^{-}_\mathbf{kz}\right)\Li_{\psi_\mathbf{z}}\lambda
\right)
\in\Li_{\widehat{g}}\Omega_{r+l^*+1}
\]
 
and we have:

\[
\sum_{\beta(j^*)} M_{\beta_1,\ldots,\beta_{j^*}}^{ \alpha_{l^*-j^*+2},\ldots,\alpha_{l^*+1}}
\sum_\mathbf{z}
\left(
\Li_{\widehat{g}^{\alpha_{l^*-j^*+1}}}\nabla\C^{-}_\mathbf{kz}+
\sum_\mathbf{z'}\nabla\C^{-}_\mathbf{kz'}\C^{\alpha_{l^*-j^*+1}}_\mathbf{z'z}+
\right.
\]
\[
\left.
\sum_\mathbf{z'}\C^{-}_\mathbf{kz'}\nabla\C^{\alpha_{l^*-j^*+1}}_\mathbf{z'z}
\right)\Li_{\psi_\mathbf{z}}\lambda=
\]

\[
\sum_{\beta(j^*)} M_{\beta_1,\ldots,\beta_{j^*}}^{ \alpha_{l^*-j^*+2},\ldots,\alpha_{l^*+1}}
\sum_\mathbf{z}
\left(
-\nabla\C^{-}_\mathbf{kz}\Li_{\widehat{g}^{\alpha_{l^*-j^*+1}}}\Li_{\psi_\mathbf{z}}\lambda+
\right.
\]
\[
\left.
\sum_\mathbf{z'}\nabla\C^{-}_\mathbf{kz'}\C^{\alpha_{l^*-j^*+1}}_\mathbf{z'z}\Li_{\psi_\mathbf{z}}\lambda+
\sum_\mathbf{z'}\C^{-}_\mathbf{kz'}\nabla\C^{\alpha_{l^*-j^*+1}}_\mathbf{z'z}\Li_{\psi_\mathbf{z}}\lambda
\right)
\]

Let us consider the last term of the above. We consider (\ref{EquationLemmaNonAbellm1}). It holds for any $l\le l^*$. As $j^*\le l^*-1$, we are allowed to set $l=j^*+1$. We have:

\[
\nabla\Li_{[\psi_\mathbf{k}]^{(\alpha_1,\ldots,\alpha_{j^*+1})}} h = \sum_{\beta(j^*)} M_{\beta_1,\ldots,\beta_{j^*}}^{ \alpha_{2},\ldots,\alpha_{j^*+1}}
\sum_\mathbf{z}
\left(\nabla\C^{ \alpha_1}_\mathbf{kz}\right)\Li_{\psi_\mathbf{z}}\lambda
\]
\[
\mod~\Omega_{r+j^*+1}+\Li_{\widehat{g}}\Omega_{r+j^*+1}
\]

Hence, from the invertibility of $M$, it follows that

\[
\sum_\mathbf{z}\left(\nabla\C^{ \alpha_1}_\mathbf{kz}\right)\Li_{\psi_\mathbf{z}}\lambda\in\Omega_{r+j^*+2}
\]

As $j^*\le l^*-1$, this proves that:

\[
\sum_\mathbf{z}\left(\nabla\C^{ \alpha_1}_\mathbf{kz}\right)\Li_{\psi_\mathbf{z}}\lambda\in\Omega_{r+l^*+1} ~~\forall\alpha_1
\]

We have:

\[
\sum_{\beta(j^*)} M_{\beta_1,\ldots,\beta_{j^*}}^{ \alpha_{l^*-j^*+2},\ldots,\alpha_{l^*+1}}
\sum_\mathbf{z}\left(
-\nabla\C^{-}_\mathbf{kz}\Li_{\widehat{g}^{\alpha_{l^*-j^*+1}}}\Li_{\psi_\mathbf{z}}\lambda
\right.
\]
\[
\left.
+\sum_\mathbf{z'}\nabla\C^{-}_\mathbf{kz'}\C^{\alpha_{l^*-j^*+1}}_\mathbf{z'z}\Li_{\psi_\mathbf{z}}\lambda
+\sum_\mathbf{z'}\C^{-}_\mathbf{kz'}\nabla\C^{\alpha_{l^*-j^*+1}}_\mathbf{z'z}\Li_{\psi_\mathbf{z}}\lambda
\right)=
\]

\[
\sum_{\beta(j^*)} M_{\beta_1,\ldots,\beta_{j^*}}^{ \alpha_{l^*-j^*+2},\ldots,\alpha_{l^*+1}}
\sum_\mathbf{z}\left(
-\nabla\C^{-}_\mathbf{kz}\Li_{\widehat{g}^{\alpha_{l^*-j^*+1}}}\Li_{\psi_\mathbf{z}}\lambda
\right.
\]
\[
\left.
+\sum_\mathbf{z'}\nabla\C^{-}_\mathbf{kz'}\C^{\alpha_{l^*-j^*+1}}_\mathbf{z'z}\Li_{\psi_\mathbf{z}}\lambda
\right)~~~\mod~\Omega_{r+l^*+1}
\]

We have:

\[
-\sum_\mathbf{z}\nabla\C^{-}_\mathbf{kz}\Li_{\widehat{g}^{\alpha_{l^*-j^*+1}}}\Li_{\psi_\mathbf{z}}\lambda
+\sum_\mathbf{zz'}\nabla\C^{-}_\mathbf{kz}\C^{\alpha_{l^*-j^*+1}}_\mathbf{zz'}\Li_{\psi_\mathbf{z'}}\lambda=
\]

\[
\sum_\mathbf{z}\nabla\C^{-}_\mathbf{kz}
\left(
-\Li_{\widehat{g}^{\alpha_{l^*-j^*+1}}}\Li_{\psi_\mathbf{z}}\lambda
+\sum_\mathbf{z'}\C^{\alpha_{l^*-j^*+1}}_\mathbf{zz'}\Li_{\psi_\mathbf{z'}}\lambda
\right)=
\]
\[
\sum_\mathbf{z}\sum_{\gamma=0}^{m_w}\nabla\C^{-}_\mathbf{kz}\nu^{\alpha_{l^*-j^*+1}}_{\gamma}
\left(
-\Li_{g^{\gamma}}\Li_{\psi_\mathbf{z}}\lambda
+\Li_{\left[g^{\gamma},~\psi_\mathbf{z}\right]}\lambda
\right)=
\]
\[
\sum_{\gamma=0}^{m_w}\sum_\mathbf{z}\nabla\C^{-}_\mathbf{kz}\nu^{\alpha_{l^*-j^*+1}}_{\gamma}
\left(
-\Li_{\psi_\mathbf{z}}\Li_{g^{\gamma}}\lambda
\right)
\]

Hence, (\ref{Equation37}) becomes:

\[
\nabla\Li_{[\psi_\mathbf{k}]^{(\alpha_1,\ldots,\alpha_{l^*},\alpha_{l^*+1})}} h = 
\]
\[
\sum_{\beta(j^*)} M_{\beta_1,\ldots,\beta_{j^*}}^{ \alpha_{l^*-j^*+2},\ldots,\alpha_{l^*+1}}
\sum_{\gamma=0}^{m_w}\sum_\mathbf{z}(-\nu^{\alpha_{l^*-j^*+1}}_{\gamma})
\nabla\C^{-}_\mathbf{kz}
\left(
\Li_{\psi_\mathbf{z}}\Li_{g^{\gamma}}\lambda
\right)=
\]
\[
\sum_{\beta(j^*)} \sum_{\gamma=0}^{m_w}M_{\beta_1,\ldots,\beta_{j^*}}^{ \alpha_{l^*-j^*+2},\ldots,\alpha_{l^*+1}}
(-\nu^{\alpha_{l^*-j^*+1}}_{\gamma})
\sum_\mathbf{z}
\nabla\C^{-}_\mathbf{kz}
\left(
\Li_{\psi_\mathbf{z}}\Li_{g^{\gamma}}\lambda
\right)
\]
\[
\mod~\Omega_{r+j^*+1}+\Li_{\widehat{g}}\Omega_{r+j^*+1}
\]

which coincides with (\ref{Equation38}) with the dummy index $\gamma$ equal to the dummy index $\beta_{j^*+1}$ in (\ref{Equation38}) and

\[
M_{\beta_1,\ldots,\beta_{j^*},\beta_{j^*+1}}^{ \alpha_{l^*-j^*+1},\ldots,\alpha_{l^*+1}}=
M_{\beta_1,\ldots,\beta_{j^*}}^{ \alpha_{l^*-j^*+2},\ldots,\alpha_{l^*+1}}(-\nu^{\alpha_{l^*-j^*+1}}_{\beta_{j^*+1}}),
\]
which remains non singular as $\nu$ is non singular.
$\blacktriangleleft$}

We have the following fundamental result:

\begin{pr}\label{PropositionFundamentalNewNonAbel}
Given $\psi_\mathbf{r}$ ($\mathbf{r}=(r,j_r)$), for any integer $l\ge1$, and any set of integers $\alpha_1,\ldots,\alpha_l$ that take the values $0,1,\ldots,m_w$, and any $q=1,\ldots,m_w$, we have:

\begin{equation}\label{EquationPropoNonAbel}
\nabla\Li_{[\psi_\mathbf{r}]^{(\alpha_1,\ldots,\alpha_l)}} \tih_q = 
\end{equation}
\[
\sum_{\beta(l-1)} M_{\beta_1,\ldots,\beta_{l-1}}^{ \alpha_{2},\ldots,\alpha_l}\sum_\mathbf{z}\left(\nabla\C^{ \alpha_1}_\mathbf{rz}\right)\Li_{\psi_\mathbf{z}}\Li_{g^{\beta_{l-1}}}\ldots\Li_{g^{\beta_1}} \tih_q
\]
\[
\mod~\Omega_{r+l}+\Li_{\widehat{g}}\Omega_{r+l}
\]
with:

\begin{itemize}

\item $M$ a suitable multi-index object, which is non singular (i.e., it can be inverted with respect to all its indices).

\item 
$\sum_{\beta(l-1)}\triangleq\sum_{\beta_1=0}^{m_w}\sum_{\beta_2=0}^{m_w}\ldots\sum_{\beta_{l-1}=0}^{m_w}$.

\item $
\sum_\mathbf{z}=\sum_{(z,j_z)}\triangleq \sum_{z=0}^r\sum_{j_z=1}^{(m_w+1)^zm_u}$.
\end{itemize}
\end{pr}

\proof{The above equality is obtained from Lemma \ref{LemmaFundamentalNewNonAbel}, with $j=l-1$ and $\mathbf{k}=\mathbf{r}$.
$\blacktriangleleft$}


\chapter{Solution in the general non canonical case}\label{ChapterSolutionNonCanonic}

Chapter \ref{ChapterSolutionCanonic} provides
the solution of the UIO problem. On the other hand, to run the algorithm that provides the observability codistribution (Algorithm \ref{AlgoNonAbel}) the system must be in canonical form. In other words, in accordance with Definition \ref{DefinitionCanonicalForm}, the output must include $m_w$ functions, $\widetilde{h}_1,\ldots,\widetilde{h}_{m_w}$, such that the tensor $\mu$ defined in (\ref{EquationTensorMSynchro}) is non singular (or, equivalently, its rank is equal to $m_w$).
But this is a special case. We want to deal with the case where this is not possible (e.g., when $m_w>p$).
In addition, we also want to deal with systems that are not canonic with respect to their unknown inputs. As we will see, in some cases, we can transform a non canonic system into a canonic system. We will call these systems {\it canonizable}.
We also want to deal with the case where the system is not canonic and not even canonizable.
In few words, we will deal with any system.

This chapter provides the automatic procedure that, in a finite number of steps, provides the observability codistribution in any case.
This chapter only provides a very concise summary of this solution. The interested reader can find its theoretical foundation in \cite{IF22}. 
The solution is here presented in the form of a new pseudocode, which is Algorithm \ref{AlgoFull}\footnote{Algorithm \ref{AlgoFull} uses the automatic procedures provided by Algorithm \ref{Algo***}, Algorithm \ref{AlgoOmegagNonAbelm}, and Algorithm \ref{AlgoDeltaNonAbelm}. Often, when we refer to Algorithm \ref{AlgoFull}, we actually mean all the four algorithms (\ref{AlgoFull}, \ref{Algo***}, \ref{AlgoOmegagNonAbelm}, and \ref{AlgoDeltaNonAbelm}).}.
This chapter provides a description of Algorithm \ref{AlgoFull}. Its scope is to make any non specialist user, with only a very basic mathematical background, able to implement this algorithm to any specific system. In particular, the implementation of Algorithm \ref{AlgoFull} does not require to entirely read \cite{IF22} (where the reader is punctually addressed to obtain some definitions and, of course,  
the interested reader can find all the theoretical foundation and proofs that demonstrate the validity and all the convergence properties of the algorithm).

Algorithm \ref{AlgoFull} is recursive. Given a system characterized by (\ref{EquationSystemDefinitionUIO}), it automatically builds the entire observation space in a finite number of iterations. In particular, it builds the so-called {\it observability codistribution}, which is, by definition, the span of the gradients of all the observable functions.
Algorithm \ref{AlgoFull} builds the generators of the observability codistribution. To this regard, 
note that, in the algorithm, statements like "${\Omega}=\textnormal{span}\left\{\nabla h_1,\ldots,\nabla h_p\right\}$" (e.g., Line \ref{AlgoLineOmegaINIT} of Algorithm \ref{AlgoFull}) require no action. The algorithm builds the codistribution by building a set of its generators. When it recursively updates the codistribution it simply computes its new generators (and this is obtained by applying certain operations directly on the generators of the previous codistribution or by computing new covectors starting from the quantities that define the original system $\Sigma$ in (\ref{EquationSystemDefinitionUIO})).
In the following sections, we provide all the ingredients necessary to implement Algorithm \ref{AlgoFull}.

\section{The algorithm}

\begin{al} 

\begin{algorithmic}[1]
~

\State Set $\Sigma$ the system in (\ref{EquationSystemDefinitionUIO}), 
$\Omega=\textnormal{span}\left\{\nabla h_1,\ldots,\nabla h_p\right\}$, boolean Continue=True.\label{AlgoLineOmegaINIT}

\Loop\Comment{{\blue This is the MAIN loop}}

\If{$\uideg \left({\Omega}\right)==m_w$ }
{\bf break} main loop.
\EndIf

\State $m=\uideg \left({\Omega}\right)$.

\State $[\tih_1,\ldots, \tih_m]=\mathcal{S}(\Sigma,~\Omega)$,~~
$\Sigma=\mathcal{R}(\Sigma,~\tih_1,\ldots, \tih_m)$.
\State \textnormal{Compute} $\mu$, $\nu$, $\widehat{g}^\alpha$ (Eqs (\ref{EquationTensorMSynchrom}-\ref{Equationgalpham})).

\State Run Algorithm \ref{Algo***}\label{AlgoLineRunAlgo***}

\If{Finish}
$\OBS=\Omega_*$, $\E=\Sigma$, Continue=false, and
{\bf break} main loop.
\EndIf

\If{$m_u>0$}

\State $\Sigma=\Sigma_*$, and Compute $\mu$, $\nu$, $\giat^\alpha$ (Eqs (\ref{EquationTensorMSynchrom}-\ref{Equationgalpham})).

\State Compute $\chi_*=\psi_{k_*-1}^{i_*}$\Comment{{\blue Carried out recursively in $\mathcal{W}$. Then, it is expressed in $\mathcal{V}$.}}
\State $\Omega=\Omega+
 \textnormal{span} \left \{\Li_{\chi_*} \nabla \tih_{q_*} \right \}$\Comment{{\blue When $k_*=1$ $\chi_*=f^{i_*}$}}

\State $\Omega=\left<\left.f^1,\ldots,f^{m_u}, \left\{ 
\begin{array}{lll}
 & \dtv{\Li}_{\giat}&\\
 \giat^{m+1}, & \ldots, & \giat^{m_w}\\
\end{array}
\right\}~\right|\Omega\right>$. \Comment{{\blue This is the nested loop of Algorithm 9 in \cite{IF22}}}\label{AlgoLineNested}

\Else

\State $\Omega=\left<\left. \left\{ 
\begin{array}{lll}
 & \dtv{\Li}_{\giat}&\\
 \giat^{m+1}, & \ldots, & \giat^{m_w}\\
\end{array}
\right\}~\right|\Omega\right>$. \Comment{{\blue This is the nested loop of Algorithm 9 in \cite{IF22}}}\label{AlgoLineNestedNonMisto}

\EndIf

\State \textnormal{Reset} $[\Sigma, ~\Omega]=\mathcal{A}^{-}(\Sigma,~\Omega)$.\Comment{{\blue This pass again in $\mathcal{W}$.}}


\EndLoop

\If{Continue}\label{AlgoLineFINIT}
\State $[\widetilde{h}_1,\ldots, \widetilde{h}_{m_w}]=\mathcal{S}(\Sigma,~\Omega)$.\label{AlgoLineFINITSelection}
\State \textnormal{Compute} $\mu$, $\nu$, $\widehat{g}^\alpha$ (Eqs (\ref{EquationTensorMSynchro}-\ref{Equationgalpha})), \textnormal{and} $\tobs$ (Eq. (\ref{EquationTOBSDef})).\label{AlgoLineFINITMuNuTOBS}
    \State $\OBS=\left<f^1,\ldots,f^{m_u}, \widehat{g}^0,\widehat{g}^1,\ldots,\widehat{g}^{m_w}~\left|~\Omega+\tobs\right.\right>$ \textnormal{and set} $\E=\Sigma$.\label{AlgoLineFINALSTEP}
\EndIf\label{AlgoLineFINITEND}

\end{algorithmic}
\label{AlgoFull}
\end{al}

\begin{al}

\begin{algorithmic}[1]
~

\If{$(m_u==0)|(m==0)$}
\State $\Omega_*=\left<\left. \giat^0, \ldots, \giat^m~\right|\Omega\right>$.
\If{($\giat^{m+1}\in\Omega_*^\bot) \And(\giat^{m+2}\in\Omega_*^\bot) \And\ldots\And(\giat^{m_w}\in\Omega_*^\bot$)}
\State Finish=True.
\Else
\State Finish=False.
\EndIf

\Else\label{AlgoLine***1}
\State Run Algorithms \ref{AlgoOmegagNonAbelm} and \ref{AlgoDeltaNonAbelm} to compute $\widehat{k}^m=s^m_x+r^m$.

\For{$k=1,\ldots,\widehat{k}^m$}

\For{$i=1,\ldots,m_u$~
$\And$~$q=1,\ldots,m$}

  \State $\Omega=\Omega+
 \sum_{\alpha_1=0}^m\ldots\sum_{\alpha_{k-1}=0}^m
 \textnormal{span} \left \{\mathcal{L}_{[f^i]^{(\alpha_1,\ldots,\alpha_{k-1})}} \nabla \tih_q \right \}$

\State $\Omega_*=\left<\left.f^1,\ldots,f^{m_u}, 
 \giat^0, \ldots, \giat^m~\right|\Omega\right>$.

\If{($\giat^{m+1}\in\Omega_*^\bot)\And(\giat^{m+2}\in\Omega_*^\bot)\And\ldots\And(\giat^{m_w}\in\Omega_*^\bot$)}
\State Reset $[\Sigma, ~\Omega]=\mathcal{A}(\Sigma,~m,~\Omega_*)$.\Comment{{\blue This adds $W^{(k-1)}_k$.}}
\State Compute $m$, $\mu$, $\nu$, $\widehat{g}^\alpha$ (Eqs (\ref{EquationTensorMSynchrom}-\ref{Equationgalpham})).
\Else
\State Finish=False. 
\State $\Sigma_*=\Sigma$, $k_*=k$, $i_*=i$, $q_*=q$.
\State {\bf return}
\EndIf

\EndFor
\EndFor

\State Finish=True.

\EndIf

\end{algorithmic}
\label{Algo***}
\end{al}

\begin{al}

\begin{algorithmic}[1]
~
\State  \textnormal{Set $\Sigma'=\Sigma$, $m$, $\tih_1,\ldots,\tih_{m}$, as at Line \ref{AlgoLineRunAlgo***} of Algorithm \ref{AlgoFull}}.
 \State \textnormal{Set $k=0$, $\Omega^m_0=\textnormal{span}\left\{\nabla\tih_1,\ldots,\nabla\tih_{m}\right\}$, and $^x\Omega^m_0=\Omega^m_0$}.
\Loop
\State \textnormal{Set} $k=k+1$.
\State \textnormal{Reset} $[\Sigma', ~\Omega^m_{k-1}]=\mathcal{A}(\Sigma',~m, ~\Omega^m_{k-1})$\label{AlgoOmegaLineRESET}
\State $\Omega^m_k=\Omega^m_{k-1}+\dt{\Li}_{g^0} \Omega^m_{k-1}+\sum_{j=1}^{m}\Li_{g^j} \Omega^m_{k-1}$.\label{AlgoOmegaLineSum}

\State $^x\Omega^m_k=\mathcal{D}_x(\Omega^m_k)$.
 \If{$^x\Omega^m_k==~^x\Omega^m_{k-1}$}
 \State \textnormal{Set $s^m_x=k$ {\bf then exit}}.
 \EndIf
 \EndLoop
\end{algorithmic}\label{AlgoOmegagNonAbelm}
\end{al}

\begin{al}
\begin{algorithmic}[1]
~
\State  \textnormal{Set $\Sigma'=\Sigma$ and $m$ as at Line \ref{AlgoLineRunAlgo***} of Algorithm \ref{AlgoFull}}.
 \State \textnormal{Set $k=0$, $\Delta_0=\textnormal{span}\left\{f^1,\ldots,f^{m_u} \right\}$}.

\Loop
\State \textnormal{Set} $k=k+1$.
\State \textnormal{Reset} $[\Sigma', ~\Delta_{k-1}]=\mathcal{A}(\Sigma',~m, ~\Delta_{k-1})$\label{AlgoDeltaLineRESET}

\State $\Delta_k=~\Delta_{k-1}+\sum_{\beta=0}^m ~\left[\Delta_{k-1}\right]^\beta$

 \If{$\Delta_k~==~\Delta_{k-1}$}
 \State \textnormal{Set $r^m=k-1$ {\bf then exit}}.
 \EndIf

\EndLoop

\end{algorithmic}\label{AlgoDeltaNonAbelm}
\end{al}

\section{Basic operations in Algorithm \ref{AlgoFull}}\label{SectionUIOBasicOperations}

\subsection{$\boldsymbol{\deg_w(\Omega)}$ operation}\label{SubSectionDeg}

In Section \ref{SectionDefinitionCanonic}, we provided the definition of unknown input degree of reconstructability from a set of scalar functions.
Given an integrable codistribution
$\Omega=\textnormal{span}
\left\{\nabla\lambda_1,\ldots, \nabla\lambda_k\right\}$,
we define the unknown input degree of reconstructability of $\Sigma$ from $\Omega$ the unknown input degree of reconstructability from a set of generators of $\Omega$ (e.g., from $\lambda_1,\ldots, \lambda_k$). We denote it by $\uideg\left(\Omega\right)$.

 The main loop of Algorithm \ref{AlgoFull} starts by checking if $\uideg(\Omega)==m_w$.
 If this is the case, the main loop is interrupted.
Then, Lines \ref{AlgoLineFINIT}-\ref{AlgoLineFINALSTEP} of Algorithm \ref{AlgoFull} are executed and Algorithm \ref{AlgoFull} ends. For the clarity sake, we devote Section \ref{SectionSolutionUIOCanonic} to describe the behaviour of Algorithm  \ref{AlgoFull} in this case. Before, we provide the definition of further operations adopted by Algorithm \ref{AlgoFull}.

\subsection{$\boldsymbol{\mathcal{S}(\Sigma,~\Omega)}$ operation}\label{SubSectionSelection}

Let us consider a system $\Sigma$ that is characterized by (\ref{EquationSystemDefinitionUIO}) and an integrable codistribution $\Omega$. Let us set $m=\uideg\left(\Omega\right)$. The $\mathcal{S}(\Sigma,~\Omega)$ operation provides a set of $m$ scalar functions, denoted by $\widetilde{h}_1,\ldots, \widetilde{h}_m$, that make the unknown input reconstructability matrix full rank.
The operation is denoted by $\mathcal{S}$ because is the {\it Selection} of the aforementioned $\widetilde{h}_1,\ldots, \widetilde{h}_m$ from the generators of $\Omega$.
The execution of this operation is immediate because the generators of $\Omega$ are always available when executing Algorithm \ref{AlgoFull}.
In this paper, we use the notation $[\widetilde{h}_1,\ldots, \widetilde{h}_m]=\mathcal{S}(\Sigma,~\Omega)$ to denote the outputs of this operation. 

\subsection{$\boldsymbol{\mathcal{R}(\Sigma,~\lambda_1,\ldots, \lambda_m)}$ operation}\label{SubSectionReorder}

Let us consider a system $\Sigma$ that is characterized by (\ref{EquationSystemDefinitionUIO}) and a set of $m<m_w$ scalar functions, $\lambda_1,\ldots, \lambda_m$, such that $\textnormal{rank}\left(\mathcal{RM}\left(\lambda_1,\ldots, \lambda_m\right)
\right)=m$. 
The $\mathcal{R}(\Sigma,~\lambda_1,\ldots, \lambda_m)$ operation provides a new system that is obtained from $\Sigma$ by reordering the unknown inputs. In particular, they are reordered in such a way that the square submatrix that consists of the first $m$ columns of the unknown input reconstructability matrix from 
$\lambda_1,\ldots, \lambda_m$, is non singular.
The operation is denoted by $\mathcal{R}$ because it is a {\it Reordering} of the unknown inputs, as explained above.
Its execution is immediate. It suffices to extract from
$\mathcal{RM}\left(\lambda_1,\ldots, \lambda_m \right)$ a set of $m$ independent columns.
In this paper, we use the notation $\Sigma'=\mathcal{R}(\Sigma,~\lambda_1,\ldots, \lambda_m)$ to denote the new ordered system.

\subsection{$\boldsymbol{\left.\left<\tau^1,\ldots,\tau^d~\right|~\Omega\right>}$ operation}\label{SubSectionOMin}

Given a codistribution $\Omega$ and a set of vector fields $\tau^1,\ldots,\tau^d$, $\left.\left<\tau^1,\ldots,\tau^d~\right|~\Omega\right>$ denotes
the smallest codistribution that contains $\Omega$ and such that, for any $\omega\in\Omega$, we have $\Li_{\tau^i}\omega\in\left.\left<\tau^1,\ldots,\tau^d~\right|~\Omega\right>$, for any $i=1,\ldots,d$.
This minimal codistribution 
can be easily computed by a simple recursive algorithm, which is \cite{Isi95}:

\begin{equation}\label{EquationAlgorithmsMinimalCod}
\left\{\begin{array}{ll}
 \Omega_0 &=  \Omega\\
 \Omega_{k+1} &=  \Omega_k+\sum_{i=1}^d\Li_{\tau^i}\Omega_k\\
\end{array}\right.,
\end{equation}

\noindent This operation is used in Algorithm \ref{AlgoFull} at Line \ref{AlgoLineFINALSTEP} (
note that, for TV systems, the Lie derivative operator along $\widehat{g}^0$, i.e., $\Li_{\widehat{g}^0}$, must be replaced by  $\dt{\Li}_{\widehat{g}^0}\triangleq \Li_{\widehat{g}^0}+\frac{\partial}{\partial t}$).

\begin{lrbox}{\mybox}
$\boldsymbol{\left<\left.\tau^1,\ldots,\tau^{d_1}, \left\{ 
\begin{array}{lll}
 \xi^1, & \ldots, & \xi^{d_2}\\
 \zeta^1, & \ldots, & \zeta^{d_3}\\
\end{array}
\right\}~\right|\Omega\right>}$
\end{lrbox}

\subsection{\usebox{\mybox} operation}\label{SubSectionOMinSpecial}

Given the codistribution $\Omega$ and the three sets of vector fields: $\tau^1,\ldots,\tau^{d_1}$, $\xi^1,\ldots,\xi^{d_2}$, and $\zeta^1,\ldots,\zeta^{d_3}$, we denote by
\begin{equation}\label{EquationOmgInvIf}
\left<\left.\tau^1,\ldots,\tau^{d_1}, \left\{ 
\begin{array}{lll}
 \xi^1, & \ldots, & \xi^{d_2}\\
 \zeta^1, & \ldots, & \zeta^{d_3}\\
\end{array}
\right\}~\right|\Omega\right>
\end{equation}
the codistribution computed by the algorithm in Equation (\ref{EquationAlgorithmsMinimalCod}), where the sum $\sum_{i=1}^d$ at the recursive step 
includes all the $d_1$ vector fields $\tau^1,\ldots,\tau^{d_1}$. In addition, by denoting with $\omega$ a generator of $\Omega_k$, the codistribution $\Omega_{k+1}$ also includes, among its generators, all the $d_2$ Lie derivatives of $\omega$ along the vector fields $\xi^1,\ldots,\xi^{d_2}$ if and only if $\Li_{\zeta^l}\omega$ vanishes for all $l=1,\ldots,d_3$. The algorithm is then interrupted at the smallest integer $j$ such that $\Omega_j=\Omega_{j-1}$.


\subsection{$\boldsymbol{\mathcal{A}(\Sigma,~m)}$ operation}\label{SubSectionSigma++}
Let us consider a system $\Sigma$ that is characterized by (\ref{EquationSystemDefinitionUIO}) and an integer $m<m_w$. 
This operation provides a new extended system, 
defined as follows. It is obtained by introducing a 
new extended state that includes the last $d=m_w-m$ unknown inputs. We have:

\begin{equation}\label{EquationAugmentationState}
x\rightarrow
[x^T, w_{m+1}, \ldots, w_{m_w}]^T
\end{equation}

Starting from (\ref{EquationSystemDefinitionUIO}) we obtain the dynamics of the above extended state. We obtain a new system that still satisfies (\ref{EquationSystemDefinitionUIO}). It is still characterized by $m_u$ known inputs and $m_w$ unknown inputs. All the $m_u$ known inputs and the first $m$ unknown inputs
coincide with the original ones. The last $d$ unknown inputs are $\dt{w}_{m+j}$, $j=1,\ldots,d$.
Regarding the new vector fields that describe the dynamics we obtain:

\begin{equation}\label{EquationAugmentationSystem}
f^i \rightarrow
\left[\begin{array}{c}
   f^i \\
   0_d \\
\end{array}
\right]
,~~
g^0 \rightarrow
\left[\begin{array}{c}
   g^0 + \sum_{l=m+1}^{m_w}  g^l w_l\\
   0_d \\
\end{array}
\right]
,\\
\end{equation}
\[
g^k \rightarrow
\left[\begin{array}{c}
   g^k \\
   0_d \\
\end{array}
\right],~~
~~
g^{m+j} \rightarrow
\left[\begin{array}{c}
   0_n \\
   e^j \\
\end{array}
\right],
\]
with
$i=1,\ldots,m_u$,
$k=1,\ldots,m$, $j=1,\ldots,d$. In addition, $0_d$ and $0_n$ denote the zero $d$-column vector and the zero $n$-column vector, respectively, and $e^j$ denotes the $d$-column vector with the $j^{th}$ entry equal to 1 and the remaining $d-1$ entries equal to 0.

The operation is denoted by $\mathcal{A}$ because it consists of a state {\it Augmentation} and the consequent re-definition of all the key vector fields that characterize the new extended system, as specified above.
In this paper, we use the notation $\Sigma'=\mathcal{A}(\Sigma,~m)$ to denote the new extended system.


In Algorithms \ref{AlgoFull}, \ref{AlgoOmegagNonAbelm}, and \ref{AlgoDeltaNonAbelm}, we also use this operation by including further outputs and further inputs. Specifically, we consider the following two cases:

\begin{enumerate}

\item $[\Sigma', ~\Omega']=\mathcal{A}(\Sigma,~m,~\Omega)$, where the second output ($\Omega'$) is trivially the augmented codistribution obtained by extending all the covectors of $\Omega$ with $d$ zero entries.

\item $[\Sigma', ~\Delta']=\mathcal{A}(\Sigma,~m, ~\Delta)$, where the second output ($\Delta'$) is trivially the augmented distribution obtained by extending all the vectors of $\Delta$ with $d$ zero entries.

\end{enumerate}

%

We can apply this operation multiple consecutive times (and this is the case of
Algorithms \ref{AlgoFull}, \ref{AlgoOmegagNonAbelm}, and \ref{AlgoDeltaNonAbelm}, as the operation appears in a loop).
The resulting systems still satisfy (\ref{EquationSystemDefinitionUIO}) and they are all characterized by $m_u$ known inputs and $m_w$ unknown inputs. In \cite{IF22}, we called these extended systems {\it Finite Unknown Inputs Extensions}. In addition, in \cite{IF22}, we introduced the following definition:

\begin{definition}[Highest UI Degree of Reconstructability]\label{DefinitionHDegUIReconstr}
Given the system in (\ref{EquationSystemDefinitionUIO}), the highest unknown input degree of reconstructability is the largest unknown input degree of reconstructability of all its finite unknown input extensions.
\end{definition}

Finally, in \cite{IF22} we called a system {\it Canonizable with respect to its unknown inputs} if its highest unknown input degree of reconstructability is equal to $m_w$.

\subsection{$\boldsymbol{\mathcal{A}^-(\Sigma,~\Omega)}$ operation}\label{SubSectionSigma--}
This operation is only executed when, in the iteration of the main loop under execution, the Boolean Finish=False, and the remaining steps of this iteration determines a new observable function that increases the unknown input degree of reconstructability. In this case, the state can be augmented by including the unknown inputs with index larger than $m$ (through the operation $\mathcal{A}$). During this iteration of the main loop, the observable function, which increases the unknown input degree of reconstructability, is added to $\Omega$.
 Let us denote this function by $\theta$.
In general, $\theta$ also depends on the quantities $v_\alpha$, defined in (\ref{Equationvi}). 
The first operation executed by $\mathcal{A}^-$ is to express $\theta$ only in terms of the unknown inputs and the original state. In other words, all the $v_\alpha$ that appear in $\theta$ are expressed in terms of the unknown inputs (this is obtained by using (\ref{Equationvi}) and its time derivatives when $\theta$ also depends on the time derivatives of $v_\alpha$).
Then, all the unknown inputs that appear in $\theta$ are added to the original state. The resulting augmented state defines the new system $\Sigma$.

\section{Algorithm \ref{AlgoFull} for systems canonic with respect to the unknown inputs}\label{SectionSolutionUIOCanonic}

 When a system is canonic with respect to its unknown inputs, Algorithm \ref{AlgoFull} ends with the execution of Lines \ref{AlgoLineFINIT}-\ref{AlgoLineFINALSTEP} (the boolean variable {\it Continue} remains set to "true"). The system could be directly in canonic form with respect to its unknown inputs
or $\Omega$ was obtained in other parts of the algorithm~
and the system has been set in canonical form after the determination of one or more observable functions that are not among the outputs.
~The execution of Line \ref{AlgoLineFINITSelection} of Algorithm \ref{AlgoFull} provides the $m_w$ scalar functions $\widetilde{h}_1,\ldots, \widetilde{h}_{m_w}$
($[\widetilde{h}_1,\ldots, \widetilde{h}_{m_w}]=\mathcal{S}(\Sigma,~\Omega)$).
Then, the algorithm
computes 
the following quantities: $\mu$, $\nu$, $\widehat{g}^\alpha$, and the codistribution $\tobs$.

\subsection{$\boldsymbol{\mu}$ and $\boldsymbol{\nu}$}\label{SubSectionMuNu}

The two-index tensor $\mu$ is defined as follows:

\begin{equation}\label{EquationTensorMSynchro}
\mu^i_j = \mathcal{L}_{g^i}\widetilde{h}_j, ~~~i,~ j=1,\ldots,m_w
\end{equation}
\[
\mu^0_0 = 1,~
\mu^i_0=0,~
\mu^0_i = \frac{\partial\widetilde{h}_i}{\partial t}  +
\mathcal{L}_{g^0}\widetilde{h}_i,~  i=1,\ldots,m_w.
\]

Note that the entries with $i,~ j=1,\ldots,m_w$ are the same entries of the unknown input reconstructability matrix from $\widetilde{h}_1,\ldots, \widetilde{h}_{m_w}$, which is full rank. As a result, the tensor $\mu$ is also non singular and 
we denote by $\nu$ its inverse. 

\subsection{$\boldsymbol{\widehat{g}^{\alpha}}$}

Starting from $\nu$ the algorithm builds the new vector fields $\widehat{g}^0,\ldots,\widehat{g}^{m_w}$, defined as follows:
\begin{equation}\label{Equationgalpha}
   \widehat{g}^{\alpha}  = \sum_{\beta=0}^{m_w}\nu^{\alpha}_{\beta} g^{\beta},~~\alpha=0,1,\ldots,m_w.
\end{equation}

\subsection{The codistribution $\boldsymbol{\tobs}$}\label{SubSectionTOBS}
The codistribution $\tobs$ is defined as follows:

\begin{equation}\label{EquationTOBSDef}
\tobs:=
\sum_{q=1}^{m_w}
\sum_{j=0}^{\NO+\ND}
\sum_{\alpha_1=0}^{m_w}\ldots\sum_{\alpha_j=0}^{m_w}
\sum_{i=1}^{m_u}
\textnormal{span}\left\{
\nabla\Li_{[f^i]^{(\alpha_1,\ldots,\alpha_j)}}\tih_q
\right\},
\end{equation}

where the second sum is up to $s+r$, and the integers $s$ and $r$ are defined as follows:

\begin{itemize}

\item  $\NO$ is the smallest integer such that $\Omega_s=\Omega_{s-1}$, where $\Omega_k$ is the codistribution at the $k^{th}$ step of the algorithm in (\ref{EquationAlgorithmsMinimalCod}), when computing:

\[
\left<g^0,g^1,\ldots,g^{m_w}~\left|~\textnormal{span}\left\{\nabla\tih_1,\ldots,\nabla\tih_{m_w}\right\}\right.\right>
\]
Note that, for TV systems, the Lie derivative operator along $g^0$, i.e., $\Li_{g^0}$, must be replaced by  $\dt{\Li}_{g^0}:=\Li_{g^0}+\frac{\partial}{\partial t}$. Note that $\NO\le n-m_w+1$.

\item  $\ND$ is the smallest integer such that $\Delta_{r+1}=\Delta_r$, where $\Delta_k$ is the distribution at the $k^{th}$ step of the algorithm that computes:

\[
\left<\left[\cdot\right]^\cdot~\left|~\textnormal{span}\left\{f^1,\ldots,f^{m_u}\right.\right\}\right>
\]
Note that $\ND\le n-1$ (it is even
$\ND\le n-\textnormal{dim}\left\{\textnormal{span}\left\{f^1,\ldots,f^{m_u} \right\}\right\}$).

\end{itemize}

\subsection{Final step}

The last operation of Algorithm \ref{AlgoFull}, when the system is canonic with respect to its unknown inputs, is the computation of the observability codistribution:

\[
\OBS=\left<f^1,\ldots,f^{m_u}, \widehat{g}^0,\widehat{g}^1,\ldots,\widehat{g}^{m_w}~\left|~\Omega+\tobs\right.\right>
\]

This is obtained by running the algorithm in (\ref{EquationAlgorithmsMinimalCod}) for the specific case (note that, for TV systems, the Lie derivative operator along $\widehat{g}^0$, i.e., $\Li_{\widehat{g}^0}$, must be replaced by  $\dt{\Li}_{\widehat{g}^0}:=\Li_{\widehat{g}^0}+\frac{\partial}{\partial t}$). The convergence
of the algorithm in (\ref{EquationAlgorithmsMinimalCod}) is attained 
at the smallest integer $j$ such that 
$\Omega_{j}=\Omega_{j-1}$. Note  that $j$ cannot exceed $n-\dim\left(\Omega+\tobs\right)+1$.

\section{Algorithm \ref{AlgoFull} for systems that are not in canonical form}\label{SectionSolutionUIONonCanonic}

Let us back to the first part of Algorithm \ref{AlgoFull}, 
at the beginning of the main loop, 
and let us consider now the case when 
$\uideg \left({\Omega}\right)<m_w$.
The system is not in canonical form with respect to its unknown inputs. 
The algorithm continues by computing
the set of scalar functions 
$\widetilde{h}_1,\ldots, \widetilde{h}_m$. They are obtained by executing the operation $\mathcal{S}(\Sigma,~\Omega)$ defined in Section \ref{SubSectionSelection}, and the unknown inputs are reordered accordingly ($\mathcal{R}(\Sigma,~\widetilde{h}_1,\ldots, \widetilde{h}_m)$, Section \ref{SubSectionReorder}).
Then, the algorithm
computes the following quantities: $\mu$, $\nu$, and $\widehat{g}^\alpha$.

\subsection{$\boldsymbol{\mu}$ and $\boldsymbol{\nu}$}\label{SubSectionMumNum}

The two-index tensor $\mu$ is defined as follows:

\begin{equation}\label{EquationTensorMSynchrom}
\mu^i_j = \mathcal{L}_{g^i}\widetilde{h}_j, ~~~i,~ j=1,\ldots,m
\end{equation}
\[
\mu^0_0 = 1,~
\mu^i_0=0,~
\mu^0_i = \frac{\partial\widetilde{h}_i}{\partial t}  +
\mathcal{L}_{g^0}\widetilde{h}_i,~  i=1,\ldots,m.
\]

Note that the entries with $i,~ j=1,\ldots,m$ are the same entries of the matrix that consists of the first $m$ columns of the unknown input reconstructability matrix from $\widetilde{h}_1,\ldots, \widetilde{h}_{m}$. By construction, this submatrix is full rank. As a result, the tensor $\mu$ is also non singular and 
we denote by $\nu$ its inverse. 

\subsection{$\boldsymbol{\widehat{g}^{\alpha}}$}

Starting from $\nu$ we build the new vector fields $\widehat{g}^0,\ldots,\widehat{g}^{m}$ and $\giat^{m+1},\ldots,\giat^{m_w}$. They are defined as follows:
\begin{equation}\label{Equationgalpham}
   \widehat{g}^{\alpha}  = \sum_{\beta=0}^m~\nu^{\alpha}_{\beta} ~g^{\beta},~~\alpha=0,1,\ldots,m,~~~~\giat^k:=g^k-\sum_{\alpha=0}^m\giat^\alpha\Li_{g^k}\tih_\alpha, ~~k=m+1,\ldots,m_w
\end{equation}

%
%

\subsection{Algorithm \ref{Algo***}}\label{SubSectionAlgo***}

This algorithm returns the Boolean Finish. When Finish==True, it means that $\Omega_*$ is the observability codistribution. When Finish==False, the execution of the algorithm only tells us that there exists an observable function that increases the unknown inputs degree of reconstructability. Let us denote this function by $\theta$. Algorithm \ref{Algo***} does not provide $\theta$. $\theta$ is determined by the remaining steps of the iteration of the main loop of Algorithm \ref{AlgoFull}, under execution. Note that, when Finish==False, Algorithm \ref{Algo***} 
returns $\Sigma_*$, which  can differ from $\Sigma$ and is the suitable system to compute $\theta$. In addition, when Finish==False, Algorithm \ref{Algo***} provides
$k_*$, $i_*$, and $q_*$, which are used in the main loop of Algorithm \ref{AlgoFull} to compute $\theta$.

\subsection{The $\boldsymbol{\giat}$ vector field}\label{SubSectionDeg}

It appears at Line \ref{AlgoLineNested} of Algorithm \ref{AlgoFull}. It is defined as follows:

\begin{equation}\label{Equationgiatdinfty}
~\giat\triangleq\sum_{\beta=0}^m ~\giat^\beta~v_\beta
\end{equation}
with ($\beta=0,1,\ldots,m$)

\begin{equation}\label{Equationvi}
v_\beta\triangleq
\left[\begin{array}{ll}
\sum_{\gamma=0}^m
~\mu^\gamma_\beta~w_\gamma + \sum_{k=m+1}^{m_w}(\Li_{g^k}\tih_\beta)~ w_k& k_*=1\\
\sum_{\gamma=0}^m
~\mu^\gamma_\beta~w_\gamma & k_*>1\\
\end{array}
\right.
\end{equation}

\noindent The operator $\dtv{\Li}_\chi$ is:

\begin{equation}\label{EquationDTV}
\dtv{\Li}_{\chi}\triangleq\left(
\Li_\chi + \sum_{\alpha=0}^m\sum_{i_\alpha=0}^\infty
v_\alpha^{(i_\alpha+1)}
\frac{\partial}{\partial v_\alpha^{(i_\alpha)}} 
\right).
\end{equation}
For time varying systems $\dtv{\Li}_{\chi}\rightarrow \dtv{\Li}_{\chi}+\frac{\partial}{\partial t}$.

\vskip.1cm
\noindent The expression of $\psi_k^i$ is obtained automatically and recursively.
We have:

\begin{equation}\label{Equationpsik+1}
\psi_0^i=f^i,~~~~~\psi_{k+1}^i=\sum_{\gamma=0}^m[g^\gamma,~\psi^i_q]w_\gamma
+\sum_{\alpha=0}^m\sum_{i_\alpha=0}^k\frac{\partial \psi_k^i}{\partial w_\alpha^{(i_\alpha)}}w_\alpha^{(i_\alpha+1)}.
\end{equation}
For time-varying systems $\psi_{k+1}^i\rightarrow\psi_{k+1}^i+\frac{\partial \psi_k^i}{\partial t}$.

Once determined, we must eliminate $w_\alpha$, $w_\alpha^{(1)}$, $\ldots$, $w_\alpha^{(i_\alpha+1)}$ by expressing them in terms of 
$v_\alpha$, $v_\alpha^{(1)}$, $\ldots$, $v_\alpha^{(i_\alpha+1)}$. This is achieved by using the inverse of Equation (\ref{Equationvi}) above, i.e.,
$w_\alpha= \sum_{\beta=0}^m~\nu_\alpha^\beta v_\beta$.

\subsection*{Algorithm \ref{AlgoOmegagNonAbelm}}
It computes the codistribution $^x\Omega^m$. 
This algorithm uses a new operation denoted by $\mathcal{D}_x(\Omega)$. The input of this operation is a codistribution that is defined in a given augmented space (in the algorithm, before this operation, the $\mathcal{A}(\Sigma,~m)$ operation is applied). The operation $\mathcal{D}_x(\Omega)$, first detects a basis of $\Omega$.
As $\Omega$ is an integrable codistribution, it detects a set of scalar functions: $\theta_1,\ldots,\theta_D$, where $D$ is the dimension of $\Omega$.
In other words, 
$\Omega=\textnormal{span}\left\{\nabla\theta_1,\ldots,\nabla\theta_D\right\}$, where $\nabla$ is the gradient with respect to the new extended state (i.e., the state defined at the last $\mathcal{A}(\Sigma,~m)$ operation). The output of the operation $\mathcal{D}_x(\Omega)$, is $^x\Omega=\textnormal{span}\left\{\partial_x\theta_1,\ldots,\partial_x\theta_D\right\}$, where $\partial_x$ is the gradient with respect to the original state, i.e., the state of $\Sigma$
defined at the first line of Algorithm \ref{AlgoOmegagNonAbelm}.

The initialization step sets the codistribution equal to the span of the gradients of the selected functions 
$\widetilde{h}_1,\ldots, \widetilde{h}_{m}$. As a result, the dimension of $\Omega^m_0$ is $m$. Then, each iteration of the loop executes the following operations:

\begin{enumerate}

\item System augmentation, as explained in Section \ref{SubSectionSigma++}.

\item Computation of $\Omega^m_k$ by adding to $\Omega^m_{k-1}$ the term $\dt{\Li}_{g^0} \Omega^m_{k-1}+\sum_{j=1}^{m}\Li_{g^j} \Omega^m_{k-1}$.

\item Computation of $^x\Omega^m_k=\mathcal{D}_x(\Omega^m_k)$.

\end{enumerate}

The convergence of Algorithm \ref{AlgoOmegagNonAbelm} occurs at the smallest integer $j$ such that $^x\Omega^m_j=~^x\Omega^m_{j-1}$ and $j\le n - m+1$. We set $s^m_x=j$.
Note that, because of the presence of the state augmentation in the loop (Line \ref{AlgoOmegaLineRESET} of Algorithm \ref{AlgoOmegagNonAbelm}), proving the validity of the above convergence property is very demanding (see Appendix E of \cite{IF22} and in particular the proof of Proposition 11 in that appendix).

\subsection*{Algorithm \ref{AlgoDeltaNonAbelm}}
It computes the distribution $\Delta$. 
The initialization step sets the distribution equal to the span of the vector fields $f^1,\ldots,f^{m_u}$. Then, the recursive step adds to $\Delta_{k-1}$ the term
$\sum_{\beta=0}^{m_w}\left[\Delta_{k-1}\right]^\beta$.
The convergence is attained at the smallest integer $j$ such that $\Delta_j=\Delta_{j-1}$ and $j\le n-\textnormal{dim}\left\{\textnormal{span}\left\{f^1,\ldots,f^{m_u} \right\}\right\}+1$.
We set $r^m+1=j$.
Note that, because of the presence of the state augmentation in the loop (Line \ref{AlgoDeltaLineRESET} of Algorithm \ref{AlgoDeltaNonAbelm}), proving the validity of the above convergence property is non trivial (see Section 6.2.1 of \cite{IF22} and in particular the proof of Proposition 1 given in Appendix D of \cite{IF22}).

\vskip .2cm

Note that, the system augmentation performed by the $\mathcal{A}$ operation executed by Algorithms \ref{AlgoOmegagNonAbelm} and \ref{AlgoDeltaNonAbelm}, does not reset the system in Algorithm \ref{AlgoFull}. 
The execution of Algorithm \ref{AlgoOmegagNonAbelm} and Algorithm \ref{AlgoDeltaNonAbelm} only returns the integer $\widehat{k}^m=s^m_x+r^m$. All remaining quantities remain unchanged.


\chapter{Unknown input reconstruction}\label{ChapterUIReconstruction}

The content of this chapter was removed. The reader can find  it in \cite{IF22}.
In addition, he/she can find a complete and exhaustive study of this problem in \cite{arXivODE}.

\chapter{Applications}\label{ChapterApplication}

This section illustrates the implementation of Algorithm \ref{AlgoFull}. 
We refer to the visual inertial sensor fusion problem. For the sake of clarity, we restrict our analysis to a $2D$ environment. The $3D$ case differs only in a more laborious calculation. 
We consider three variants of this sensor fusion problem depending on the available inertial signals.
The description of these variants is given in Section \ref{SectionSystem}.
In Section \ref{SectionApplicationStateObservability}, we use Algorithm \ref{AlgoFull} to obtain the observability properties. In particular, for one of the considered variants, we provide all the details of the computation, following, step by step, the implementation of Algorithm \ref{AlgoFull} (Section \ref{SubSectionObservabilitySystem2}).

\section{The system}\label{SectionSystem}

\begin{figure}[htbp]
\begin{center}
\includegraphics[width=.8\columnwidth]{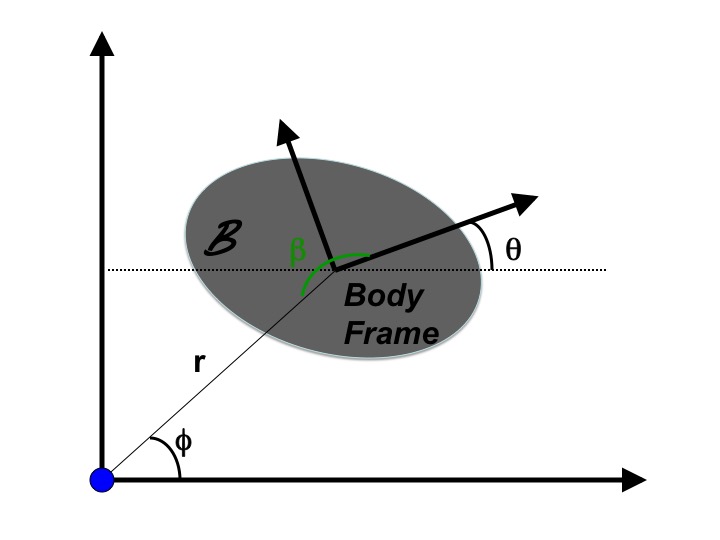}
\caption{The global frame, the body frame and the observation provided by the \VS~ sensor ($\beta$).} \label{FigVISFM2DCalibrated}
\end{center}
\end{figure}

We consider a rigid body ($\mathcal{B}$) equipped with a visual sensor (\VS) and an inertial measurement unit (\IMU). The body moves on a plane. A complete \IMU~ measures the body acceleration and the angular speed. In $2D$, the acceleration is a two dimensional vector and the angular speed is a scalar. The visual sensor provides the bearing angle of the features in its own local frame.
We assume that the \VS~ frame coincides with the \IMU~ frame. We call this frame, the body frame. In addition, we assume that the inertial measurements are unbiased.
Figure \ref{FigVISFM2DCalibrated} depicts our system.


\noindent It is very convenient to work in polar coordinates. 
Hence, we define the state: 

\begin{equation}\label{EquationApplicationVISFM2CalibratedState}
x=[r,~\phi,~v,~\alpha,~\theta]^T
\end{equation}
where $r$ and $\phi$ characterize the body position, $\theta$ its orientation (see Fig. \ref{FigVISFM2DCalibrated} for an illustration), and $v$ and $\alpha$ the body speed in polar coordinates. In particular,  $v= \sqrt{v_x^2+v_y^2}$ and $\alpha= \arctan\left(\frac{v_y}{v_x}\right)$, where $[v_x,~v_y]^T$ is the body speed in Cartesian coordinates. The dynamics are:

\begin{equation}\label{EquationApplicationVISFM2Calibrated}
\left[\begin{array}{ll}
  \dot{r} ~&= v \cos(\alpha-\phi)\\
  \dot{\phi} ~&= \frac{v}{r} \sin(\alpha-\phi)\\
  \dot{v} ~&= A_x \cos(\alpha-\theta)  +  A_y \sin(\alpha-\theta)\\
 \dot{\alpha}~&= -\frac{A_x}{v} \sin(\alpha-\theta)  + \frac{A_y}{v} \cos(\alpha-\theta)\\
  \dot{\theta} ~~&= \omega \\
\end{array}\right.
\end{equation}
where $[A_x, ~A_y]^T$ is the body acceleration in the body frame and $\omega$ the angular speed.
Without loss of generality, we assume that the feature is positioned at the origin of the global frame. The \VS~ sensor provides the angle $\beta=\pi-\theta+\phi$. Hence, we can perform the observability analysis by using the output (we ignore $\pi$):

\begin{equation}\label{EquationApplicationVISFM2CalibratedOutput}
y=h(x)=\phi-\theta
\end{equation}

 We consider the following three variants, which differ in a different setting of the \IMU~ sensor:

\begin{enumerate}

\item Variant 1: \IMU~ provides the body acceleration ($[A_x, ~A_y]^T$) and the angular speed ($\omega$).

\item Variant 2: \IMU~ consists of a single gyroscope that only provides 
the angular speed ($\omega$).

\item Variant 3: \IMU~ consists of a single accelerometer that only provides 
the component of the acceleration along the $x-$axis of the body frame ($A_x$).

\end{enumerate}

\section{State observability}\label{SectionApplicationStateObservability}

We use Algorithm \ref{AlgoFull} to obtain the observability properties of our system.
For the brevity sake, for Variant 3, we do not provide the details of the computation (which are similar to the ones for Variant 2) but only the result.

\subsection{State observability for Variant 1}\label{SubSectionObservabilitySystem1}

When $A_x$, $A_y$, and $\omega$ are known,
by comparing (\ref{EquationApplicationVISFM2Calibrated}) and (\ref{EquationApplicationVISFM2CalibratedOutput}) with (\ref{EquationSystemDefinitionUIO}) we have: $n=5$, $m_u=3$, $m_w=0$,
$p=1$,
$u_1=\omega$, $u_2=A_x$, $u_3=A_y$,

\[
g^0=\left[\begin{array}{c}
v \cos(\alpha-\phi)\\
\frac{v}{r} \sin(\alpha-\phi)\\
 0 \\
 0  \\
 0  \\
\end{array}
\right],~
f^1=\left[\begin{array}{c}
 0 \\
 0  \\
 0  \\
 0  \\
 1  \\
\end{array}
\right],
\]
\[
f^2=\left[\begin{array}{c}
 0 \\
 0  \\
 \cos(\alpha-\theta)   \\
-\frac{1}{v} \sin(\alpha-\theta)\\
 0  \\
\end{array}
\right],~
f^3=\left[\begin{array}{c}
 0 \\
 0  \\
 \sin(\alpha-\theta)   \\
\frac{1}{v} \cos(\alpha-\theta)\\
 0  \\
\end{array}
\right]
\]

As all the inputs are known, we can compute the observability codistribution by using the standard observability rank condition, which is a special case of Algorithm \ref{AlgoFull} when $m_w=0$.

We obtain, at Line \ref{AlgoLineOmegaINIT}, $\Omega=\textnormal{span}\left\{\nabla h\right\}=\textnormal{span}\left\{[0~1~0~0~-1]\right\}$.
As $m_w=0$, $\uideg(\Omega)=0$. As a result, the main loop is interrupted. 
The execution of Algorithm \ref{AlgoFull} continues with Line \ref{AlgoLineFINIT}, which provides no function ($m_w=0$). When $m_w=0$, the tensors $\mu$ and $\nu$ have the single component $\mu^0_0=\nu^0_0=1$ (see Eqs. (\ref{EquationTensorMInertial2}) and (\ref{EquationTensorNInertial})). As a result, $\widehat{g}^0=g^0$   
(see Eq. (\ref{Equationgalpha})). Finally, when $m_w=0$, the codistribution $\tobs$ vanishes (see Eq. (\ref{EquationTOBSDef})).

We obtain (Line \ref{AlgoLineFINALSTEP}):
$\OBS=\left<f^1,f^2,f^3, g^0~\left|~\Omega \right.\right>$.
This is obtained by running the algorithm in Equation (\ref{EquationAlgorithmsMinimalCod}) for the specific case.

At the first iterative step we obtain $\textnormal{span}\left\{\nabla h, ~\nabla\mathcal{L}_{g^0}h\right\}$ and its dimension is $2>1$. At the second iterative step we obtain $\textnormal{span}\left\{\nabla h, ~\nabla\mathcal{L}_{g^0}h, ~\nabla\mathcal{L}_{g^0}\mathcal{L}_{g^0}h, ~\nabla\mathcal{L}_{f^2}\mathcal{L}_{g^0}h \right\}$ and its dimension is $4>2$. At the next step it remains the same, meaning that the algorithm has converged. Therefore:

\[
\OBS=\textnormal{span}\left\{\nabla h, ~\nabla\mathcal{L}_{g^0}h, ~\nabla\mathcal{L}_{g^0}\mathcal{L}_{g^0}h, ~\nabla\mathcal{L}_{f^2}\mathcal{L}_{g^0}h \right\}
\]

We compute the orthogonal distribution.
By a direct computation, we obtain: 

\[
(\OBS)^{\bot}=\textnormal{span}\{
[0~1~0~1~1]^T\}.
\]

The generator of $(\OBS)^{\bot}$ expresses the following invariance ($\epsilon$ is an infinitesimal parameter):

\begin{equation}\label{EquationIntroductionTransfomationIndSetinfRotforGeneral}
\left[\begin{array}{c}
r \\
\phi \\
 v \\
\alpha\\
\theta\\
\end{array}\right]
 \rightarrow
 \left[\begin{array}{c}
r' \\
\phi' \\
 v' \\
\alpha'\\
\theta'\\
\end{array}\right]=\left[\begin{array}{c}
r \\
\phi \\
 v \\
\alpha\\
\theta\\
\end{array}\right]+\epsilon
\left[\begin{array}{c}
0 \\
1 \\
 0 \\
1 \\
1 \\
\end{array}\right],
\end{equation}
which is an infinitesimal rotation about the vertical axis.

In conclusion, for this system (Variant 1), only the first and the third component of the state are observable (i.e., $r$ and $v$), while the three remaining components, $\phi$, $\alpha$, and $\theta$, are not. However, the difference between any two of these angles (e.g., $\phi-\theta$, or $\alpha-\theta$) is observable.

\subsection{State observability for Variant 2}\label{SubSectionObservabilitySystem2}

When $A_x$, $A_y$ are unknown and $\omega$ is known,
by comparing (\ref{EquationApplicationVISFM2Calibrated}) and (\ref{EquationApplicationVISFM2CalibratedOutput}) with (\ref{EquationSystemDefinitionUIO}) we have: $n=5$, $m_u=1$, $m_w=2$,
$p=1$,
$u_1=\omega$, $w_1=A_x$, $w_2=A_y$,

\[
g^0=\left[\begin{array}{c}
v \cos(\alpha-\phi)\\
\frac{v}{r} \sin(\alpha-\phi)\\
 0 \\
 0  \\
 0  \\
\end{array}
\right],~
f^1=\left[\begin{array}{c}
 0 \\
 0  \\
 0  \\
 0  \\
 1  \\
\end{array}
\right],
\]
\[
g^1=\left[\begin{array}{c}
 0 \\
 0  \\
 \cos(\alpha-\theta)   \\
-\frac{1}{v} \sin(\alpha-\theta)\\
 0  \\
\end{array}
\right],~
g^2=\left[\begin{array}{c}
 0 \\
 0  \\
 \sin(\alpha-\theta)   \\
\frac{1}{v} \cos(\alpha-\theta)\\
 0  \\
\end{array}
\right]
\]
As $p=1<2=m_w$ the system is certainly not in canonical form with respect to $w_1$ and $w_2$. We even do not know if it is canonic with respect to them.
We apply Algorithm \ref{AlgoFull}.

\noindent We obtain, at Line \ref{AlgoLineOmegaINIT}, $\Omega=\textnormal{span}\left\{\nabla h\right\}=\textnormal{span}\left\{[0~1~0~0~-1]\right\}$.
Let us compute $\uideg(\Omega)$:
we have $\mathcal{L}_{g^1}h=\mathcal{L}_{g^2}h=0$. As a result, $\textnormal{rank}\left(\RM\left(h\right)
\right)=0$ and $\uideg(\Omega)=0$. 
We have $m=0$, $\mu$ and $\nu$ have the single component $\mu^0_0=~\nu^0_0=1$, $\widehat{g}^0=g^0$, $\giat^1=g^1$, and $\giat^2=g^2$ (see Eqs (\ref{EquationTensorMSynchrom})-(\ref{Equationgalpham})). 
Let us run Algorithm \ref{Algo***}. For $k=1$ we obtain:

\[
\Omega=
\left<\left.f^1, \left\{ 
\begin{array}{ll}
 g^0 & \\
 g^1, & g^2\\
\end{array}
\right\}~\right|\Omega\right>=\textnormal{span}\left\{\nabla h, \nabla\Li_{g^0}h
\right\}
\]
Hence, Finish=False, $\Sigma_*=\Sigma$, and the iteration of the main loop of Algorithm \ref{AlgoFull} returns the same $\Omega$ (this is always the case when $m=0$).
We have
$\uideg(\Omega)=1$.  We obtain:
$m=1$, $\widetilde{h}_1=\Li_{g^0}h=\frac{v}{r} \sin(\alpha-\phi)$. 
The tensor $\mu$ has 4 components:

\[
\mu=\left[\begin{array}{cc}
1&0\\
-\frac{v^2 \sin(2 \alpha - 2 \phi)}{ r^2}& -\frac{\sin(\phi - \theta)}{ r}\\
\end{array}
\right],
\]
where we set $\mu_1^1=\Li_{g^1}\widetilde{h}_1$ (we could also choose $\mu_1^1=\Li_{g^2}\widetilde{h}_1$, as both do not vanish). Its inverse is:

\[
\nu=\left[\begin{array}{cc}
1&0\\
-\frac{v^2 \sin(2 \alpha - 2 \phi)}{ r \sin(\phi - \theta)}& 
 -\frac{ r}{\sin(\phi - \theta)}\\
\end{array}
\right].
\]
Finally:

\[
\widehat{g}^0=\left[
\begin{array}{c}
v \cos(\alpha - \phi)\\
v \sin(\alpha - \phi)/ r\\
 -v^2 \sin(2 \alpha - 2 \phi) \cos(\alpha - \theta)/( r \sin(\phi - \theta))\\
 v \sin(2 \alpha - 2 \phi) \sin(\alpha - \theta)/( r \sin(\phi - \theta))\\
 0\\
\end{array}
\right],~~
\widehat{g}^1=\left[
\begin{array}{c}
 0\\
 0\\
 - r \cos(\alpha - \theta)/\sin(\phi - \theta)\\
  r \sin(\alpha - \theta)/(v \sin(\phi - \theta))\\
 0\\
\end{array}
\right],
\]
\[
\giat^2=\left[
\begin{array}{c}
 0\\
 0\\
 \cos(\alpha - \phi)/\sin(\phi - \theta)\\
 -\sin(\alpha - \phi)/(v \sin(\phi - \theta))\\
 0\\
\end{array}
\right].
\]

We run Algorithm \ref{Algo***}. 
%

It provides Finish=False. In addition, $\Sigma_*$ is the system where the state is augmented, by including the second UI ($A_y$).

%
We back to the main loop and we obtain 
\begin{equation}\label{EquationApplicationVariant2Omega}
\Omega=\textnormal{span}\left\{
\phi - \theta,~~
\frac{v}{r} \sin(\alpha-\phi),~~
A_y\frac{\cos(\phi-\theta)}{ r},~~
A_y\frac{\sin(\phi-\theta)}{ r}
\right\},
\end{equation}

We have now $\uideg(\Omega_1)=2$. As a result, the system has been set in canonical form. The execution continues with Line \ref{AlgoLineFINIT}, which provides:

\[
\widetilde{h}_1=\frac{v}{r} \sin(\alpha-\phi),~~
\widetilde{h}_2=A_y\frac{\cos(\phi-\theta)}{ r}
\]

In addition:

\[
\mu=\left[\begin{array}{ccc}
1&0&0\\
\frac{A_y r\cos(\phi - \theta)-v^2 \sin(2 \alpha - 2 \phi)}{ r^2}& -\frac{\sin(\phi - \theta)}{ r}&0\\
-\frac{A_yv\cos(\alpha - 2\phi + \theta)}{ r^2}&0&\frac{\cos(\phi - \theta)}{ r}
\end{array}
\right],
\]

\[
\nu=\left[\begin{array}{ccc}
1&0&0\\
\frac{-v^2 \sin(2 \alpha - 2 \phi)+A_y r\cos(\phi - \theta)}{ r \sin(\phi - \theta)}& 
 -\frac{ r}{\sin(\phi - \theta)}&0\\
 \frac{A_yv\cos(\alpha - 2\phi + \theta)}{ r\cos(\phi - \theta)}&0&\frac{ r}{\cos(\phi - \theta)}
\end{array}
\right],
\]

\[
\widehat{g}^0=\left[
\begin{array}{c}
v \cos(\alpha - \phi)\\
\frac{v\sin(\alpha - \phi)}{ r}\\
 A_y\sin(\alpha - \theta)+\cos(\alpha - \theta)\frac{A_y r\cos(\phi - \theta)-v^2 \sin(2 \alpha - 2 \phi)}{ r \sin(\phi - \theta)}\\
 \frac{A_y\cos(\alpha - \theta)}{v}+\sin(\alpha - \theta)\frac{v^2\sin(2\alpha - 2\phi) - A_y r\cos(\phi - \theta)}{ r v\sin(\phi - \theta)}\\
 0\\
 -A_yv\frac{\sin(\alpha - 2\phi + \theta)}{ r\sin(\phi - \theta)}\\
\end{array}
\right],
\]

\[
\widehat{g}^1=\left[
\begin{array}{c}
 0\\
 0\\
 - r \cos(\alpha - \theta)/\sin(\phi - \theta)\\
  r \sin(\alpha - \theta)/(v \sin(\phi - \theta))\\
 0\\
 0\\
\end{array}
\right],~~
\widehat{g}^2=\left[
\begin{array}{c}
 0\\
 0\\
 0\\
 0\\
 0\\
  r/\sin(\phi - \theta)\\
\end{array}
\right],
\]
and $\tobs$. In particular,
$\tobs$ is obtained by using (\ref{EquationTOBSDef}) and its computation requires
the integers $s$ and $r$, which are obtained by running Algorithm \ref{AlgoOmegagNonAbel} and Algorithm \ref{AlgoDeltaNonAbel}, respectively. We obtain: $s=r=2$. For this specific case, we obtain that $\tobs$ coincides with $\Omega$ in (\ref{EquationApplicationVariant2Omega}).
The last step (Line \ref{AlgoLineFINALSTEP}) is the computation of:

\[
\OBS=\left<f^1,~\widehat{g}^0,~\widehat{g}^1,~\widehat{g}^2,~\left|~\Omega \right.\right>,
\]
where $\Omega$ is the codistribution given in (\ref{EquationApplicationVariant2Omega}).
This is obtained by running the algorithm in Equation (\ref{EquationAlgorithmsMinimalCod}), for the specific case. The convergence is achieved at the first step.
Therefore:

\begin{equation}\label{EquationApplicationVariant2OBS}
\OBS=\Omega=\textnormal{span}\left\{\nabla
\left[
\phi - \theta\right],~
\nabla\left[\frac{v}{r} \sin(\alpha-\phi)\right],~
\nabla\left[ A_y\frac{\cos(\phi-\theta)}{ r}\right],~
\nabla\left[ A_y\frac{\sin(\phi-\theta)}{ r}\right]
\right\}.
\end{equation}

\noindent We compute the orthogonal distribution.
By a direct computation, we obtain: 

\[
(\OBS)^{\bot}=\textnormal{span}\{
[0~1~0~1~1 ~0]^T,~~
[r~0~v~0~0~A_y]^T
\}.
\]
The first generator of $(\OBS)^{\bot}$ expresses the same invariance found 
in the case discussed in section \ref{SubSectionObservabilitySystem1}. This is the rotation around the vertical axis (see equation (\ref{EquationIntroductionTransfomationIndSetinfRotforGeneral})). This fact is not surprising. By removing the accelerometer we lose information and, consequently, all the degrees of unobservability remain.
 
The second generator expresses the following new invariance
($\epsilon$ is an infinitesimal parameter):

\begin{equation}\label{EquationIntroductionTransfomationScale}
\left[\begin{array}{c}
r \\
\phi \\
 v \\
\alpha\\
\theta\\
A_y\\
\end{array}\right]
 \rightarrow
 \left[\begin{array}{c}
r' \\
\phi' \\
 v' \\
\alpha'\\
\theta'\\
A'_y\\
\end{array}\right]=\left[\begin{array}{c}
r \\
\phi \\
 v \\
\alpha\\
\theta\\
A_y\\
\end{array}\right]+\epsilon
\left[\begin{array}{c}
r \\
0 \\
 v \\
0 \\
0 \\
A_y\\
\end{array}\right],
\end{equation}
which is an infinitesimal scale transform.
The presence of this new degree of unobservability is also not surprising. All the information is provided by the measurements that now only consist of angular measurements (the sensor \VS~ only provides the angle $\beta$ in Fig. \ref{FigVISFM2DCalibrated} and the sensor \IMU~only provides the angular speed). As a result, the system does not have any source of metric information.

In conclusion, for this system (Variant 2), no component of the state is observable.
However, as for Variant 1, the difference between any two of the three angles $\phi$, $\alpha$, and $\theta$,  (e.g., $\phi-\theta$, or $\alpha-\theta$) is observable. In addition, regarding the metric quantities in the new state (i.e., $r$, $v$, and $A_y$), the ratio between any two of them (e.g., $\frac{v}{r}$, or $\frac{A_y}{r}$) is observable.

Note that, even if the above computation could seem laborious, it is automatic. In particular, it is carried out by a simple code that uses symbolic computation.

\subsection{State observability for Variant 3}\label{SubSectionObservabilitySystem3}

By running Algorithm \ref{AlgoFull}, we obtain that the observability codistribution coincides with the one obtained in Section \ref{SubSectionObservabilitySystem1}, i.e., in the case of a complete \IMU~ (Variant 1). In other words, in this case, the system has a single degree of unobservability that is the system invariance against a rotation around the vertical axis.

\chapter{Conclusion}\label{ChapterConclusion}

This paper provided the general analytical solution of a fundamental open problem introduced long time ago (in the middle of the 1960's). The problem is the possibility of introducing an analytical and automatic test able to check the observability of the state that characterizes a system whose dynamics are also driven by inputs that are unknown. This is the extension of the well known {\it Observability Rank Condition} introduced in the 1970's and that does not account for the presence of unknown inputs. This paper provided this extension. This problem arises in a large class of domains, ranging from mechanical engineering, robotics, computer vision, up to biology, chemistry and economics.
Observability is a fundamental structural property of any dynamic system
and describes the possibility of reconstructing the state that characterizes the system from observing its inputs and outputs.
The dynamics of most real systems are driven by inputs that are usually unknown.
Very surprisingly, the complexity of the solution here introduced is comparable to the complexity of the observability rank condition. 
Given any nonlinear system characterized by any type of nonlinearity, driven by both known and unknown inputs, the state observability is obtained automatically, i.e., by the usage of a very simple code that uses symbolic computation. This is a fundamental practical (and unexpected) advantage.

To obtain this solution, the paper used several important new concepts introduced very recently in \cite{SIAMbook}.
Note that also in \cite{SIAMbook} a solution of the unknown input observability problem was introduced.
The novelties of the complete solution introduced by this paper, with respect to the solution given in \cite{SIAMbook} are:

\begin{enumerate}

\item Full characterization of the concept of {\it canonicity with respect to the unknown inputs}, given in Chapter \ref{ChapterSystem} and in the first part of Chapter \ref{ChapterSolutionNonCanonic}.

\item Algorithm \ref{AlgoFull}, which is the general solution that holds even in the non canonic case and not even canonizable. In particular, when the system is not canonizable, Algorithm \ref{AlgoFull} returns a new system with the highest unknown input degree of reconstructability, together with the observability codistribution.


\item A new criterion of convergence of the solution in the canonic case. 
In particular, 
the criterion proposed in \cite{SIAMbook}, which is based on the computation of the tensor $\mathcal{T}$, can fail. The new criterion here introduced, which extended the one introduced in \cite{SARAFRAZI} to the general case with drift, multiple unknown inputs and TV, holds always (and is even simpler). In addition, 
the algorithm that solves the problem was written in a new manner (see Algorithm \ref{AlgoNonAbel}), where the initialization step includes all the terms of Algorithm \ref{AlgoBook} that make the convergence criterion of Algorithm \ref{AlgoBook} non trivial.

\end{enumerate}

Finally, as a simple consequence of the results here obtained, the paper provided a preliminary answer to the problem of unknown input reconstruction, which is intimately related to the problem of state observability. A final answer to this fundamental problem can be found in  \cite{arXivODE,arXivTAC,arXivErratum}.

The solution was illustrated with a simple application.
We studied the observability properties of a nonlinear system in the framework of visual-inertial sensor fusion. The dynamics of this system are driven by two unknown inputs and one known input and they are also characterized by a nonlinear drift. The system is not in canonical form with respect to its unknown inputs. However, by following the steps of Algorithm \ref{AlgoFull}, it was set in canonical form, and all the observability properties were automatically obtained.

\end{document}